\documentclass[brochure, french]{smfbourbaki}
\pdfoutput=1        

\usepackage[T1]{fontenc}
\usepackage{lmodern,amssymb,bm, mathrsfs}
\usepackage{graphicx}

\usepackage{babel}

\usepackage[utf8]{inputenc}

\usepackage[colorlinks=true, linkcolor=blue, citecolor=red, urlcolor=blue]{hyperref}

\usepackage[stretch=10,shrink=10,step=2,kerning=true,protrusion=true,expansion=true,final]{microtype}

\usepackage[
backend=bibtex,
style=authoryear, 
citestyle=authoryear-comp,
maxnames=7,
sortcites=false 
]{biblatex}
\usepackage{csquotes}

\DefineBibliographyExtras{french}{\restorecommand\mkbibnamefamily}

\DeclareDelimFormat{nameyeardelim}{\addcomma\space}

\DeclareNameAlias{sortname}{given-family}

\renewcommand{\bibnamedash}{\leavevmode\raise3pt\hbox to3em{\hrulefill}\space}

\AtEveryBibitem{%
  \clearfield{issn} 
  \clearfield{isbn} 
  \clearfield{doi} 
  \clearlist{language} 
  \ifentrytype{online}{}{
  \ifentrytype{unpublished}{}{
    \clearfield{url}
  }
  }
}

\renewbibmacro{in:}{%
    \ifentrytype{article}{}{\printtext{\bibstring{in}\intitlepunct}}}

\DeclareFieldFormat[article,periodical,inreference]{number}{\mkbibparens{#1}}
\DeclareFieldFormat[article,periodical,inreference]{volume}{\mkbibbold{#1}}
\renewbibmacro*{volume+number+eid}{%
    \printfield{volume}%
    \setunit*{\addthinspace}
    \printfield{number}%
    \setunit{\addcomma\space}%
    \printfield{eid}}

\DeclareFieldFormat[article,inbook,incollection]{title}{\enquote{#1}\addcomma} 

\usepackage{enumitem}
\setlist[enumerate]{label={\roman*}}
\newlist{proprietes}{enumerate}{1}
\setlist[proprietes]{label={(\roman*)}}
\newlist{conditions}{enumerate}{1}
\setlist[conditions]{label={(\roman*)}, font={\normalfont}}
\newlist{assertions}{enumerate}{1}
\setlist[assertions]{label={\alph*)}, font={\normalfont}}

\usepackage{xspace}

\renewcommand{\iemes}{\nobreakdash-i\`emes\xspace}
\newcommand{\numero}{n\textsuperscript o\kern+.2em}
\newcommand{\numeros}{n\textsuperscript {os}\kern+.2em}

\newcommand{\cf}{\emph{cf.}\nobreak\xspace}
\newcommand{\resp}{\textup{resp.}\nobreak\xspace}

\newcommand{\loccit}{\emph{loc. cit.}\nobreak\xspace}

\newcommand{\pv}{\textup{;}}
\newcommand{\pt}{\textup{.}}
\newcommand{\ass}[1]{\emph{{#1}})}

  \DeclareMathSymbol{A}{\mathalpha}{operators}{`A}%
  \DeclareMathSymbol{B}{\mathalpha}{operators}{`B}%
  \DeclareMathSymbol{C}{\mathalpha}{operators}{`C}%
  \DeclareMathSymbol{D}{\mathalpha}{operators}{`D}%
  \DeclareMathSymbol{E}{\mathalpha}{operators}{`E}%
  \DeclareMathSymbol{F}{\mathalpha}{operators}{`F}%
  \DeclareMathSymbol{G}{\mathalpha}{operators}{`G}%
  \DeclareMathSymbol{H}{\mathalpha}{operators}{`H}%
  \DeclareMathSymbol{I}{\mathalpha}{operators}{`I}%
  \DeclareMathSymbol{J}{\mathalpha}{operators}{`J}%
  \DeclareMathSymbol{K}{\mathalpha}{operators}{`K}%
  \DeclareMathSymbol{L}{\mathalpha}{operators}{`L}%
  \DeclareMathSymbol{M}{\mathalpha}{operators}{`M}%
  \DeclareMathSymbol{N}{\mathalpha}{operators}{`N}%
  \DeclareMathSymbol{O}{\mathalpha}{operators}{`O}%
  \DeclareMathSymbol{P}{\mathalpha}{operators}{`P}%
  \DeclareMathSymbol{Q}{\mathalpha}{operators}{`Q}%
  \DeclareMathSymbol{R}{\mathalpha}{operators}{`R}%
  \DeclareMathSymbol{S}{\mathalpha}{operators}{`S}%
  \DeclareMathSymbol{T}{\mathalpha}{operators}{`T}%
  \DeclareMathSymbol{U}{\mathalpha}{operators}{`U}%
  \DeclareMathSymbol{V}{\mathalpha}{operators}{`V}%
  \DeclareMathSymbol{W}{\mathalpha}{operators}{`W}%
  \DeclareMathSymbol{X}{\mathalpha}{operators}{`X}%
  \DeclareMathSymbol{Y}{\mathalpha}{operators}{`Y}%
  \DeclareMathSymbol{Z}{\mathalpha}{operators}{`Z}%

\let\epsilon\varepsilon
\let\phi\varphi
\let\leq\leqslant
\let\geq\geqslant

\let\phi\varphi
\let\eps\varepsilon

\let\tilde\widetilde

\def\GL{\mathop{\mathbf{GL}}}
\def\SL{\mathop{\mathbf{SL}}}
\def\diag{\operatorname{diag}}

\def\Tr{\operatorname{Tr}}

\def\dual#1#2{\langle#1,#2\rangle}

\def\norm#1{\mathopen\|#1\mathclose\|}

\def\abs#1{\mathopen|{#1}\mathclose|}

\def\Card{\operatorname{Card}}

\def\goodwedge{\bm{\wedge}}

\def\im{\operatorname{Im}}

\def\crochets#1{\mathopen{[}#1\mathclose{]}}

\renewcommand{\setminus}{\mathchoice
  {\mathbin{\vrule height .72ex width 1.61ex depth -.38ex}}
  {\mathbin{\vrule height .72ex width 1.61ex depth -.38ex}}
  {\mathbin{\vrule height .50ex width 0.85ex depth -.28ex}}
  {\mathbin{\vrule height .20ex width 0.570ex depth -.24ex}}
}

\newenvironment{smallrema}{\list{}%
   {\advance\leftmargin\parindent
    \advance\labelwidth\labelsep
    \itemindent\parindent
	\small }\item\relax}{\endlist \ignorespacesafterend}

\newcommand*\separateur{\begin{center} \includegraphics[width=0.12\textwidth]{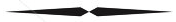}\end{center} }

\def\B{{\mathbf B}}
\def\C{{\mathbf C}}
\def\R{{\mathbf R}}
\def\N{{\mathbf N}}
\def\Z{{\mathbf Z}}
\def\Q{{\mathbf Q}}
\def\M{{\mathbf M}}
\def\P{{\mathbf P}}

\def\J{{\mathbf J}}
\def\L{{\mathbf L}}
\def\T{{\mathbf T}}
\def\G{{\mathbf G}}
\def\F{{\mathbf F}}

\def\bS{{\mathbf S}}

\newcommand{\lie}[1]{\mathfrak{#1}}
\def\Lie{\mathrm{Lie}}

\def\Greg{{G}_{\mathrm{reg}}}
\def\Gtopss{G_{\mathrm{top.ss}}}
\def\Gtopnil{G_{\mathrm{top.nil}}}
\def\Gder{\mathbf{G}_{\mathrm{d}}}

\def\ad{\mathrm{ad}}

\def\Ind{\mathrm{ind}}
\def\cInd{\textrm{c-ind}}
\def\Irr{\mathrm{Irr}}
\def\IrrCusp{\mathrm{Irr}_{\textrm{cusp}}}

\def\Gfini{\mathsf{G}}
\def\Sfini{\mathsf{S}}

\def\Spm{S_{\pm}}
\def\Sreg{S_{\mathrm{reg}}}

\def\EspHeis{V^{\mathrm{Heis}}}
\def\DeLus{\mathsf{R}}

\def\rhalf{\frac{r}{2}}

\def\premier{\mathfrak{p}}
\def\entiers{\mathfrak{o}}
\def\kprem{k_{F,p}}
\def\ClotureF{\overline{F}}
\def\Fsep{F^{\mathrm{s}}}
\def\Fpmalpha{F_{\pm\alpha}}
\def\Falpha{F_{\alpha}}
\def\Gal{\mathrm{Gal}}

\def\dualG{\widehat{G}}
\def\Paquet{\Pi}
\def\Langlands{\mathrm{Lan}}
\def\Weil{\mathbf{W}}
\def\WeilDeligne{\mathbf{WD}}
\def\Inertie{\mathbf{I}}
\def\InertieSauvage{\mathbf{P}}
\def\Plongements{\mathscr{J}}

\newcommand{\ch}[1]{\negthinspace\negthinspace\negthinspace\phantom{a}^\vee #1}
\newcommand{\lgr}[1]{\negthinspace\negthinspace\negthinspace\phantom{a}^L #1}

\def\Racines{\mathsf{\Phi}}
\def\RacinesSym{\Racines_{\mathrm{sym}}}

\def\Imm{\mathscr{B}}
\def\Immred{\mathscr{B}_{\textrm{red}}}
\def\reduit#1{\mathopen[#1\mathclose]}

\def\Char{\operatorname{\mathsf{X}}}
\def\Charspec{\operatorname{\mathsf{X}_{\mathrm{spec}}}}
\def\Charspecreg{\operatorname{\mathsf{X}_{\mathrm{spec}, \mathrm{reg}}}}

\date{Novembre 2023}
\bbkannee{76\textsuperscript{e} année, 2023--2024}  
\bbknumero{1211}                                      

\title{Progrès récents\\ sur les représentations supercuspidales}

\author{Alexandre Afgoustidis}
\address{CNRS \& Institut \'Elie Cartan de Lorraine}
\email{alexandre.afgoustidis@math.cnrs.fr}

\begin{document}
\maketitle



\section{Introduction}


\subsubsection*{\textbf{\textup{1.1}}}
Soient~$F$ un corps local et~$\G$ un groupe algébrique réductif connexe défini sur~$F$.
C'est aux représentations complexes du groupe~$G=\G(F)$ que l'on s'intéresse ici. 

Lorsque $F$ est archimédien, donc isomorphe à $\R$ ou~$\C$, le groupe~$G$ est un groupe de Lie avec un nombre fini de composantes connexes.
L'étude de ses représentations est le grand~œuvre d'Harish-Chandra, commencé immédiatement après 1945.
Pour les questions qui nous occupent ci-dessous, l'essentiel était compris à la fin des années 1970 ; 
la classification  de~\textcite{Langlands_CIRRAG} des représentations irréductibles «\,{}admissibles\,{}» a~d'ailleurs tout juste cinquante ans. 

Lorsque~$F$ n'est pas archimédien, le groupe~$G$ est totalement discontinu.
L'étude de ses représentations «\,{}lisses\,{}» a débuté vers 1960~\parencite{Mautner_1958, Bruhat}.
Elle a pris son essor une dizaine d'années plus tard, portée par son lien avec les formes automorphes,
et se poursuit depuis sans relâche --- avec, à l'horizon, les conjectures de Langlands et leur promesse arithmétique.


\subsubsection*{\textbf{\textup{1.2}}}\label{sec:reduction_jacquet}
L'analogie avec les groupes réels a joué un rôle central pour établir les fondements de la théorie non archimédienne, et pour la conception des conjectures de~Langlands. 
Les bases de cette théorie s'expriment dans des termes très proches du cas réel, plus proches qu'on ne le devinerait d'après la structure seule des groupes~---
au point que cette analogie a pu paraître miraculeuse.
Harish-Chandra empruntait à la géométrie algébrique l'expression « principe de Lefschetz » pour parler de ces rapprochements. 

Mais les chemins ont divergé dès le début des années 1970. Dans le cas des groupes~réels, la description des représentations irréductibles se ramenait, par inductions successives,
à la recherche de représentations \emph{de carré intégrable}  de certains sous-groupes réductifs de~$G$.
Vers~\cite*{HarishChandra_SerieDiscrete}, Harish-Chandra avait réussi à classifier ces représentations,
en les paramétrant par des caractères de sous-groupes de Cartan compacts et en donnant une formule pour leurs caractères globaux
(voir le \S\,{}\ref{sec:cas_reel} et l'exposé de~\cite{Duflo_1977_seminaire}).
Ce paramétrage par des caractères de tores était un ingrédient essentiel de la classification de Langlands.
La situation s'est révélée très différente dans le cas non archimédien : la classification des représentations de carré intégrable reste hors d'atteinte à ce jour (sauf pour certains groupes classiques),
et les informations explicites dont on dispose sur les caractères restent minces. 

Cependant, à la fin des années 1950, \textcite{Mautner_1964} avait découvert
que si~$F$ n'est pas archimédien, il existe des représentations de $\mathbf{PGL}(2,F)$ dont les coefficients matriciels sont \emph{à support compact}.
Or, pour les représentations de groupes réductifs réels, aucun coefficient matriciel ne peut avoir cette propriété, sauf dans des cas triviaux.
Il est vite apparu que les représentations \emph{supercuspidales}, celles dont les coefficients matriciels sont à support compact modulo le centre,
joueraient un rôle crucial dans l'organisation de la théorie non archimédienne : toute représentation lisse irréductible de~$G$ peut se plonger dans l'induite parabolique d'une représentation supercuspidale de sous-groupe de Levi de~$G$ \parencite{Jacquet_1971, Godement_Jacquet}.
On pouvait donc adopter la stratégie suivante pour étudier les représentations de~$G$ dans le cas non archimédien :
\begin{conditions}
\item Décrire le dual supercuspidal de~$G$, c'est-à-dire l'ensemble des classes d'équivalence de représentations supercuspidales irréductibles ; 
\item En supposant connu le dual supercuspidal pour les sous-groupes de Levi de~$G$, décrire plus précisément les autres représentations de~$G$. 
\end{conditions}
Les progrès sur ces deux problèmes ont été continus depuis un demi-siècle.
C'est le premier des deux qui fait l'objet de cet exposé.
Je n'aborderai pas le second problème, bien qu'il ait vu lui aussi d'importants progrès ces dernières années : pour quelques indications bibliographiques, voir le \S\,{}\ref{sec:beyond}.
Pour le cas de $\GL(n)$, où la réduction au cas supercuspidal fonctionne sans obstacle, on pourra consulter l'exposé de \textcite{Rodier_1981_seminaire}.


\subsubsection*{\textbf{\textup{1.3}}}
Parmi les grands courants qui irriguent encore l'étude des représentations supercuspidales, certains prennent leur source dans la période  mentionnée ci-dessus.
Par exemple, il a semblé dès cette époque que le problème pouvait être plus abordable dans des cas~\mbox{«\,{}modérés\,{}»}
où la caractéristique résiduelle~$p$ de~$F$ est suffisamment grande relativement à~$\G$ (la~bonne condition à imposer à $p$ fait  partie de la question).
Si~$\G$ est le groupe $\GL(n)$ et si~$p$ ne~divise pas~$n$, Howe a décrit en~\cite*{Howe_1977}
une construction de représentations supercuspidales à partir de caractères de tores maximaux elliptiques de~$\G$,
qui s'accorde très bien avec le cas des groupes réels.
Moy a montré en~\cite*{Moy_1986} que cette construction donne tout le dual supercuspidal dans ce cas,
et en a déduit une première proposition de correspondance de Langlands  «\,{}modérée\,{}» pour~$\GL(n)$.

Cet exposé, qui est entièrement consacré au cas modéré, décrit en quelque sorte l'actualité de ce courant d'idées lorsque~$\G$ est un groupe réductif quelconque.


\subsubsection*{\textbf{\textup{1.4}}}\label{sec:BK}
Un autre courant, toujours porteur, mène vers des méthodes qui permettent d'aller au-delà du cas modéré --- au prix d'une plus grande sophistication combinatoire.
Ces méthodes portent généralement sur les groupes classiques et leurs formes intérieures.
Pour la détermination du dual cuspidal de $\GL(n)$, c'est ce qu'ont fait~\textcite{Carayol_1984} pour~$n$ premier, puis~\textcite{Bushnell_Kutzko} pour~$n$ quelconque ;
voir l'exposé d'\textcite{Henniart_1991_seminaire}. 
Les idées de Bushnell et Kutzko ont été généralisées avec beaucoup de succès pour les groupes classiques, en utilisant le cas de~$\GL(n)$ comme camp de base  :
elles permettent de construire toutes les représentations supercuspidales des groupes classiques en caractéristique résiduelle $\neq 2$, et celles des formes intérieures de $\GL(n)$.
Voir le \S\,{}\ref{sec:beyond} pour quelques indications bibliographiques.


\subsubsection*{\textbf{\textup{1.5}}}\label{sec:cas_bien_modere}
Parallèlement, l'étude du cas modéré s'est poursuivie, avec en vue la terre promise des groupes réductifs généraux.
Grâce à des idées introduites par Moy et Prasad en~\cite*{Moy_Prasad}, \textcite{Adler_Pacific} puis \textcite{Yu_2001} ont décrit une construction très générale de représentations supercuspidales, que j'expose dans le \S\,{}\ref{sec:Yu_et_Fintzen}. La filtration de~Moy--Prasad et la construction de Yu ont permis d'avancer dans des voies jusqu'alors presque inaccessibles, notamment vers des renseignements explicites sur les formules de caractères  \parencite{Adler_Spice, DeBacker_Spice, Spice_Compositio, Spice_2}.

 Les progrès autour des représentations obtenues par la construction de Yu se sont beaucoup accélérés ces dernières années.
 Les avancées les plus spectaculaires concernent le cas où~$\G$ se déploie sur une extension modérément ramifiée de~$F$ et où la caractéristique résiduelle~$p$ ne divise pas l'ordre~$\abs{W}$ du groupe de Weyl de~$\G$.
Dans ce cas, que j'appellerai «\,{}bien modéré\,{}» dans ce texte, le tableau général des représentations supercuspidales est maintenant complet.
Comme nous le verrons, on dispose aussi d'informations nouvelles qui rendent une partie du tableau beaucoup plus accessible. 


\subsubsection*{\textbf{\textup{1.6}}}\label{sec:Kim_Fintzen}
L'avancée la plus simple à énoncer est que dans le cas «\,{}bien modéré\,{}», la construction de Yu est exhaustive (\cite{Fintzen_Annals}).
Grâce à des travaux de~\textcite{Hakim_Murnaghan} complétés par~\textcite{Kaletha_regular_supercuspidals},
cela induit une véritable \emph{classification} des représentations supercuspidales. 

Ju-Lee Kim avait déjà montré, en~\cite*{Kim_exhaustion}, que la construction de Yu fournit toutes les représentations supercuspidales irréductibles
lorsque~$F$ est de caractéristique nulle et~$p$ est «\,{}très grand\footnote{Il serait en principe possible d'expliciter une condition suffisante sur~$p$ pour le résultat de Kim, mais je n'ai connaissance d'aucun exemple où cela ait été fait. }\,{}».
Le résultat de Fintzen étend celui de Kim ;
pour mesurer la portée de l'amélioration, mentionnons seulement que si~$\G$ est semi-simple et~$p$ ne divise pas~$\abs{W}$, alors~$\G$ se déploie toujours sur une extension modérément ramifiée ;
et que si le système de racines de~$\G$ est de type~$\mathbf{E_8}$, alors $\abs{W}=696729600=2^{14}\cdot 3^{5} \cdot 5^2 \cdot 7$,
ce qui devrait montrer l'attrait de la nouvelle condition sur~$p$. 


\subsubsection*{\textbf{\textup{1.7}}} \label{sec:idee_regulieres}
Comme on le verra, la construction de Yu se fait à partir de données en apparence sophistiquées. Ces données comprennent, entre autres, des représentations de sous-groupes de~$G$ qui proviennent essentiellement de représentations cuspidales de groupes réductifs finis.
Voilà qui semble bien différent des paramétrages simples, par des caractères de tores, que donnaient \textcite{HarishChandra_SerieDiscrete} pour la série discrète si~\mbox{$F=\R$},
et \textcite{Howe_1977} pour le dual cuspidal de $\GL(n)$ si~$p\nmid n$. 

Dans le cas bien modéré, les travaux récents de \textcite{Kaletha_regular_supercuspidals, Kaletha_Lpackets} éclairent d'un jour nouveau  les représentations de Yu.
Poursuivant une idée de Murnaghan~\parencite*{Murnaghan_2011}, ils  établissent que  «\,{}presque toute\,{}» représentation supercuspidale peut être construite à partir d'un caractère de tore maximal elliptique de~$G$.
Il s'agit de remarquer que la plupart des représentations cuspidales des groupes finis sont construites à partir de caractères de tores (théorie de Deligne--Lusztig),
et de s'astreindre à utiliser ces dernières comme ingrédients «\,{}finis\,{}» dans la construction de Yu.
Les représentations ainsi obtenues, que Kaletha nomme «\,{}régulières\,{}» et «\,{}non~singulières\,{}», forment l'essentiel du dual supercuspidal.
On les imagine volontiers plus maniables que les autres. 

De fait, Kaletha remarque que pour ces représentations, les travaux sur les formules de caractère mentionnés ci-dessus (Adler, DeBacker, Spice...)
peuvent être complétés jusqu'à obtenir une formule du caractère \emph{entièrement explicite} (sur une partie adéquate de~$G$),
qui s'avère avoir la même structure que la formule d'Harish-Chandra pour les caractères des séries discrètes des groupes réels ---
et coïncide même avec cette formule si l'on en interprète correctement les ingrédients.
Cette coïncidence est d'autant plus frappante qu'à ma connaissance, le principe de Lefschetz, qui était l'une des motivations de  la construction de~\textcite{Howe_1977},
avait presque entièrement disparu des travaux évoqués aux \nos \hyperref[sec:BK]{1.4}--\hyperref[sec:Kim_Fintzen]{1.6}. 


\subsubsection*{\textbf{\textup{1.8}}}
L'application la plus spectaculaire des idées du \no \hyperref[sec:idee_regulieres]{1.7}
est probablement la construction d'une correspondance de Langlands, dans le cas «\,{}bien modéré\,{}»,
pour les paramètres de Langlands supercuspidaux.
Sans entrer pour le moment dans les détails, la conjecture de Langlands locale suggère qu'on peut répartir «\,{}naturellement\,{}» les représentations lisses irréductibles de~$G$
en paquets finis indexés par certains homomorphismes continus~$\phi\colon \WeilDeligne_F \to \lgr{G}$,
où~$\WeilDeligne_F$ est un groupe de nature arithmétique attaché à~$F$ et où~$\lgr{G}$ est le dual de Langlands de~$G$.
On peut imaginer que savoir répartir  les représentations supercuspidales en paquets de Langlands soit un pas crucial dans cette voie.
Pour le groupe~$\GL(n)$, la construction de la correspondance de Langlands s'y ramène entièrement :
voir les exposés de \textcite{Rodier_1981_seminaire} et de \textcite{Carayol_2000_seminaire}. 

Pour les groupes réductifs généraux, en revanche, on s'attend à ce que pour certains paramètres~$\phi$,
le paquet~$\Pi(\phi)$ contienne une représentation supercuspidale mais ne soit pas formé entièrement de représentations supercuspidales.
On peut se demander quels sont les paramètres~$\phi$ donnant lieu à un paquet «\,{}entièrement supercuspidal\,{}», et comment construire les paquets attachés à ces paramètres.
C'est à ces questions que \textcite{Kaletha_Lpackets}  répond complètement dans le cas bien modéré --- 
par une méthode purement locale, directement inspirée de la construction de~\textcite{Langlands_CIRRAG}.

Dans le cas des groupes réels, la correspondance construite par  Langlands est tout entière basée sur le lien entre séries discrètes de~$G$ et caractères de tores maximaux :
les tores maximaux dans~$\G$ sont liés à ceux de (la composante neutre de)~$\lgr{G}$,
et Langlands combine ces deux liens pour passer d'un paramètre~$\phi$ au paquet~$\Paquet(\phi)$.
La construction de Kaletha, exposée dans le \S\,{}\ref{sec:LLC}, suit le même chemin,
en utilisant le lien entre représentations supercuspidales non singulières et caractères de tores maximaux~elliptiques.

Cette construction explicite est bien éloignée des méthodes géométriques, et globales,
qui ont mené aux retentissants travaux de~\textcite{Genestier_Lafforgue} puis de~\textcite{Fargues_Scholze}.
Espérons qu'il devienne bientôt possible de les comparer. 


\subsubsection*{\textbf{\textup{1.9}}} \label{sec:intro_cover}
La construction de Langlands, comme la première version de celle de Kaletha, a un désagrément :
la route qui relie les paramètres de Langlands $\phi$ à des caractères de tores maximaux
est, selon le mot de \textcite{Carayol_1984}, tortueuse.
Dans le cas archimédien, Adams et Vogan ont montré en~\cite*{AV1} qu'on pouvait simplifier la construction de Langlands
en y faisant intervenir certains revêtements doubles de tores maximaux de~$G$ ;
cela éclaire aussi les formules de caractère.
Récemment, \textcite{Kaletha_Covers} a montré comment généraliser au cas non archimédien le revêtement d'Adams et Vogan.
Il a esquissé les simplifications considérables que cela apporte
(disparition des facteurs de transfert des formules de caractère, existence de $L$-plongements canoniques...).


\subsubsection*{\textbf{\textup{1.10}}}  
Voici le plan de ce texte. 
Les \S\,{}\ref{sec:generalites} et \ref{sec:cas_reel} sont préparatoires ;
j'y consacre quelques pages au cas réel pour donner un cas simple de l'usage de revêtements doubles,
et pour que l'on puisse mieux voir le rôle du principe de Lefschetz dans les \S\,{}\ref{sec:supercuspidales_nonsingulieres}--\ref{sec:LLC}. 
Ce~n'est qu'au~\mbox{\S\,{}\ref{sec:Yu_et_Fintzen}} que j'entre dans le vif du sujet,
avec une description de la construction de Yu et des résultats de classification récents dans le cas « bien modéré ».
La suite de l'exposé porte sur les représentations régulières ou non singulières :
le \S\,{}\ref{sec:supercuspidales_nonsingulieres} décrit ces représentations et leurs caractères,
et le \S\,{}\ref{sec:LLC} décrit la correspondance de Langlands supercuspidale.
Dans ces deux paragraphes, j'ai choisi d'utiliser dès le départ le revêtement double de Kaletha pour décrire les constructions et les résultats. Cela permet d'assez nombreux raccourcis et clarifications, notamment dans la description de la correspondance de Langlands.
Enfin, le \S\,{}\ref{sec:beyond} donne des indications bibliographiques sur des progrès récents liés à notre sujet, qui vont au-delà du cas modéré ou des  représentations~supercuspidales.


\subsubsection*{\textbf{\textup{1.11}}}
Plusieurs des résultats décrits ici ont été très bien exposés récemment. 
La construction de Yu est magnifiquement expliquée dans les notes de~\textcite{Fintzen_Harvard}, que le \S\,{}\ref{sec:Yu_et_Fintzen} suit d'assez près. 
Le texte de~\textcite{Kaletha_ICM} pour le Congrès international donne un aperçu concis et très peu technique du sujet.
Par ailleurs, les exposés de Fintzen, Kaletha ou Spice que l'on peut voir en ligne m'ont été fort utiles 
(je pense à des exposés donnés en visioconférence en 2020, ou à ceux d'une école d'été de l'IHÉS en 2022).
Enfin, les représentations des groupes réductifs ont fait l'objet de nombreux exposés au séminaire Bourbaki ;
il serait absurde de les citer tous, mais en rédigeant ces notes, j'ai utilisé ceux 
de van~\textcite{VanDijk} et \textcite{Cartier_1977_seminaire} sur les caractères, de~\textcite{Serre_seminaire} sur la théorie de Deligne--Lusztig pour les groupes finis,
de~\textcite{Duflo_1977_seminaire} sur la série discrète des groupes réels, de~\textcite{Rodier_1981_seminaire} sur les travaux de Bernstein et~Zelevinsky,
 d'\textcite{Henniart_1991_seminaire} sur les représentations supercuspidales de $\GL(n)$, et de \textcite{Carayol_2000_seminaire} sur la correspondance de Langlands pour ce groupe. 

\separateur

Les éclaircissements et les conseils de Jeffrey Adams, Anne-Marie Aubert, Benoît Cadorel, Charlotte Chan, Jessica Fintzen, Auguste Hébert, Sergey Lysenko,  Masao Oi, Guy Rousseau et, particulièrement, Tasho Kaletha, m'ont été  précieux pour préparer cet exposé. Les remarques de Guy Henniart, Monica Nevins et Maarten Solleveld ont permis de corriger et d'améliorer ce texte. Merci !


\section{Généralités sur les représentations}\label{sec:generalites}

Pour plus d'informations, ou des démonstrations, on pourra consulter
le livre \mbox{de~Renard~\parencite*{Renard}},
les notes de~\textcite{Casselman}
et l'indémodable \parencite*{Corvallis}.


\subsection{Corps locaux, groupes réductifs et sous-groupes compacts}


\subsubsection{}
Soit~$F$ un corps local. Si~$F$ est archimédien, alors c'est~$\R$ ou~$\C$.
Si~$F$ n'est pas archimédien, notons~$k_F$ le corps résiduel de~$F$, quotient de l'anneau d'entiers~$\entiers_F$ par son unique idéal maximal.
Le corps~$k_F$ est fini et on note~$p$ sa caractéristique.
Si~$F$ est de caractéristique nulle, alors~$F$ est isomorphe  à une extension finie du corps~$\Q_p$ des nombres $p$-adiques.
Si la caractéristique de~$F$ est non nulle, elle est égale à~$p$, et~$F$ est isomorphe au corps des séries de Laurent formelles à coefficients dans~$k_F$.


\subsubsection{} \label{sec:sg_compacts}
Soient~$\G$ un groupe réductif algébrique connexe défini sur~$F$ et~$G$ le groupe de points rationnels~$\G(F)$.
C'est un groupe localement compact.
Les sous-groupes compacts de~$G$ jouent un rôle essentiel dans ce qui suit.
Si~$F$ est archimédien, les sous-groupes compacts maximaux de~$G$ sont tous conjugués ;
de plus, si~$K$ est un tel sous-groupe, la variété différentiable~$G/K$ est difféomorphe à son espace tangent en l'identité.

Si~$F$ n'est pas archimédien, alors il y a un nombre fini de classes de conjugaison de sous-groupes compacts maximaux,
mais ce nombre de classes de conjugaison est~\mbox{$>1$} en~général.
De plus, le groupe~$G$ admet « beaucoup » de sous-groupes compacts :
la famille des sous-groupes compacts ouverts de~$G$ est une base de voisinages de l'identité~$1_G$.
Les relations d'incidence entre  sous-groupes compacts ouverts sont importantes pour le \S\,{}\ref{sec:Yu_et_Fintzen} ;
mentionnons que si~$K_1, K_2$ sont des sous-groupes compacts ouverts de~$G$ avec $K_1 \subset K_2$, alors $K_1$ est d'indice fini dans~$K_2$. 

On s'intéressera souvent aux sous-groupes \emph{compacts ouverts modulo le centre de~$G$} :
ce sont ceux qui sont ouverts, contiennent le centre~$Z(G)$, et dont le quotient par~$Z(G)$ est compact.


\subsection{Représentations supercuspidales et séries discrètes}

Dans ce numéro, on suppose~$F$ non archimédien ;
nous verrons au \S\,{}\ref{sec:gen_arch} comment modifier les énoncés qui suivent dans le cas archimédien.

\subsubsection{Représentations lisses}
Soit~$\pi$ une représentation de~$G$ sur un espace vectoriel complexe~$V$,
c'est-à-dire un homomorphisme de~$G$ dans~$\GL(V)$.
On dit que~$\pi$ est \emph{lisse} si le stabilisateur de tout élément de~$V$ est ouvert dans~$G$
(cela revient à demander que l'application de~$G \times V$ dans~$V$ provenant de~$\pi$ soit continue lorsqu'on munit~$V$ de la topologie discrète) ;
et que~$\pi$ est \emph{irréductible} si~$V$ est non nul et si ses seuls sous-espaces stables par $\pi(G)$  sont~$\{0\}$ et~$V$.
Si~$\pi$ est lisse et irréductible, alors elle est \emph{admissible} :
pour tout sous-groupe compact ouvert~$K$ de~$G$, l'espace $V^K$ des vecteurs $\pi(K)$-invariants est de dimension finie.
On note~$\Irr(G)$  l'ensemble des classes d'équivalence de représentations lisses irréductibles de~$G$, pour la notion usuelle d'équivalence (« dual lisse de~$G$ »).

\subsubsection{Induction  compacte} \label{sec:induction}
Si l'on se donne un sous-groupe fermé~$H$ de~$G$ et une représentation lisse~$\rho$ de~$H$ sur un espace vectoriel~$W$,
 on peut construire une représentation lisse de~$G$, comme suit. 
Soit~$V$ l'espace des fonctions $f\colon G \to W$ vérifiant : 
\begin{conditions}
\item Il existe un sous-groupe compact ouvert~$K$ de~$G$ tel qu'on ait $f(hgk)=\rho(h) f(g)$ pour tout $(h,g,k) \in H \times G \times K$ ;
\item Il existe une partie compacte~$C$ de~$G$ telle que~$f$ soit nulle en dehors de~$H\,{}C$. 
\end{conditions}
Le groupe~$G$ agit sur~$V$ par translations à droite ; on obtient ainsi une représentation lisse de~$G$ sur~$V$ (« induite compacte »),
notée $\Ind_H^G(\rho)$ ou $\cInd_H^G(\rho)$. Nous utiliserons principalement le procédé d'induction compacte dans le cas où~$H$ est compact ouvert modulo le centre ; la seconde notation est alors plus habituelle. 

\subsubsection{Représentations supercuspidales irréductibles}
Considérons cependant le cas où l'on induit à partir d'un sous-groupe \emph{parabolique} de~$G$.
Soit~$\P$ un sous-groupe parabolique de~$\G$ défini sur~$F$
\footnote{Rappelons qu'un sous-groupe de Borel de~$\G$ est un sous-groupe algébrique  de~$\G$ qui est fermé, connexe, résoluble, et maximal pour ces propriétés ;
et qu'un sous-groupe algébrique~$\P$ de~$\G$ est \emph{parabolique} s'il contient un sous-groupe de Borel.
Si~$\N$ est le radical unipotent de~$\P$, un \emph{facteur de Levi} de~$\P$ est un sous-groupe algébrique réductif connexe~$\L$ de $\P$, contenu dans le normalisateur~$\mathrm{Norm}_{\P}(\N)$, qui vérifie $\L \cap \N = \{1\}$ et $\P=\L\N$.
Si~$\P$ est défini sur~$F$, c'est aussi le cas de~$\L$ et~$\N$ ;
le groupe~$N=\N(F)$ est le radical unipotent de $P=\P(F)$, et le groupe $L=\L(F)$ vérifie alors $P=LN$ et $L\cap N=\{1_G\}$ («\,{}décomposition de Levi\,{}»).
Un \emph{sous-groupe de Levi} de~$\G$ est un sous-groupe réductif algébrique connexe de~$\G$  qui apparaît comme facteur de Levi d'un sous-groupe parabolique de~$\G$. S'il est défini sur~$F$, on dit que~$L=\L(F)$ est un sous-groupe de Levi de~$G$.
}
; notons $\N$ le radical unipotent de~$\P$ et considérons les groupes $P=\P(F)$ et~$N=\N(F)$.
Soit~$\rho$  une représentation lisse irréductible de~$P$ triviale sur~$N$ ;
elle s'identifie donc à une représentation lisse d'un sous-groupe de Levi~$L$ de~$G$.
La représentation~$\Ind_P^G(\rho)$ est alors de longueur finie (et l'espace $G/P$ est compact, donc la condition (ii) ci-dessus est vide). 

On dit qu'une représentation lisse irréductible~$\pi$ de~$G$ est \emph{supercuspidale}
si elle ne peut apparaître comme sous-quotient irréductible d'une telle induite $\Ind_P^G(\rho)$ pour~\mbox{$P \neq G$}.

On sait que~$G$ admet toujours des représentations supercuspidales.
Il fut longtemps difficile de donner une référence, mais dans le cas où~$F$ est de caractéristique~$0$, une preuve très courte a été rédigée par \textcite{BeuzartPlessis_Pacific}.

Comme évoqué au \no\hyperref[sec:reduction_jacquet]{1.2},
pour toute représentation lisse irréductible~$\pi$ de~$G$,
il existe un couple $(P, \sigma)$, où $P$ est un sous-groupe parabolique de~$G$ et~$\sigma$ est une représentation supercuspidale irréductible d'un facteur de Levi~$L$ de~$P$,
tel que~$\pi$ soit équivalente à un sous-module irréductible de $\Ind_P^G(\sigma)$. (L'induction se fait ici à partir de la représentation de~$P$ qui étend trivialement~$\sigma$ le long de la décomposition de Levi $P=LN$.)

\subsubsection{Lien avec les coefficients matriciels}
Si~$\pi$ est une représentation lisse de~$G$ sur un espace vectoriel~$V$, si~$\beta$ est une forme linéaire sur~$V$ et si~$v$ est un élément de~$V$, 
on note~$f_{\beta, v}$ la fonction de~$G$ dans~$\C$ définie par \mbox{$f_{\beta, v}(g) = \dual{\beta}{\pi(g)v}$}.

Pour obtenir ainsi des fonctions continues, on impose une condition sur~$\beta$.
Soient~$V^\ast$ l'espace dual de~$V$ et $\pi^\ast\colon G\to \GL(V^\ast)$ la représentation $g \mapsto {}^t{\pi(g^{-1})}$.
Cette dernière n'est pas lisse en général.
Pour tout sous-groupe compact ouvert~$K$ de~$G$, notons~$(V^\ast)^K$ l'espace des vecteurs~$\pi^\ast(K)$-invariants de~$V^\ast$,
et formons l'espace~\mbox{$\widetilde{V} = \bigcup_{K} (V^\ast)^K$}, la réunion étant prise sur l'ensemble des sous-groupes compacts ouverts de~$G$.
L'espace~$\widetilde{V}$ est stable par~$\pi^\ast$ ; on obtient ainsi une représentation lisse de~$G$.
Pour tout~$\beta$ dans~$\widetilde{V}$ et pour tout~$v$ dans~$V$, la fonction~$f_{\beta, v}\colon G \to \C$ est alors continue, et même localement~constante.
On dit que les fonctions $f_{\beta, v}$ sont les \emph{coefficients matriciels} de~$\pi$. 

Si~$\pi$ est une représentation lisse et irréductible de~$G$,
alors elle est supercuspidale si et seulement si chacun de ses coefficients matriciels est à support compact modulo~$Z(G)$,
c'est-à-dire nul en dehors d'une partie de~$G$ de la forme $C\cdot Z(G)$ avec~$C$ compacte.

Plus généralement, je dirai qu'une représentation lisse~$\pi$ de~$G$ est \emph{supercuspidale}
si elle est admissible et si tous ses coefficients matriciels sont à support compact modulo~$Z(G)$.
Si~$\pi$ est réductible, certains auteurs parlent plutôt de \mbox{«\,{}représentation compacte\,{}».}

\subsubsection{Série discrète} \label{sec:def_discrete}
Si~$\pi\colon G\to \GL(V)$ est irréductible, alors sa restriction au centre de~$G$ est de la forme~$z \mapsto \omega_\pi(z) 1_V$
où~$\omega_\pi$ est un caractère de~$Z(G)$ (c'est-à-dire un homomorphisme continu de~$Z(G)$ dans~$\C^\ast$), appelé le \emph{caractère central} de~$\pi$. 

Supposons que~$\omega_\pi$ soit unitaire, c'est-à-dire à valeurs dans le cercle unité.
Alors pour tout~$v$ dans~$V$ et tout~$\beta$ dans~$\widetilde{V}$, la fonction $\abs{f_{\beta,v}}$ induit une fonction sur~$G/Z(G)$.
On fixe une mesure de Haar sur le groupe unimodulaire~$G/Z(G)$, et on dit que~$\pi$ est \emph{de série discrète} s'il existe un couple $(v, \beta)$ tel que $\abs{f_{\beta,v}}^2$ soit non nulle et intégrable sur~$G/Z(G)$.
Cela ne dépend pas du choix de mesure de Haar sur~$G/Z(G)$
et la fonction~$\abs{f_{\beta,v}}$ est alors de carré intégrable sur~$G/Z(G)$ pour tout couple $(v, \beta) \in V \times \widetilde{V}$.

\begin{smallrema}
Dans ce cas, on peut identifier~$\pi$ à une représentation unitaire.
Fixons~\mbox{$\beta \in \widetilde{V}\setminus\{0\}$} et une mesure de Haar sur~$G/Z(G)$ ;
posons $\langle v_1, v_2 \rangle = \int_{G/Z(G)} f_{\beta, v_1} \overline{f_{\beta, v_2}}$ pour $v_1, v_2 \in V$.
Cela définit un produit scalaire $G$-invariant sur~$V$.
En complétant, on obtient une représentation unitaire irréductible de~$G$. On peut donc voir~$\pi$ comme une représentation de carré intégrable modulo le centre,
au sens où l'entend \mbox{\textcite[\S\,{}2,\no8]{Bourbaki_TSV}.}
\end{smallrema}

Si l'on ne suppose plus~$\omega_\pi$ unitaire, alors il existe un unique caractère $\chi\colon G \to \R^\ast_{+}$ tel que le caractère central de~$\pi \otimes \chi^{-1}$ soit unitaire,
et on dit que~$\pi$ est de série discrète si c'est le cas pour $\pi \otimes \chi^{-1}$.

La classification des séries discrètes reste mystérieuse dans le cas non archimédien.
Hormis le cas de~$\GL(n)$ (Bernstein et Zelevinsky ; voir l'exposé de \mbox{\cite{Rodier_1981_seminaire}}),
les meilleurs résultats connus sont probablement ceux de~\textcite{Moeglin}, puis \mbox{\textcite{Moeglin_Tadic}}, pour les groupes classiques déployés.
Ces résultats décrivent les séries discrètes à partir des représentations supercuspidales des groupes~\mbox{$\GL(m)$, $m \leq n$.}

\subsection{Tores, racines, groupes de Weyl}

\subsubsection{} \label{sec:tores}
Un \emph{tore} est un groupe algébrique isomorphe à un produit de groupes $\GL(1)$.
Si~$\T$ est un tore contenu dans~$\G$ et défini sur~$F$, nous parlerons parfois de tore  dans~\mbox{$G=\G(F)$} pour le groupe~$T = \T(F)$.
Les tores maximaux de~$\G$ ont tous la même dimension.
Si l'on fixe une clôture algébrique~$\ClotureF$ de~$F$, alors tous les tores maximaux de~$G$ sont conjugués sous~$\G(\ClotureF)$, mais pas nécessairement sous~$G$.

\subsubsection{Tores et groupes déployés} 
Soit~$\T$. On dit que~$\T$ est \emph{déployé} sur~$F$, ou \emph{se déploie} sur~$F$, s'il est défini sur~$F$ et si~$T=\T(F)$ est isomorphe à $(F^\times)^r$, où~$r \geq 1$ est l'entier tel que~$\T$ soit isomorphe (en tant que groupe algébrique) à $\GL(1)^r$.
On dit que~$\G$ se déploie sur~$F$ si l'un de ses tores maximaux est déployé sur~$F$. 

\subsubsection{Tores maximaux elliptiques} \label{sec:def_elliptique}
Si~$\T$ est un tore maximal de~$\G$, on dit que~$\T$ est \emph{$F$-elliptique}, ou plus simplement \emph{elliptique}, s'il est défini sur~$F$ et si $T=\T(F)$ est compact modulo le centre de~$G$.
Si~$F$ n'est pas archimédien, alors~$\G$ admet toujours un tore maximal elliptique.
Si~$F=\C$, il n'en admet jamais, sauf si~$\G$ est un tore.
Si~\mbox{$F=\R$}, il existe des groupes~$\G$ admettant des tores maximaux $F$-elliptiques (comme~$\SL(2)$)
et d'autres qui n'en admettent pas (comme $\SL(n)$ si $n > 2$).

\subsubsection{Caractères de tores ; racines}
Si~$\T$ est un tore défini sur~$F$, on note $\Char(T)$ le groupe des caractères de $T=\T(F)$ (homomorphismes continus~\mbox{$T \to \C^\ast$}).

On utilise aussi des notions plus algébriques : on note $\Char^\ast(\T)$ le groupe des caractères algébriques de~$\T$ (morphismes algébriques \mbox{$\T \to \GL(1)$})
et  $\Char_\ast(\T)$ le groupe des cocaractères algébriques de~$\T$ (morphismes algébriques \mbox{$\GL(1)\to\T$}).
Ce sont des groupes abéliens libres de type fini.
Si~$F=\R$, on peut identifier~$\Char_\ast(\T)$ à une partie de l'algèbre de Lie~$\lie{t}_\C$ de~$\T$, et~$\Char^\ast(\T)$ à une partie de son dual~$\lie{t}_\C^\ast$. 

Si~$\T$ est un tore maximal de~$\G$ défini sur~$F$, alors les racines de $\T$ dans~$\G$ forment une partie finie de~$\Char^\ast(\T)$. 
Chaque racine $\alpha$ vient avec une coracine~$\alpha^\vee \in \Char_\ast(\T)$.

\subsubsection{Groupes de Weyl}\label{sec:weyl}
Soient~$\T$ un tore maximal de~$\G$ et~$\mathrm{Norm}_{\G}(\T)$ le normalisateur de~$\T$ dans~$\G$.
Le quotient~$W(\G, \T)=\mathrm{Norm}_{\G}(\T)/\T$ est un groupe algébrique fini, appelé \emph{groupe de Weyl absolu} de~$(\G, \T)$.
Son ordre ne dépend pas du choix de~$\T$ et s'il n'y a pas d'ambiguïté, on écrit~$W$ plutôt que $W({\G, \T})$.
Pour toute racine~$\alpha$ de~$\T$ dans~$\G$, on dispose d'un élément~$s_\alpha$ d'ordre 2 dans~$W$, et ces réflexions engendrent~$W$.
On note~\mbox{$\eps\colon W\to \{\pm 1\}$} la \emph{signature},  unique homomorphisme de~$W$ dans~$\C^\ast$ vérifiant~\mbox{$\eps(s_\alpha)=-1$} pour chacune des réflexion~$s_\alpha$. 

Si~$\T$ est défini sur~$F$, on peut aussi considérer le groupe $W(G,T) = \mathrm{Norm}_{G}(T)/T$.
Il s'identifie à un sous-groupe de~$W(\G,\T)$, on peut donc lui restreindre la signature.

\subsection{Caractères des représentations (cas non archimédien)}\label{sec:defs_caracteres}

\subsubsection{Caractère-distribution}
Si~$\pi$ est une représentation lisse de~$G$ sur un espace vectoriel~$V$,
alors~$V$ est  généralement de dimension infinie, et~$\pi$ n'a pas de caractère au sens usuel :
les opérateurs~$\pi(g)$, $g \in G$, n'ont pas de trace en général. 

Cependant, fixons une mesure de Haar sur~$G$.
Notons $C^\infty_c(G)$ l'ensemble des fonctions de~$G$ dans~$\C$ qui sont localement constantes et à support compact.
Alors pour tout élément~$f$ de~$C^\infty_c(G)$, on peut définir l'endomorphisme $\pi(f)$ par la formule
\mbox{$\pi(f)(v) = \int_{G} \pi(g)(v) f(g) dg$} pour~$v \in V$ (l'intégrale se réduit à une somme finie).
Si~$\pi$ est admissible (notamment si~$\pi$ est irréductible), alors $\pi(f)$ est de rang fini pour tout~$f$.
On peut donc considérer la forme linéaire \mbox{$\Theta_\pi^D\colon f \mapsto \Tr(\pi(f))$ sur~$C^\infty_c(G)$} ; c'est le \emph{caractère-distribution} de~$\pi$.
Il ne dépend que de la classe d'équivalence de~$\pi$. Si~$\pi, \pi'$ sont des représentations lisses de même caractère-distribution, elles sont équivalentes.

\subsubsection{Caractères sur les éléments réguliers}\label{sec:caractere_sur_reg}
Si~$\pi$ est une représentation lisse admissible, alors la distribution~$\Theta_\pi^D$ peut être représentée par une fonction localement constante,
sinon sur~$G$ tout entier, du moins sur l'ensemble des éléments \emph{réguliers} de~$G$.

On introduit une fonction~$\Delta\colon G \to \R$, le discriminant de Weyl,  en considérant l'action adjointe de~$G$ sur son algèbre de Lie.
Notons $\lie{g}=\mathrm{Lie}(\G)(F)$ ; c'est un espace vectoriel de dimension finie sur~$F$.
Pour~$g \in G$, on considère les puissances extérieures de l'endomorphisme~$\mathrm{Ad}(g)$ de~$F$ et on note $\Delta(g)=\abs{\Tr(\goodwedge^{\!d} \,{}\mathrm{Ad}(g))}_F$
où~$d = \dim(\G)-\mathrm{rg}(\G)$, le rang~$\mathrm{rg}(\G)$ étant la dimension commune des tores maximaux de~$\G$ (\no\ref{sec:tores}),
et où on utilise la valeur absolue normalisée sur~$F$.
On appelle \emph{ensemble régulier }de~$G$, et on note~$\Greg$, l'ensemble $G \setminus \Delta^{-1}(0)$.
C'est un ouvert dense de~$G$, invariant par conjugaison.

 Harish-Chandra montre qu'il existe une unique fonction~$\Theta_\pi\colon \Greg \to \C$
qui soit localement constante et vérifie $\dual{\Theta_\pi^D}{f} = \int_G \Theta_\pi f$ pour toute fonction $f \in C^\infty_c(\Greg)$.
Si~$\pi$ est de série discrète, elle se prolonge d'ailleurs en une fonction localement intégrable sur~$G$ de telle façon que cette identité soit satisfaite pour toute  $f \in C^\infty_c(G)$ :
Harish-Chandra n'avait démontré cela que pour~$F$ de caractéristique~$0$ (voir les exposés de \mbox{van~\cite{VanDijk}} et~\cite{Cartier_1977_seminaire}),
mais Beuzart-Plessis vient d'annoncer une preuve dans le cas général. 


\subsection{Le cas archimédien}\label{sec:gen_arch}
 Si~$F$ est archimédien, la notion adéquate de représentation lisse irréductible est moins commode et moins algébrique.
 Sans détails, disons que~$\pi$ est \emph{lisse}
 si~$V$ est un espace de Fréchet et~$\pi$ est~$C^\infty$, admissible, à croissance modérée \parencite{Casselman_croissance} ; 
 disons que~$\pi$ est \emph{irréductible} si~$V \neq \{0\}$ et n'a pas de sous-espace invariant \emph{fermé} non trivial.
 On note encore~$\Irr(G)$ l'ensemble des classes d'équivalence de représentations lisses irréductibles.
 En pratique, sauf au \S\,{}\ref{sec:LLC_generalites} pour des généralités sur la correspondance de~Langlands, nous ne nous intéresserons qu'aux représentations de la série discrète.
 Celles-ci sont définies comme au \no\ref{sec:def_discrete}, mais en considérant les coefficients matriciels $K$-finis.
 Il n'en existe que si~$G$ admet un tore maximal elliptique, et donc jamais si~$F=\C$ en dehors de cas triviaux.
 Les caractères sont définis comme au \S\,{}\ref{sec:defs_caracteres}, à ceci près que $C^\infty_c(G)$ désigne l'espace des fonctions $C^\infty$ à support compact au sens usuel.
 Pour $f \in C^\infty_c(G)$, l'endomorphisme~$\pi(f)$ est nucléaire.
 Je renvoie pour tout cela au livre de~\textcite{Wallach}.


\section{Caractères et construction des séries discrètes des groupes réels}\label{sec:cas_reel}

Esquissons le tableau de la série discrète si~$F=\R$.
Rappelons que~\mbox{$G=\G(\R)$} n'est pas connexe en général, mais possède un nombre fini de composantes connexes.
On fixe un tore maximal~$\T$ de~$\G$, défini sur~$\R$, et on pose~$T = \T(\R)$. 

\subsection{Formule du caractère de Weyl et revêtements doubles de tores}

Supposons d'abord que le groupe~$G$ soit compact.
On rappelle ici la formule de H.~Weyl pour les caractères des représentations irréductibles de~$G$.
Nous verrons qu'elle prend une forme particulièrement simple si l'on introduit un \emph{revêtement double} de~$T$.

\subsubsection{} \label{Weyl_entier}
Pour toute fonction~$f\colon T\to \C$, définissons une fonction $J_f\colon T \to \C$ par 
\begin{equation}\label{def_J} J_{f}(t) = \sum \limits_{w \in W(G,T)} \epsilon(w) \cdot f(w^{-1}t)\end{equation}
où~$\eps\colon W(G,T) \to\{\pm1\}$ est la signature (\no\ref{sec:weyl}).
Nous allons voir que les caractères des représentations irréductibles de~$G$ s'expriment comme des quotients de telles fonctions. 

Soit~$\B$ un sous-groupe de Borel de~$\G$ admettant~$\T$ comme facteur de Levi.
Soit~$\N$ le radical nilpotent de~$\B$ ; c'est un sous-groupe distingué de~$\B$.
L'action adjointe de~$\B$ sur~$\N$ induit une action de~$T$ sur l'algèbre de Lie~$\lie{n}$ de~$\N$ ;
notons $\delta_{\B}$ ou~$\delta$ la fonction de~$T$ dans~$ \C^\ast$ qui, à un élément~$t$ de~$T$, associe le déterminant $\det(\mathrm{Ad}(t)_{|\lie{n}})$.
Alors~$\delta$ est un élément de~$\Char(T)$, qui n'est autre que la restriction à~$T$ de la fonction module de~$\B$.
Dans ce numéro, faisons comme~\textcite[\S\,{}7, \no4]{Bourbaki_Lie9} l'hypothèse suivante : 
\begin{equation} \label{hyp_integralite} \text{\emph{Il existe un caractère continu $\rho\colon  T \to \C^\ast$ vérifiant $\rho^2=\delta$}}.\end{equation} 

La fonction~$J_\rho$ ne s'annule pas sur l'ensemble~$T_{\mathrm{reg}}=T\cap\Greg$ des éléments réguliers de~$T$.
Pour~$\lambda$ dans~$\Char(T)$, on peut donc former la fonction $\Theta_\lambda\colon T_{\mathrm{reg}} \to \C$ définie par 
\begin{equation} \label{formule_caractere} \Theta_{\lambda}(t)=   J_{\lambda\rho}(t)/J_{\rho}(t).\end{equation}
Jusqu'à la fin de ce numéro, supposons de plus que~$G$ soit \emph{connexe}.
Alors les fonctions~${\Theta}_{\lambda}$ sont les restrictions à $T_\mathrm{reg}$ des caractères des représentations irréductibles de~$G$ (\cf~\cite[\S\,{}7]{Bourbaki_Lie9}).
Plus précisément : 
\begin{assertions}
\item Pour tout élément~$\lambda$ de~$\Char(T)$, il existe une représentation irréductible~$\pi_{\lambda}$ de~$G$, unique à équivalence près, dont le caractère coïncide sur~$T_{\mathrm{reg}}$ avec la fonction~$\Theta_{\lambda}$ ;
\item Toute représentation irréductible de~$G$ est équivalente à une représentation~$\pi_{\lambda}$.
\end{assertions}
Ces résultats ne donnent pas une idée très concrète des représentations irréductibles. 
Pour cela, il faudrait au moins compléter~\ass{a}
en expliquant comment  construire concrètement la représentation~$\pi_{\lambda}$ à partir d'un caractère~$\lambda$ («\,théorème de Borel--Weil\,»).
Il faudrait aussi, c'est plus facile, compléter~\ass{b}
en associant concrètement à tout élément~$\pi$ de~$\Irr(G)$ un élément~$\lambda$ de~$\Char(T)$ vérifiant~$\pi\simeq\pi_\lambda$ («\,{}théorie du plus haut poids\,{}»). 

À propos de~\ass{b}, contentons-nous de dire qu'ayant fixé un sous-groupe de Borel~$\B$ de~$\G$, on dispose d'un ordre partiel sur~$\Char(T)$.
\begin{smallrema}
En effet, le choix de~$\B$ détermine un cône  (chambre de Weyl) dans le dual~$\lie{t}^\ast$ de l'algèbre de Lie de~$T$.
Ce cône définit un ordre partiel sur~$\lie{t}^\ast$,
et on compare les éléments de $\Char(T)$ en comparant leurs différentielles, vues comme des éléments de $\lie{t}^\ast$.  
\end{smallrema}
Si~$\pi$ est une représentation irréductible de~$G$,
l'ensemble des éléments de~$\Char(T)$ apparaissant dans la restriction de $\pi$ à $T$ admet alors un plus grand élément~$\lambda_{\mathrm{dom}}$,
le \emph{plus haut poids} de~$\pi$.
On a $\pi \simeq \pi_{\lambda_{\mathrm{dom}}}$ (\cf\cite[\S\,{}7, \no1]{Bourbaki_Lie9}).

\subsubsection{} \label{Weyl_general}
Indiquons maintenant comment se passer de l'hypothèse~\eqref{hyp_integralite}.
Il n'existe certes pas toujours de caractère continu de~$T$ donnant une racine carrée de~$\delta$ ;
mais pour pouvoir extraire une telle racine, il suffit de remplacer~$T$ par un \emph{revêtement double}. 

Notons $T_{\pm}$ l'ensemble des couples $(t,z)\in T \times \C^\ast$ vérifiant $\delta(t) = z^2$.
C'est un sous-groupe du produit direct $T \times \C^\ast$, et la projection sur le premier facteur fait de~$T_\pm$ un revêtement à deux feuillets de~$T$.

Les caractères continus de~$T_{\pm}$ sont de deux sortes :
ceux qui sont triviaux sur le noyau de la projection~\mbox{$p_{\pm}\colon T_{\pm} \to T$}, et sont donc de la forme~$\lambda \circ p_{\pm}$ où~$\lambda$ est un caractère de~$T$ ;
et ceux qui ne sont pas triviaux sur~$\ker(p_{\pm})$, qu'on appelle caractères \emph{spécifiques} de~$T_{\pm}$.
On note $\Charspec(T_{\pm})$ l'ensemble des caractères spécifiques de~$T_{\pm}$. 

Pour~$\tilde{s}=(t, z)$ dans~$T_\pm$, on pose~$\rho(\tilde{s}) = z$.
On définit ainsi un caractère spécifique~$\rho$ de~$T_{\pm}$ ; c'est une racine carrée de $\delta\circ p_{\pm}$.
Si~$\vartheta$ appartient à~$\Charspec(T_{\pm})$, alors le produit~$\vartheta \rho$ se~factorise par~$p_\pm$ ; la multiplication par~$\rho$ induit  une bijection entre~\mbox{$\Charspec(T_{\pm})$~et~$\Char(T)$.}

Pour donner un sens à la formule~\eqref{def_J} en présence du revêtement~$T_{\pm}$, on relève à~$T_\pm$ l'action de~$W(G,T)$ sur~$T$.
Si~$w$ est un élément de~$W(G,T)$, et si~$(t,z)$ appartient à~$T_{\pm}$, alors $(w^{-1} t, z)$ n'appartient pas nécessairement à~$T_{\pm}$.
Cependant, il existe un caractère $\zeta_w$ de~$T$ tel que $(w^{-1} t, \zeta_w(t) z)$ appartienne à~$T_{\pm}$ pour tout~$(t,z) \in T_{\pm}$ :
en effet, le quotient~$\delta_{w^{-1} \B w}/\delta_\B$ est le carré d'un produit de racines de~$(\G,\T)$, qui définit le caractère $\zeta_w$ annoncé.
On obtient ainsi une action de~$W(G,T)$ sur~$T_{\pm}$.

Pour toute fonction~$f\colon T_{\pm} \to \C$, la formule~\eqref{def_J} définit alors une fonction $J_f\colon T_{\pm}\to \C$,
et~$J_{\rho}$ ne s'annule pas sur l'image réciproque $p_{\pm}^{-1}(T_{\mathrm{reg}})$.
Si~$\vartheta$ est un caractère spécifique de~$T_{\pm}$,
la fonction $\tilde{s} \mapsto J_{\vartheta}(\tilde{s})/J_{\rho}(\tilde{s})$ se factorise par la projection~$p_\pm$,
et définit donc une fonction~$\Theta_{\vartheta}$ sur~$T_{\mathrm{reg}}$.
Les résultats du~\S\ref{Weyl_entier} se reformulent pour~$G$~connexe, et s'étendent en général, ainsi : 
\begin{prop}\label{prop:cas_compact}
\begin{assertions}
\item Pour tout élément~$\vartheta$ de~$\Charspec(T_\pm)$, il existe une unique représentation irréductible~$\pi_{\vartheta}$ de~$G$ dont le caractère coïncide sur~$T_{\mathrm{reg}}$ avec~$\Theta_{\vartheta}$\,\pv
\item L'application~$\vartheta \mapsto \Theta_{\vartheta}$ induit une bijection de $\Charspec(T_{\pm})/W(G,T)$ sur l'ensemble des restrictions à~$T_{\mathrm{reg}}$ de caractères de représentations irréductibles de~$G$\pt
\end{assertions}
\end{prop}


\subsection{Caractères des séries discrètes des groupes réels}\label{sec:caractere_HC}

Abordons à présent le cas où~$G$ n'est pas compact.
La série discrète de~$G$ est non vide si et seulement si~$\G$ admet un tore maximal \emph{elliptique} ;
on suppose désormais que~$\T$ est elliptique.
Si~$\Gder$ est le groupe dérivé de~$\G$, alors il existe un unique sous-groupe compact maximal~$K$ de~$G$ tel que~$T\cap \Gder$ soit un sous-groupe de Cartan de~$K$ (centralisateur d'un tore maximal de~$K$). L'espace symétrique~$G/K$ est de dimension~paire. 

Étant donné un sous-groupe de Borel~$\B$ de~$\G$ admettant~$\T$ comme facteur de Levi, on peut toujours former le revêtement double~$T_{\pm}$
associé à la fonction module~$\delta$ de~$\B$,
et considérer l'ensemble $\Charspec(T_{\pm})$ des caractères spécifiques de~$T_{\pm}$.
Il contient un élément~$\rho$ donnant une racine carrée de~$\delta$.
On peut toujours associer à  toute fonction~\mbox{$f\colon T_{\pm}\to\C$} une fonction $J_f\colon T_{\pm}\to\C$, par la formule~\eqref{def_J}.
La fonction~$J_{\rho}$ ne s'annule pas.
Si~$\vartheta$ est un caractère spécifique de~$T_{\pm}$, alors~\mbox{$(J_{\vartheta}/J_{\rho})$} se factorise par la projection $T_{\pm}\to T$ ;
on peut donc considérer la fonction 
\begin{equation} \label{formule_caractere_2} \Theta_{\vartheta}(s)=   (-1)^{q} (J_{\vartheta}/J_{\rho})(s)\quad \text{pour~$s \in T_{\mathrm{reg}}$,}\end{equation}
où~$q$ est la moitié de l'entier pair $\dim(G/K)$.

Les groupes $T_{\pm}$ et $T$ ayant la même algèbre de Lie,
on peut identifier la différentielle d'un caractère de $T_{\pm}$ ou de $T$
à un élément de l'espace dual $\lie{t}^\ast$ de l'algèbre de Lie de~$T$.
On dit qu'un caractère spécifique de~$T_{\pm}$ est \emph{régulier} si sa différentielle, ainsi considérée comme un élément de $\lie{t}^\ast$,
ne s'annule sur aucune des coracines de~$(\G,\T)$.
Notons~$\Charspecreg(T_{\pm})$ l'ensemble des caractères spécifiques réguliers de~$T_{\pm}$.

\begin{theo}[\cite{HarishChandra_SerieDiscrete}]\label{HC_parametrage}
\begin{assertions}
\item
Si~$\vartheta$ est un caractère spécifique régulier de~$T_\pm$, alors il existe une représentation de série discrète de~$G$, unique à équivalence près, dont le caractère coïncide sur~$T_{\mathrm{reg}}$ avec la fonction~$\Theta_{\vartheta}$\,\pv
\item
L'application~$\vartheta \mapsto \Theta_{\vartheta}$ induit une bijection
de $\Charspecreg(T_{\pm})/W(G,T)$ sur l'ensemble des restrictions à~$T_{\mathrm{reg}}$ de caractères de représentations de série discrète de~$G$\pt
\end{assertions}
\end{theo}
Mentionnons une différence importante entre ce théorème et la prop.~\ref{prop:cas_compact}, qui en est un cas particulier. Si~$\pi$ est une représentation de série discrète de~$G$, alors le th.~\ref{HC_parametrage}  montre l'existence une formule simple pour le caractère de~$\pi$, mais seulement sur l'ensemble~$G_{\mathrm{ell}}$ des éléments de~$G$ conjugués à un élément de~$T$. La restriction à~$G_{\mathrm{ell}}$ du caractère de~$\pi$ suffit à la caractériser au sein de la série discrète.
Si~$G$ est compact on a \mbox{$G_{\mathrm{ell}}=G$},  mais en général, il n'y a pas de formule simple pour le caractère de~$\pi$ sur~$G \setminus G_{\mathrm{ell}}$. 


\subsection{Construction des représentations de série discrète}


Le th.~\ref{HC_parametrage} fournit un paramétrage de la série discrète de~$G$ par les caractères spécifiques réguliers de~$T_{\pm}$. 
Mais si $\vartheta$~est un caractère spécifique régulier de~$T_{\pm}$, le résultat d'Harish-Chandra ne fournit pas de réalisation concrète de la représentation~$\pi_\vartheta$ correspondante. 
Plusieurs réalisations ont été découvertes dans les années 1970 ; l'exposé de \textcite{Duflo_1977_seminaire} faisait le point sur ces constructions.
L'une d'elles dépend d'un lien entre séries discrètes et représentations de sous-groupes compacts,
qu'il me semble utile d'esquisser en vue du \S\,{}\ref{sec:Yu_et_Fintzen}.
Pour simplifier et utiliser le \no\ref{Weyl_entier}, on suppose~$G$~connexe. 

Fixons un sous-groupe compact maximal~$K$ de~$G$ comme dans le \S\,{}\ref{sec:caractere_HC}. 
Si~$\pi$ est une représentation de série discrète de~$G$, alors la restriction~$\pi_{|K}$ est somme hilbertienne de représentations irréductibles.
Supposons toujours fixé un sous-groupe de Borel $\B$ de $\G$ et considérons l'ensemble des plus hauts poids des éléments de~$\Irr(K)$ qui apparaissent dans $\pi_{|K}$.
Cet ensemble admet un plus petit élément (\cf\cite{Duflo_1977_seminaire}) ; la représentation correspondante~$\mu_{\pi}$ apparaît avec multiplicité~$1$ dans~$\pi_{|K}$.
On dit que~$\mu_{\pi}$ est le \emph{$K$-type minimal} de~$\pi$. Il ne dépend que de~$\pi$ et~$K$ et non  du choix de~$\B$.
Si~$\pi$ et~$\pi'$ sont des représentations de série discrète de~$G$ qui ne sont pas équivalentes, alors~\mbox{$\mu_{\pi} \neq \mu_{\pi'}$}.

Notons $\Irr_{\mathrm{reg}}(K)$ l'ensemble des éléments de~$\Irr(K)$ apparaissant comme le $K$-type minimal d'une représentation de série discrète.
Si~$\mu$ est un élément de~$\Irr_{\mathrm{reg}}(K)$, alors il est possible de donner une réalisation géométrique de la représentation de série discrète correspondante
en utilisant l'équation de Dirac sur la variété riemannienne~$G/K$,
suivant une méthode initiée par~\textcite{Parthasarathy}, amplifiée par~\textcite{AtiyahSchmid}. 


\begin{smallrema}
Nous n'aurons pas besoin ci-dessous du détail de la construction : \cf~\textcite{Duflo_1977_seminaire}.
Disons seulement qu'en munissant~$G/K$ d'une métrique riemannienne~$G$-invariante, on obtient une variété spinorielle ;
à partir du fibré des spineurs~$\mathscr{S}$ sur~$G/K$ et de la représentation~$\mu$, on construit un fibré équivariant~$\mathscr{E}$ sur~$G/K$,
dont la fibre à l'origine est une représentation de dimension finie de~$K$, qui contient~$\mu$ avec multiplicité~$1$.
Le fibré~$\mathscr{E}$ est canoniquement muni d'un opérateur de Dirac $G$-équivariant, qui agit sur les sections~$C^\infty$ de~$\mathscr{E}$.
L'ensemble des sections de~$\mathscr{E}$ qui sont~$C^\infty$, de carré intégrable, et  annulées par l'opérateur de Dirac, est alors une réalisation géométrique de~$\pi$. 
\end{smallrema}


En résumé, si l'on se donne un couple $(K, \mu)$ où~$K$ est un sous-groupe compact maximal de~$G$ et~$\mu$ appartient à $\Irr_{\mathrm{reg}}(K)$,
alors la construction d'Atiyah--Schmid fournit une représentation de série discrète~$\pi$ de~$G$.
Il n'est pas difficile de déterminer les couples~$(K, \mu)$ pertinents ;
en revanche, la construction de~$\pi$ à partir de~$\mu$, et la démonstration du fait qu'on obtient ainsi toutes les séries discrètes,
utilisent des ingrédients relativement sophistiqués (opérateur de Dirac, théorie de l'indice $L^2$).


\section{Classification des représentations supercuspidales dans le cas modéré}\label{sec:Yu_et_Fintzen}

Dans ce paragraphe on suppose~$F$ non archimédien, de caractéristique résiduelle~\mbox{$p>2$}.
Étant donné un couple~$(K, \mu)$ où~$K$ est un sous-groupe compact ouvert modulo le centre de~$G$ et $\mu$ est une représentation lisse irréductible de~$K$,
le procédé d'induction compacte du \no\ref{sec:induction} fournit une représentation lisse de~$G$. 
Cela donne une manière beaucoup plus simple qu'au \S\,{}\ref{sec:cas_reel} d'associer à~$(K, \mu)$ une représentation de~$G$ ; mais les remarques suivantes montrent qu'on peut espérer obtenir ainsi des représentations~supercuspidales :  
\begin{assertions}
\item Si~$\pi = \cInd_{K}^{G}(\mu)$ est irréductible, alors elle est nécessairement supercuspidale ; 
\item\label{critere_irred} Pour que~$\pi$ soit irréductible, il suffit que pour $g \in G \setminus K$,
les restrictions à $K \cap g K g^{-1}$ des représentations $\mu$ et $x \mapsto \mu(g^{-1} x g)$ n'aient aucun constituant commun.
\footnote{Cette formulation est empruntée à l'exposé d'\textcite{Henniart_1991_seminaire}.} 
\end{assertions}

On conjecture depuis longtemps 
que toute représentation supercuspidale~$\pi$ de~$G$ peut s'écrire sous la forme $\pi=\cInd_{K}^{G}(\mu)$ pour un certain couple $(K, \mu)$.
Ce qui est difficile, c'est de trouver des couples~$(K,\mu)$ pertinents.
Nous verrons que de tels couples peuvent être construits à partir de la structure fine de l'immeuble de Bruhat--Tits de~$G$.
La construction la plus complète est due à~\textcite{Yu_2001} et généralise des constructions antérieures : la première est celle de Howe~\parencite*{Howe_1977} pour $\GL(n)$, dont la portée fut établie par \textcite{Moy_1986},
tandis qu'\textcite{Adler_Pacific} fournissait une construction déjà très générale à l'aide de la filtration de Moy--Prasad.

On suppose dans tout ce paragraphe que~$\G$ se déploie sur une extension modérément ramifiée de~$F$. 
Si de plus~$p$ ne divise pas l'ordre du groupe de Weyl  absolu de~$\G$ (\no\ref{sec:weyl}),
alors~\textcite{Fintzen_Annals} montre  que toute représentation supercuspidale irréductible~$\pi$ de~$G$ s'obtient à partir de la construction de Yu.
Sa démonstration fait même voir comment extraire de~$\pi$ un couple~$(K, \mu)$ vérifiant $\pi=\cInd_{K}^{G}(\mu)$.


\subsection{Les immeubles du groupe et la notion de profondeur des représentations}


\subsubsection{Immeubles}\label{sec:defs_immeubles}
Les informations les plus intéressantes sur les sous-groupes compacts ouverts de~$G$ sont encodées par la géométrie des \emph{immeubles} associés à~$G$.
Je ne définirai pas ici ces immeubles, et ne décrirai pas leur structure en détail. 
Des références canoniques sont les textes de \textcite{Bruhat_Tits_I} et le résumé de Tits~(\cite*{Tits_Corvallis},~\cite*{Corvallis}).
Le livre récent de~\textcite{Kaletha_Prasad} couvre largement les besoins de ce texte.
Enfin, les notes de~\textcite{Fintzen_Harvard} contiennent une présentation assez originale, conçue pour aborder notre sujet avec le moins possible de prérequis.

Soit~$\Imm(G)$ l'immeuble (étendu) de~$G$.
C'est un espace métrique sur lequel~$G$ agit par isométries.
Pour tout élément~$x$ de~$\Imm(G)$, le stabilisateur~$G_x$ est un sous-groupe compact ouvert de~$G$.
Réciproquement, tout sous-groupe compact ouvert maximal de~$G$ est égal à~$G_x$ pour un certain point~$x$ de~$\Imm(G)$.
L'immeuble~$\Imm(G)$ a des~\emph{sommets} : ce sont des points~$x$ tels que~$G_x$ soit un sous-groupe compact maximal.

L'immeuble étendu permet donc de construire de nombreux sous-groupes compacts ouverts de~$G$.
L'action de~$G$ sur~$\Imm(G)$ n'est pas transitive en général ;
si~$x$ et~$y$ sont deux points de~$\Imm(G)$, les stabilisateurs $G_x$ et~$G_y$ ne sont donc pas nécessairement conjugués. 

Nous aurons besoin de l'immeuble \emph{réduit} de~$G$.
Si~$\G_{\ad}$ est le quotient de~$\G$ par son centre $\Z_{\G}$, et si~$G_\ad$ est le groupe~$\G_\ad(F)$,
alors l'immeuble réduit~$\Immred(G)$ n'est autre que l'immeuble~$\mathscr{B}(G_\ad)$.
On peut écrire~$\Imm(G)$ comme produit direct de~$\Immred(G)$ et de l'espace affine $\Char_\ast(\Z_\G)_F\otimes_\Z \R$,
où $\Char_\ast(\Z_\G)_F$ est le groupe des cocaractères algébriques de~$\Z_\G$ définis sur~$F$.
Cela induit une surjection de $\Imm(G)$ sur~$\Immred(G)$ ;
si~$x$ est un élément de~$\Imm(G)$, on note~$\reduit{x}$ l'image de~$x$ dans l'immeuble réduit~$\Immred(G)$. 
Le groupe~$G$ agit sur~$\Immred(G)$. Si~$x$ est un point de~$\Imm(G)$, alors le stabilisateur~$G_{\reduit{x}}$ de~$\reduit{x}$ n'est plus nécessairement compact,
mais il est {compact modulo le centre}. 

\subsubsection{Filtration de Moy--Prasad} \label{sec:la_filtration}
Fixons un point~$x$ de~$\Imm(G)$.
\textcite{Moy_Prasad} ont décrit une filtration de~$G_x$ par des sous-groupes distingués qui jouent un rôle essentiel dans la construction de représentation supercuspidales.

La filtration de Moy--Prasad est une famille décroissante $(G_{x,r})_{r \in \R_+}$ de sous-groupes (compacts) ouverts de~$G_x$.
L'intersection $\bigcap_{r \geq 0} G_{x,r}$ est réduite à l'identité. 

 Pour tout~$r \geq 0$, on note $G_{x,r^+} = \bigcup_{s > r} G_{x,s}$.
 On a $G_{x,r^+}=G_{x,r}$ en dehors d'une suite $(r_{n})_{n \in \N}$ vérifiant~$r_0 = 0$.
 Connaître la filtration $(G_{x,r})_{r \in \R_+}$ revient donc à connaître une suite décroissante
\[ G_{x,0} \supset G_{x,r_1} \supset \dots \supset G_{x,r_n} \supset \dots\]
de sous-groupes de~$G_x$.
Le sous-groupe~$G_{x,0}$ avait été introduit par Bruhat et~Tits~\parencite*{Bruhat_Tits_I} ; c'est un sous-groupe \emph{parahorique} de~$G$.
Si~$\G$ est semi-simple et simplement connexe, alors~$G_{x,0}$ est égal à $G_x$.
En général, $G_{x,0}$ est l'intersection de~$G_x$ et du noyau de l'homomorphisme de Kottwitz (\cf~\cite[Chap.~11]{Kaletha_Prasad}).
Le sous-groupe~$G_{x,0^+}$ est le radical pro-$p$-unipotent de~$G_{x,0}$, tandis que le stabilisateur~$G_{\reduit{x}}$ du \no\ref{sec:defs_immeubles} est le normalisateur de~$G_{x,0}$ dans~$G$.

Les relations entre les sous-groupes de la filtration de Moy--Prasad jouent un rôle crucial dans ce qui suit.
La plus importante concerne les commutateurs entre éléments de groupes de la filtration.
Si~$A$ et~$B$ sont des groupes, notons $\crochets{A,B}$ l'ensemble des éléments de la forme $aba^{-1}b^{-1}$ pour $a \in A$ et $b \in B$.
Alors  : 
\begin{assertions}
\item[a)] Pour~$r \geq 0$ et~$s \geq 0$, on a $\crochets{G_{x,r}, G_{x,s}} \subset G_{x,r+s}$.
\end{assertions}
On en déduit aisément les propriétés suivantes :
\begin{assertions}
\item[b)] Pour~$r \geq 0$ et~$s \geq r$, le sous-groupe $G_{x,s}$ est distingué dans~$G_{x,r}$. 

En particulier, $G_{x,r^+}$ est distingué dans~$G_{x,r}$ ; et~$G_{x,r}$ est distingué dans~$G_{x,0}$.
\item[c)] Pour~$r > 0$, le quotient~$G_{x,r}/G_{x,r^+}$ est abélien.

\end{assertions}
Les groupes~$G_{x,r}/G_{x,r^+}$ sont tous \emph{finis} d'après le \no\ref{sec:sg_compacts}. Leur structure  est connue :  
\begin{assertions}
\item[d)] Pour $r>0$, le groupe abélien~$G_{x,r}/G_{x,r^+}$ peut être muni d'une structure d'espace vectoriel sur le corps résiduel~$k_F$.

\begin{smallrema}
En fait, Moy et Prasad introduisent une filtration $(\lie{g}_{x,r})_{r \in \R}$ de l'algèbre de Lie~$\lie{g} = \mathrm{Lie}(\G)(F)$ par des réseaux entiers,
et décrivent un isomorphisme explicite de~$G_{x,r}/G_{x,r^+}$ sur $\lie{g}_{x,r}/\lie{g}_{x,r^+}$ pour~$r >0$.
La structure de $k_F$-espace vectoriel sur~$G_{x,r}/G_{x,r^+}$ est obtenue par transport de structure à partir de celle de~$\lie{g}_{x,r}/\lie{g}_{x,r^+}$.
Plus généralement, pour tout réel~$s \in [r, 2r]$, on dispose d'un isomorphisme de $G_{x,r}/G_{x,s}$ sur $\lie{g}_{x,r}/\lie{g}_{x,s}$ («\,{}isomorphisme de Moy--Prasad\,{}»).
\end{smallrema} 
\item[e)] Le quotient $G_{x,0}/G_{x,0^+}$ est  un groupe réductif fini  : il s'identifie~à~$\Gfini_x(k_F)$, où  $\Gfini_x$ est un groupe réductif connexe défini sur~$k_F$. 
\end{assertions}


\subsubsection{Représentations de profondeur nulle} \label{sec:profondeur_nulle}
Les données qui permettent de construire des représentations supercuspidales font grand usage des quotients ci-dessus.
Par exemple, soit~$x$ un point de~$\Imm(G)$ tel que~$\reduit{x}$ soit un \emph{sommet} de $\Immred(G)$. Supposons donnée une représentation lisse irréductible de~$G_{\reduit{x}}$ triviale sur~$G_{x,0^+}$,
et supposons que~$\rho_{|G_{x,0}}$ contienne une représentation irréductible \emph{cuspidale}
\footnote{Si on identifie $\rho_{|G_{x,0}}$ à une représentation irréductible de $\Gfini_x(k_F)$, cela signifie, naturellement, qu'elle ne se plonge dans aucune  induite parabolique de représentation d'un sous-groupe de Levi~propre.}
 du groupe réductif fini~$G_{x,0}/G_{x,0^+}$.
Alors on vérifie aisément que la représentation $\cInd_{G_{\reduit{x}}}^{G}(\rho)$ est irréductible, donc supercuspidale (\cf \cite[th.~4.4.1]{Fintzen_Harvard}). 

Les représentations que nous venons de construire sont caractérisées par leur  \emph{profondeur},
une notion introduite par Moy et Prasad à partir de la filtration du \no\ref{sec:la_filtration}.

Les sous-groupes compacts ouverts~$G_{x,r}$, pour~$x \in \Imm(G)$ et $r \geq 0$, forment une base de voisinages de l'identité de~$G$.
Par conséquent, si~$\pi$ est une représentation irréductible lisse de~$G$ sur un espace vectoriel~$V$,
et si l'on considère l'ensemble~$P_\pi$ des réels~$r \geq 0$ tels qu'il existe~$x \in \Imm(G)$ vérifiant $V^{G_{x, r^+}} \neq \{0\}$,
alors~$P_\pi$ est non vide.
L'ensemble~$P_\pi$ admet en fait un plus petit élément~\parencite[\S\,{}5]{DeBacker_Michigan}, appelé \emph{profondeur} de~$\pi$.

Les représentations supercuspidales irréductibles de profondeur nulle sont précisément celles que l'on obtient par le procédé ci-dessus \mbox{\parencite[prop.~6.8]{Moy_Prasad_2}}.



\subsection{La construction de Yu}


\subsubsection{Données cuspidales à un étage et construction dans ce cas}\label{sec:Yu_1pas}

La construction de~Yu s'appuie sur les représentations de profondeur nulle, mais on considère ces représentations pour \emph{sous-groupe de Levi tordu}~$L$ de~$G$.
Par définition, il s'agit du groupe~\mbox{$L=\L(F)$} des points rationnels d'un sous-groupe~$\L$ de~$\G$ qui est défini sur~$F$
et facteur de Levi d'un sous-groupe parabolique~$\P$ de~$\G$,
le groupe~$\P$ étant supposé défini sur une extension modérément ramifiée de~$F$ (mais pas nécessairement sur~$F$).
\smallskip

Dans ce numéro, on décrit un cas particulier de la construction. Supposons donnés : 
\begin{itemize}
\item Un sous-groupe de Levi tordu~$L$ de~$G$, avec $L \neq G$ ;
\item Une représentation supercuspidale irréductible~$\sigma$ de~$L$,  de profondeur nulle ;
\item Un caractère $\phi$ de~$L$, que l'on suppose de profondeur non nulle. 
\end{itemize}

Sous certaines conditions sur~$L$ et~$\phi$, nous allons construire un couple~$(K,\mu)$
où~$K$ est un sous-groupe compact modulo le centre et $\mu$ une représentation lisse irréductible de~$K$.
L'induite compacte $\cInd_K^G(\mu)$ fournira une représentation supercuspidale irréductible. 

\paragraph*{$\bullet$ Le sous-groupe~$K$ et la forme générale de~$\mu$}
D'après le~\no\ref{sec:profondeur_nulle}, il existe un couple $(x, \rho)$ où $x$ est un point de~$\Imm(L)$ tel que~$\reduit{x}$ soit un sommet de~$\Immred(L)$,
et où~$\rho$ est une représentation lisse irréductible de~$L_{\reduit{x}}$ triviale sur~$L_{x,0^+}$,
tels qu'on ait~\mbox{$\sigma = \cInd_{L_{\reduit{x}}}^L(\rho)$.}
Soit~$r$ l'unique réel $\geq 0$ tel que~$\phi$ soit trivial sur~$L_{x,r^+}$ mais non sur~$L_{x,r}$.
La profondeur de~$\phi$ ne peut excéder~$r$, donc $r>0$.
Nous ajouterons bientôt des conditions de «\,{}généricité\,{}» sur~$\phi$ qui impliquent que sa profondeur est précisément~$r$. 

Par ailleurs, on dispose d'un plongement de~$\Imm(L)$ dans~$\Imm(G)$,
qui est canonique à l'action près du groupe  $\Char^\ast(\Z_\L)_F\otimes_\Z \R$ des $L$-automorphismes de $\Imm(L)$.
Cela permet de voir~$x$ comme un point de~$\Imm(G)$ et de former les groupes~$G_{x,t}$, pour $t \geq 0$. Définissons
\[ K = L_{\reduit{x}} \,{} G_{x,\rhalf}.\]
C'est un sous-groupe de~$G$.
Si~$Z_{\L,\G}=Z(\L)/Z(\G)$ est anisotrope, c'est-à-dire si~$Z_{\L,\G}(F)$ est compact,
alors~$K$ est compact modulo le centre.
Faisons désormais cette hypothèse.  Nous allons construire, à partir de~$\rho$ et~$\phi$, une représentation irréductible~$\mu$ de~$K$.

On peut  étendre la représentation~$\rho$ en une représentation~$\rho^\sharp$ de~$K$ en faisant agir trivialement les éléments de~$G_{x, \rhalf}$,
car  l'intersection~$L_{\reduit{x}} \cap G_{x, s}$ est contenue dans $L_{x, s}$ pour tout $s \geq 0$. 
Pour définir~$\mu$,  on construit une représentation lisse~$\kappa$ de~$K$ à partir du caractère~$\phi$, puis on pose \mbox{$\mu=\rho^\sharp \otimes \kappa$}.
Reste à décrire cette représentation~$\kappa$. 


\paragraph*{$\bullet$ Construction de $\kappa$ : l'étape abélienne}
Voyons d'abord comment le caractère~$\phi$ de~$L_{\reduit{x}}$ peut être étendu au groupe $K_+ =L_{\reduit{x}}\,{} G_{x,\rhalf+}.$ 

À partir des isomorphismes de Moy--Prasad évoqués au \no\ref{sec:la_filtration},
il n'est pas  difficile de construire un morphisme $\Theta\colon G_{x,\rhalf^+}/G_{x,r^+} \to L_{x,\rhalf^+}/L_{x,r^+}$
dont la restriction à~$L_{x,\rhalf^+}$ est induite par la projection de~$L_{x,\rhalf^+}$ sur $L_{x,\rhalf^+}/L_{x,r^+}$. 
\begin{smallrema}
Plus précisément, on introduit un supplémentaire explicite~$\lie{r}$ de $\lie{l}=\mathrm{Lie}(\L)(F)$ dans~$\lie{g}$
(\emph{grosso modo} une somme d'espaces radiciels ; \cf \cite[\S\,{}4.2]{Fintzen_Harvard}).
La projection $\lie{g}\to \lie{l}$ parallèlement à~$\lie{r}$ induit un morphisme \mbox{$\lie{g}_{x,\rhalf^+}/\lie{g}_{x,r^+}\to \lie{l}_{x,\rhalf^+}/\lie{l}_{x,r^+}$} ;
 les isomorphismes de Moy--Prasad fournissent alors le morphisme~$\Theta$.
\end{smallrema}
Par ailleurs, on sait que~$\phi$ induit un caractère de $L_{x,\rhalf^+}/L_{x,r^+}$ (par définition de~$r$).
En composant avec le morphisme $\Theta$, on obtient un caractère $\tilde{\phi}$ de~$G_{x, \rhalf^+}$, trivial sur~$G_{x, r^+}$.
Sur l'intersection~$G_{x, \rhalf^+} \cap  L_{\reduit{x}} = L_{x, \rhalf^+}$, le caractère~$\tilde{\phi}$ coïncide avec~$\phi$. 

Il existe donc un unique caractère~$\widehat{\phi}$ de~$K_+ =   L_{\reduit{x}}\,{}G_{x,\rhalf^+}$
qui coïncide avec~$\tilde{\phi}$ sur~$G_{x, \rhalf^+}$ et avec~$\phi$ sur~$L_{\reduit{x}}$.
Le caractère~$\widehat{\phi}$ est l'extension de~$\phi$ à~$K_+$ que nous voulions construire. 


\paragraph*{$\bullet$ Construction de $\kappa$ : entrée du groupe de Heisenberg}
Reste à construire, à partir du caractère~$\widehat{\phi}$ de~$K_+$, une représentation lisse~$\kappa$ de~$K=L_{\reduit{x}}G_{x, \rhalf}$. 

Commençons par remarquer que le groupe dérivé de~$G_{x, \rhalf}$ est contenu dans~$G_{x, \rhalf^+}$ d'après la propriété~\ass{a} du \no\ref{sec:la_filtration}.
Le caractère $\widehat{\phi}$ permet donc de définir une application $\omega\colon  G_{x,\rhalf} \times G_{x, \rhalf}\to \C$ par~$\omega(x,y) = \widehat{\phi}(xyx^{-1}y^{-1})$.
Posons~$J = G_{x, \rhalf}$ et~\mbox{$J_{+} = L_{x, \rhalf}G_{x, \rhalf^+}$}.
Le quotient 
\[ W = J/ J_{+}\]
peut donc être muni d'une structure d'espace vectoriel sur~$k_F$ ou, si l'on préfère, sur le sous-corps premier~$\kprem$ de~$k_F$.
Ce dernier est isomorphe à~$\F_p$ ; si on fixe une identification entre le groupe additif de~$\kprem$ et le groupe des racines $p$\iemes de l'unité dans~$\C$,
alors~$\omega$ peut être vue comme une application de~$W \times W$ dans~$\kprem$ ;
on constate alors que~$\omega$ est une \emph{forme symplectique} sur~$W$.
Dans la suite de ce numéro, nous supposerons que le caractère~$\phi$ est tel que la forme~$\omega$ soit \emph{non dégénérée}.

C'est par la forme symplectique~$\omega$ qu'arrive le \emph{deus ex machina} des constructions connues de représentations supercuspidales :
le couple formé du groupe de Heisenberg et de la représentation de Weil du groupe symplectique.
Dans le cas des corps finis, ce couple semble naître dans les travaux de~\textcite{Gerardin_1975, Gerardin_1977} sur notre sujet. 

Soit~$\mathscr{H}(W)$ le groupe de Heisenberg de~$W$,
d'ensemble sous-jacent $W \times \kprem$ et de loi $(u_1, \alpha_1)(u_2,\alpha_2)=(u_1+u_2,\alpha_1+\alpha_2+ \frac12\omega(\alpha_1,\alpha_2))$.
C'est un $p$-groupe. Soit~$N$ le sous-groupe~$\ker(\widehat{\phi}) \cap J_+$ de~$J$. Alors le centre de~$J/N$ est~$J_+/N$. De plus, il existe un sous-groupe~$H$ de~$J/N$ et un isomorphisme~$\iota$ de~$\mathscr{H}(W)$ sur~$H$, tous deux explicites, tels qu'on ait \mbox{$J/N = H\cdot(J_+/N)$} et que~$\iota$ envoie le centre $\{1_H\} \times \kprem$ de~$\mathscr{H}(W)$ sur l'intersection $Z=H \cap (J_+/N)$, qui est donc isomorphe à~$\kprem$. Le caractère~$\widehat{\phi}$ induit un caractère non trivial~$\psi$ de~$J_{+}/N$, donc un caractère non trivial~$\psi^\flat$ du centre de~$\mathscr{H}(W)$.

Or il existe (à équivalence près) une unique représentation irréductible de~$\mathscr{H}(W)$ de caractère central~$\psi^\flat$ («  théorème de Stone--von Neumann »). On en déduit une représentation irréductible de~$H$ sur un espace~$\EspHeis$.
En faisant agir~$J_+/N$ sur~$\EspHeis$ par le caractère~$\psi$, on obtient alors une représentation de~$J/N$, puis de~$J$, sur~$\EspHeis$. 


\paragraph*{$\bullet$ Construction de $\kappa$ : conclusion avec la représentation de Weil}
À ce stade, on dispose d'une représentation de~$J=G_{x, \rhalf}$ sur~$\EspHeis$ et du caractère~$\phi$ de~$K^0=L_{\reduit{x}}$.
Nous allons les combiner en une représentation de $K=K^0 J$ sur~$\EspHeis$. 

Observons pour cela que~$K^0$ agit sur~$J$ par conjugaison en préservant~$J_+$ ;
on en déduit un morphisme de~$K^0$ dans~$\mathrm{GL}(W)$, qui est en fait à valeurs dans le groupe symplectique~$\mathrm{Sp}(W)$.
Or, la représentation de Weil fournit une action de~$\mathrm{Sp}(W)$ sur~$\EspHeis$
qui s'accorde si bien avec l'action de~$\mathscr{H}(W)$ qu'elle induit une représentation irréductible sur~$\EspHeis$
du produit semi-direct (externe) $ \mathscr{H}(W) \rtimes \mathrm{Sp}(W)$.
Les morphismes de~$H$ dans~$\mathscr{H}(W)$ et de~$K^0$ dans~$\mathrm{Sp}(W)$ induisent alors une représentation irréductible sur~$\EspHeis$
du produit semi-direct externe $J \rtimes K^0$ associé à l'action de~$K^0$ sur~$J$ par conjugaison.
Cela définit un morphisme de $J \rtimes K^0$ dans~$\mathrm{GL}(\EspHeis)$, que l'on peut tordre par le caractère~$\phi_{|{K^0}}$ de~$K^0$. 
On obtient ainsi une représentation irréductible de~$J \rtimes K^0$ sur~$\EspHeis$, qui se factorise par l'application évidente de~$J \rtimes K^0$ dans~$K=K^0J$.
Elle définit donc une représentation irréductible~$\kappa$ de $K$ sur~$\EspHeis$, ce qui achève la construction.

Sous l'hypothèse que~$Z_{\L,\G}$ est anisotrope et que~$\phi$ est tel que $\omega$ soit non dégénérée, on~a donc construit un couple~$(K, \mu)$ comme annoncé. Si le caractère~$\phi$ satisfait de plus une condition de \emph{généricité} (\cf \cite[\S\,{}2.1]{Fintzen_Michigan}), alors la représentation~$\cInd_K^G(\mu)$ est irréductible, donc supercuspidale. Cela conclut la description de notre cas particulier.


\subsubsection{Données cuspidales arbitraires}\label{sec:donnees_cuspidales}
Passons au cas général. Supposons donnés :
\begin{enumerate}
\item[(1)] Une suite croissante $L^0 \subset \cdots \subset L^n=G$ de sous-groupes de Levi tordus de~$G$ ;
\item[(2a)] Un point~$x$ de~$\Imm(L^0)$ ;
\item[(2b)] Une représentation lisse~$\rho$ de~$L^0_{\reduit{x}}$ ;
\item[(3a)] Des nombres réels $r_0, \dots r_{n-1}$ avec $0<r_0 < \cdots < r_{n-1}$ ;
\item[(3b)] Pour tout~$i \in \{0, \dots, n-1\}$, un caractère~$\phi_i$ de~$L^i$.
\end{enumerate}
On dit qu'ils définissent une \emph{donnée cuspidale}~$\mathscr{D}$ pour~$G$ lorsqu'ils vérifient :
\begin{enumerate}
\item[(H1)]  Le groupe $Z(\L^0)/Z(\G)$ est anisotrope ; 
\item[(H2)] Le point~$\reduit{x}$ est un sommet de~$\Immred(L^0)$, la représentation~$\rho$ est triviale sur~$L^0_{x,0^+}$,
et la restriction~$\rho_{|L^{0}_{x, 0}}$ induit une représentation du groupe réductif $L^{0}_{x, 0}/L^{0}_{x, 0^+}$ qui contient une somme finie de représentations irréductibles cuspidales ;
\item[(H3)] Pour tout~$i \in \{0, \dots, n-1\}$, le caractère~$\phi_i$ est de profondeur~$r_i$
et satisfait l'hypothèse de généricité de~\textcite[\S\,{}2.1]{Fintzen_Michigan} pour le couple $(L^i, L^{i+1})$.
\end{enumerate}
Si la représentation~$\rho$ de (2b) est \emph{irréductible}, nous dirons que la donnée~$\mathscr{D}$ est irréductible.
Si~$n=0$, il n'y a pas d'objets (3a)--(3b) et on retrouve les données qui permettent de construire les représentations de profondeur nulle (\no\ref{sec:profondeur_nulle}).

\begin{smallrema}
Remarquons que se donner un couple~$(x,\rho)$ comme en~(2a-b) satisfaisant~(H2)
revient à se donner une représentation supercuspidale de profondeur nulle~$\sigma^0$ de~$L^0$, qu'on appelle le \emph{socle} de la donnée~$\mathscr{D}$.
Par ailleurs, sous l'hypothèse~(H3), les caractères~$\phi_i$ déterminent les~$r_i$.
Fournir une donnée cuspidale~$\mathscr{D}$ revient donc à se donner
une suite croissante $L^0 \subset \cdots \subset L^n=G$ de sous-groupes de Levi tordus,
une représentation supercuspidale de profondeur nulle de~$L^0$,
et des caractères des~$L^i$ dont la profondeur va strictement croissant,
vérifiant (H1)--(H3). Si~$n=1$, on retrouve les données du~\no\ref{sec:Yu_1pas}.

\end{smallrema} 

Expliquons comment construire une représentation supercuspidale~$\sigma(\mathscr{D})$ à partir d'une donnée cuspidale~$\mathscr{D}$
(\cf~\cite{Fintzen_Compositio} et~\cite{Aubert_Tits} ; la présentation ci-dessous est proche de celle d'Aubert, tandis que celle de Fintzen est «\,{}dé-récursifiée\,{}»).
Posons
\begin{equation}\label{def_groupe_compact} K = (L^0)_{\reduit{x}} (L^1)_{x, \frac{r_{0}}{2}}\cdots (L^n)_{x, \frac{r_{n-1}}{2}} \end{equation}
en utilisant le fait que~$L^0$ est un sous-groupe de Levi tordu non seulement de~$G$, mais de tous les~$L^i$,
ce qui fournit des plongements de~$\Imm(L^0)$ dans chacun des immeubles~$\Imm(L^i)$ et permet de définir chacun des facteurs.
On vérifie que~$K$ est un sous-groupe de~$G$, compact modulo le centre grâce à l'hypothèse~(H1),
et que la représentation~$\rho$ s'étend, comme dans le \no\,{}précédent, en une représentation irréductible~$\rho^\sharp$ de~$K$.
Nous allons construire une représentation irréductible lisse~$\kappa$ de~$K$, puis poser 
\begin{equation}\label{def_sigma_D} \sigma(\mathscr{D})=\cInd_{K}^G(\rho^\sharp \otimes \kappa).\end{equation}
La construction de~$\kappa$ est récursive et utilise le \no\ref{sec:Yu_1pas}.
Pour~$i \in \{0, \dots, n\}$, notons
\[ K^i = (L^0)_{\reduit{x}} (L^1)_{x, \frac{r_{0}}{2}}\cdots (L^i)_{x, \frac{r_{i-1}}{2}} 
\qquad \text{et} \qquad 
K^i_+ = (L^0)_{\reduit{x}} (L^1)_{x, \frac{r_{0}}{2}^+}\cdots (L^i)_{x, \frac{r_{i-1}}{2}^+}.\] 
Posons~$\kappa_0=\phi_0$ ; c'est une représentation  de dimension~$1$ de~$K^0$.
Supposons donnée une représentation~$\kappa_{i-1}$ de~$K^{i-1}$.
On étend d'abord~$\phi_i$ en un caractère~$\widehat{\phi_i}$ de~$K^{i}_+$ en suivant l'étape abélienne du \no\,{}précédent.
Pour étendre~$\widehat{\phi_i}$ à~$K^{i}$, on reprend la méthode ci-dessus.
On peut définir des sous-groupes $J^i$ et~$J^{i}_+$ de telle façon que le quotient s'identifie à un~$\kprem$-espace vectoriel
sur lequel le caractère~$\widehat{\phi_i}$ induit une forme symplectique
(voir \cite[p.~9]{Aubert_Tits}, pour la définition de~$J^i$ et~$J^{i}_+$ ; le groupe~$J^i$ est l'un des groupes notés~$(G_i)_{x, \tilde{r}, \tilde{r}'}$ par \cite[\S\,{}2.5]{Fintzen_Compositio}).
Sous l'hypothèse~(H3), cette forme symplectique est non dégénérée.
On peut alors reprendre la construction du \no\ref{sec:Yu_1pas} en~y remplaçant~$(K, K^0, \phi, J, J_+)$ par $(K^{i}, K^{i-1}, \kappa_{i-1}, J^{i}, J^{i}_+)$.
Cela fournit une représentation~$\kappa_i$ de~$K^i$ sur un espace~$V_i$.
On a $K^n = K$ et la représentation~$\kappa=\kappa_n$ de~$K$ est irréductible~\parencite[\S\,{}3]{Fintzen_Compositio}.
C'est elle que l'on insère dans~\eqref{def_sigma_D} pour définir~$\sigma(\mathscr{D})$. 
\begin{theo}[\cite{Fintzen_Compositio}]\label{th:construction_sans_twist}
Si la donnée~$\mathscr{D}$ est irréductible,
alors la représentation~$\sigma(\mathscr{D})$ de~$G$ est  irréductible et supercuspidale.
\end{theo}
Si~$\mathscr{D}$ est une donnée cuspidale arbitraire, la représentation~$\sigma(\mathscr{D})$ est donc somme directe finie d'irréductibles supercuspidales. Dans tous les cas, sa profondeur est~$r_{n-1}$. 


\subsubsection{Le caractère correcteur}\label{sec:twist}
Le th.~\ref{th:construction_sans_twist} était le résultat principal de~\textcite{Yu_2001}
et la démonstration de Yu s'appuyait sur les propriétés générales de l'induction compacte évoquées ci-dessus (p.~\pageref{critere_irred}, a) et b)),
en vérifiant des résultats de non entrelacement. 

Malheureusement, Yu utilisait un résultat de~\textcite{Gerardin_1977} dont l'énoncé imprimé contenait une coquille
(c'est l'impression qui est en cause, et non la preuve de Gérardin). 

Heureusement, \textcite{Fintzen_Compositio} a trouvé une autre preuve du  th.~\ref{th:construction_sans_twist},
qui utilise le résultat de Gérardin corrigé et un argument différent de celui de Yu (et plus court).

Malheureusement, les résultats intermédiaires au cœur de la preuve de Yu sont bel et bien faux (\cf~\cite{Fintzen_Compositio}, \S\,{}4) et ces résultats avaient été largement utilisés,
notamment dans la recherche de formules du caractère pour les représentations~$\sigma(\mathscr{D})$~\parencite{Adler_Spice, Spice_Compositio, Spice_2}.
C'est ce qui a permis de repérer le problème en 2018.

Heureusement,~\textcite{Fintzen_Kaletha_Spice} ont montré qu'une modification subtile mais simple de la construction de~Yu permet de rétablir l'ensemble de ces résultats.
Cette modification induit d'ailleurs des simplifications non négligeables dans les formules de caractère.
C'est un ingrédient important des travaux discutés aux \S\,{}\ref{sec:supercuspidales_nonsingulieres}--\ref{sec:LLC}. 

Fixons une donnée cuspidale~$\mathscr{D}$ et reprenons les notations du \no\,\ref{sec:donnees_cuspidales}. La modification salvatrice consiste à définir un caractère~$\varepsilon\colon K\to \{\pm1\}$ du groupe~\eqref{def_groupe_compact} et à poser 
\begin{equation}\label{def_pi_D} \pi(\mathscr{D})=\cInd_{K}^G(\rho^\sharp \otimes \kappa \otimes \varepsilon).\end{equation}

\begin{theo}[\cite{Fintzen_Kaletha_Spice}]\label{th:construction_avec_twist}
Si la donnée~$\mathscr{D}$ est irréductible, alors la représentation~$\pi(\mathscr{D})$ de~$G$ est  irréductible et supercuspidale.
\end{theo}

La définition du caractère quadratique~$\varepsilon$  \parencite[déf.~3.1]{Fintzen_Kaletha_Spice} est extrêmement technique
et combine des ingrédients de nature algébrique dont le rapport avec la construction de Yu n'est pas évident.
Il est  heureux qu'elle permette tout à la fois de «\,{}réparer\,{}» la construction de Yu et d'apporter des simplifications essentielles pour les applications qu'on en verra ci-dessous.
Je recommande avec enthousiasme un exposé de~\textcite{Kaletha_expose_3signs} qui explique très bien les ressorts du caractère~$\varepsilon$.


\subsection{Résultats de classification}


\subsubsection{}
Si~$\G$ est le groupe $\GL(n)$ et $p > n$, la construction de Yu (sans caractère correcteur) coïncide avec celle de~\textcite{Howe_1977}, qui fournit toutes les représentations supercuspidales irréductibles de $\GL(n,F)$~\parencite{Moy_1986}.
Lorsque~$F$ est de caractéristique~0, Kim~\parencite*{Kim_exhaustion} avait montré que la construction de Yu est exhaustive si~$p$ est «\,{}très grand\,{}».

\begin{theo}[\cite{Fintzen_Annals}] \label{th:fintzen_exhaustion}
Supposons que~$\G$ se déploie sur une extension modérément ramifiée et que~$p$ ne divise pas l'ordre du groupe de Weyl absolu.
Pour toute représentation irréductible supercuspidale~$\pi$ de~$G$, il existe une donnée cuspidale irréductible~$\mathscr{D}$
telle que~$\pi$ soit équivalente à la représentation~$\pi(\mathscr{D})$.
\footnote{Fintzen prouve le résultat pour la construction $\mathscr{D}\mapsto \sigma(\mathscr{D})$, c'est-à-dire pour celle de Yu sans la modification du \no\ref{sec:twist}.
Mais cela implique le résultat pour la construction « tordue »~\mbox{$\mathscr{D}\mapsto \pi(\mathscr{D})$} :
si~$\mathscr{D}$ est fixée, on peut écrire $\sigma(\mathscr{D}) \simeq  \pi(\mathscr{D}')$ pour une donnée~$\mathscr{D}'$ adéquate.}
\footnote{La construction de Yu et le th.~\ref{th:fintzen_exhaustion} s'étendent au-delà des représentations complexes et s'appliquent aux représentations lisses « modulaires » de~$G$ sur des $\mathscr{C}$-espaces vectoriels, dès que~$\mathscr{C}$ est un corps de caractéristique $\neq p$  (\cf \cite{Fintzen_Michigan}). De plus, \parencite{Fintzen_Annals} ne se limite pas au cas supercuspidal et montre l'existence de « types », au sens de Bushnell et Kutzko, pour toutes les représentations~lisses.  }
\end{theo}

La méthode de Fintzen est, pour l'essentiel, constructive. Elle est aussi entièrement algébrique, là où celle de Kim utilisait la formule de Plancherel. 

\begin{smallrema}
Voici de brèves indications sur la preuve de Fintzen.
Plutôt que de reconstruire directement une donnée cuspidale~$\mathscr{D}$ à partir d'une représentation~$\pi$, elle introduit une notion de « squelette » pour~$\pi$,
initialement plus faible que celle de donnée cuspidale :
c'est un uplet $\mathscr{S}=((L^i), x, \rho, (r_i), (X_i))$ où, \emph{grosso modo}, les objets $L^i, x, \rho, r_i$ sont comme en (1)--(3a),
mais les $X_i$ sont des formes linéaires sur~$\lie{g}$ dont les centralisateurs sont les $L^i$.
De plus, l'hypothèse (H2) est considérablement affaiblie (on ne demande pas que~$x$ soit un sommet et on n'impose pas à $\rho$ de condition de cuspidalité),
et l'hypothèse (H3) est remplacée par une condition analogue sur les~$X_i$.
Si~$\mathscr{S}$ est relié à la représentation~$\pi$ par une condition qui apparaît chez Yu, et s'il vérifie une certaine condition de « maximalité »,
Fintzen décrit un algorithme qui construit à partir de~$\mathscr{S}$ une donnée~$\mathscr{D}$ vérifiant~$\pi \simeq \sigma(\mathscr{D})$.
Reste à prouver l'existence d'un squelette ayant  « la bonne » relation à la représentation~$\pi$.
C'est le cœur de l'argument, qui montre comment construire récursivement un squelette adéquat en partant de~$\pi$ et d'un point arbitraire de~$\Imm(G)$,
en passant par une ingénieuse descente~infinie.
\end{smallrema}


\subsubsection{}
D'après les th.~\ref{th:construction_avec_twist} et~\ref{th:fintzen_exhaustion}, dans le cas bien modéré,
l'application~$\mathscr{D} \mapsto \pi(\mathscr{D})$ induit une surjection de l'ensemble des données cuspidales irréductibles
sur le dual supercuspidal~$\IrrCusp(G)$.
Pour avoir une véritable classification de~$\IrrCusp(G)$, il reste à savoir à quelle condition deux données cuspidales mènent au même élément de~$\IrrCusp(G)$.
Cette condition a été décrite par~\textcite{Hakim_Murnaghan},
sous une hypothèse dont Kaletha a récemment montré qu'elle n'est pas nécessaire \parencite*[\S\,{}3.5]{Kaletha_regular_supercuspidals}.

Dans ce numéro, on suppose que~$\G$ se déploie sur une extension modérée et que~$p \neq 2$.

Soit~$\mathscr{D}$ une donnée cuspidale.
Reprenons les notations du \no\ref{sec:donnees_cuspidales}
et écrivons $\mathscr{D}$ comme une famille $((L^0, \dots, L^n), x, \rho, (r_0, \dots, r_{n-1}), (\phi_0, \dots, \phi_{n-1}))$,
où j'espère transparente la nature des objets.
Nous dirons qu'une donnée~$\mathscr{D}'$ s'obtient à partir de~$\mathscr{D}$  :
\begin{itemize}
\item par \emph{transformation triviale du socle} si~$\mathscr{D'}$ est de la forme 
\[ ((L^0, \dots, L^n), x', \rho', (r_0, \dots, r_{n-1}), (\phi_0, \dots, \phi_{n-1}))\] 
où~$x' \in \Imm(L^0)$ vérifie~$\reduit{x'}=\reduit{x}$ et la représentation~$\rho'$ de~$(L^0)_{\reduit{x'}}=(L^0)_{\reduit{x}}$ est équivalente à~$\rho$ ;  
\item par \emph{conjugaison} s'il existe un élément~$g$ de~$G$ tel que~$\mathscr{D}'$ soit de la forme 
\[ ((gL^0g^{-1}, \dots, gL^ng^{-1}), g \cdot x, \,{}^{g}\rho, (r_0, \dots, r_{n-1}), (\,{}^{g}\phi_0, \dots, \,{}^{g}\phi_{n-1}))\pv\] 
\item par \emph{refactorisation} 
si~$\mathscr{D'}$ est de la forme 
\[ ((L^0, \dots, L^n), x, \rho', (r_0, \dots, r_{n-1}), (\phi'_1, \dots, \phi'_{n-1}))\] 
et si les représentations~$\rho, \rho'$ et caractères $\phi_i$, $\phi_i'$ vérifient les conditions suivantes :
\begin{conditions}
\item les restrictions à~$L^0_{\reduit{x}}$ de $(\rho \otimes \phi_0 \otimes \cdots \otimes \phi_{n-1})$ et  $(\rho' \otimes \phi'_0 \otimes \cdots \otimes \phi'_{n-1})$ sont~identiques ;
\item posons~$r_{-1}=0$ ; alors pour tout~$i \in \{0, \dots, n-1\}$, les restrictions à $(L^i)_{x, r_{i-1}^+}$ de~$\phi_i \cdots \phi_{n-1}$ et $\phi'_i \cdots \phi'_{n-1}$ sont identiques.
\end{conditions}
\end{itemize}

Ces opérations  sur l'ensemble des données cuspidales sont réversibles et commutent. Si~$\mathscr{D}_1, \mathscr{D}_2$ sont des données cuspidales, écrivons~$\mathscr{D}_1 \sim \mathscr{D}_2$
lorsqu'on peut obtenir~$\mathscr{D}_2$ à partir de~$\mathscr{D}_1$ par une suite finie de transformations triviales du socle, de conjugaisons et de refactorisations.
Cela définit  une relation d'équivalence entre données~cuspidales. 

\begin{theo}[\cite{Hakim_Murnaghan,Kaletha_regular_supercuspidals}] \label{th:hakim_murnaghan}\noindent Soient~$\mathscr{D}$, $\mathscr{D}'$ des données cuspidales. On a~$\pi(\mathscr{D})\simeq \pi(\mathscr{D}')$ si et seulement si~$\mathscr{D}\sim \mathscr{D}'$. 
\end{theo}


\section{Représentations supercuspidales régulières et non~singulières}\label{sec:supercuspidales_nonsingulieres}

On fixe toujours un corps local~$F$ et on note toujours~$G=\G(F)$ où~$\G$ est un groupe réductif connexe défini sur~$F$. 
Pour le moment, autorisons~$F=\R$ ou~$F$ non archimédien.

Si~$\bS$ est un tore maximal de~$\G$ qui est \emph{$F$-elliptique} (\no\ref{sec:def_elliptique}), et si~$S = \bS(F)$,
alors Kaletha construit un revêtement double~$S_{\pm}$ de~$S$, puis décrit  des fonctions spéciales sur~$\Spm$.
Ces fonctions coïncident avec celles du \S\,{}\ref{sec:cas_reel} dans le cas $F=\R$
et fournissent des formules de caractère longtemps attendues dans le cas non-archimédien.


\subsection{Le revêtement double de Kaletha et les fonctions spéciales associées}\label{sec:revetement_double}

La construction de Kaletha repose sur des ingrédients de nature galoisienne bien connus dans la théorie de l'endoscopie.
Décrivons ces ingrédients.
Dans toute la suite, on fixe une clôture séparable~$\Fsep$ de~$F$ et on note~$\Gamma$ le groupe de Galois absolu $\Gal(\Fsep\,{}|\,{}F)$.


\subsubsection{} \label{sec:fonctions_signe}
Soit~$\bS$ un tore défini sur~$F$.
Nous pensons au cas où~$\bS$ est un tore maximal de~$\G$ ; mais en vue des applications à la correspondance de Langlands au \S\,{}\ref{sec:LLC},
supposons plutôt qu'on dispose d'un plongement~$j\colon \bS \to \G$ qui est défini sur~$F$ et dont l'image est un tore maximal de~$\G$.
À partir de l'ensemble des racines de~$j(\bS)$ dans~$\G$,
on obtient par transport de structure une partie finie~$\Racines$ du réseau~$\Char^\ast(\bS)$ des caractères algébriques de~$\bS$ ;
par abus de langage, nous appelons \emph{racines} les éléments de~$\Racines$.
Si~$\alpha$ est un élément de~$\Racines$, alors~$(-\alpha)$ appartient aussi à~$\Racines$. 
L'action naturelle de~$\Gamma$ sur~$\Char^\ast(\bS)$ préserve~$\Racines$. 
 
On dit qu'une racine~$\alpha$ est \emph{symétrique} si l'orbite de~$\alpha$ sous~$\Gamma$ contient~$(-\alpha)$.
Notons~$\RacinesSym$ la partie de~$\Racines$ formée des racines symétriques ;
lorsque~$F=\R$ et~$\Fsep=\C$, on parle souvent de \emph{racines imaginaires}.

Si~$\alpha$ est une racine symétrique, alors le stabilisateur~$\Gamma_\alpha$ de~$\alpha$ est d'indice~$2$ dans le stabilisateur~$\Gamma_{\pm\alpha}$ de~$\{\pm\alpha\}$. Considérons les sous-corps~$F_\alpha, \Fpmalpha$ des points fixes de~$\Gamma_\alpha$ et~$\Gamma_{\pm\alpha}$ dans~$\Fsep$ ;
alors~$\Falpha$ est une extension quadratique de~$F_{\pm\alpha}$. 
Soit~$\kappa_\alpha\colon \Fpmalpha^\times\to\{\pm1\}$ le caractère quadratique associé à cette extension :
si~$N_{\Falpha/\Fpmalpha}\colon \Falpha^\times \to\Fpmalpha^\times$ est l'application norme (\cf~\cite[\S\,{}8]{Bourbaki_A5}),
alors~$\kappa_\alpha$ est l'unique morphisme non trivial dont le noyau contient l'image de~$N_{\Falpha/\Fpmalpha}$.
Pour tous~$\alpha \in \RacinesSym$ et~$\sigma \in \Gamma$, on a~$\kappa_{\sigma(\alpha)} = \kappa_\alpha \circ \sigma^{-1}$.

Supposons donné, pour chaque racine symétrique~$\alpha$, un élément~$\delta_\alpha$ de~$\Falpha^\times$ ;
autrement dit, considérons un élément $(\delta_\alpha)_{\alpha \in \RacinesSym}$ du produit $\prod_{\alpha \in \RacinesSym} F_\alpha^\times$.
Nous dirons que c'est une \emph{famille équivariante d'éléments des~$F_\alpha$}
si  $\delta_{\sigma(\alpha)} = \sigma(\delta_\alpha)$ pour tous~$\sigma \in \Gamma$ et $\alpha \in \RacinesSym$.  

Si $(\eta_\alpha)_{\alpha \in \RacinesSym}$ est une telle famille équivariante et si  de plus  $\eta_\alpha \in \Fpmalpha^\times$ pour tout~$\alpha$,
alors lorsque~$\alpha$ parcourt une orbite de~$\Gamma$ dans~$\RacinesSym$, la valeur de~$\kappa_{\alpha}(\eta_{\alpha})$ reste constante.
On la note~$\kappa_{\Omega}((\eta_\alpha))$,
et on note $\kappa\crochets{(\eta_\alpha)_{\alpha \in \RacinesSym}}$ le produit des~$\kappa_{\Omega}((\eta_\alpha))$ quand~$\Omega$ parcourt~$\RacinesSym/\Gamma$.


\subsubsection{} \label{sec:definition_du_revetement}
Ces ingrédients en place, on définit un revêtement double~$\Spm$ de~$S=\bS(F)$.
Nous supposerons~$\RacinesSym$ non vide, sinon on prend le revêtement trivial $S \times \{\pm 1\}$.

Soit d'abord~$\Sigma_{\pm}$ l'ensemble des couples~$(s, (\delta_\alpha)_{\alpha \in \RacinesSym})$,
où  $s$ est un élément de~$S$ et~$(\delta_\alpha)_{\alpha \in \RacinesSym}$ est une famille équivariante d'éléments des~$\Falpha^\times$, vérifiant la condition : 
\begin{equation} \label{def_sigma} \text{pour tout~$\alpha \in \RacinesSym$, on a $\delta_{\alpha}/\delta_{-\alpha} = \alpha(s)$.}\end{equation}
Munissons $\Sigma_{\pm}$ de la relation d'équivalence~$\sim$ la moins fine telle qu'on ait
\[ (s, (\delta_\alpha)_{\alpha \in \RacinesSym}) \sim (s, (\eta_\alpha\delta_\alpha)_{\alpha \in \RacinesSym})\]
pour toute famille équivariante $(\eta_{\alpha})_{\alpha \in \RacinesSym}$
vérifiant~$\eta_\alpha \in \Fpmalpha$ pour tout~$\alpha \in \RacinesSym$ et telle que $\kappa\crochets{(\eta_\alpha)_{\alpha \in \RacinesSym}}=1$.
Posons enfin~$\Spm= \Sigma_{\pm}/\!\!\sim$ et munissons-le de la loi de groupe induite par le produit terme à terme. 

Alors la projection $(s, (\delta_\alpha)_{\alpha \in \RacinesSym}) \mapsto s$ induit un morphisme surjectif de~$\Spm$ dans~$S$.
La fibre au-dessus de l'identité de~$S$ est de cardinal~$2$,
engendrée par la classe de tout élément de la forme~$(1, (\eta_\alpha))$
où $(\eta_\alpha)$ est une famille équivariante d'éléments des~$\Fpmalpha^\times$ vérifiant $\kappa\crochets{(\eta_\alpha)}=-1$.
La projection $\Spm \to S$ est donc un revêtement à deux feuillets. 

Si $F=\R$ et si~$\bS$ est un tore maximal elliptique de~$\G$,
alors on dispose d'un isomorphisme canonique de~$\Spm$ sur le revêtement double de nature plus algébrique utilisé dans le \S\,{}\ref{sec:cas_reel} (\cf\cite[\S\,{}5.2]{Kaletha_Covers}).

Comme dans le cas réel, on s'intéresse aux caractères \emph{spécifiques} de~$\Spm$, c'est-à-dire aux caractères de~$\Spm$ qui ne sont pas constants sur les fibres du revêtement.
On note~$\Charspec(\Spm)$ l'ensemble des caractères spécifiques de~$\Spm$.


\subsubsection{} \label{sec:protoformules_caracteres}
Nous allons associer à chaque caractère spécifique de~$\Spm$ une fonction à valeurs complexes sur~$S$.
Les fonctions obtenues nous mettront sur la voie des caractères de «\,{}presque toutes\,{}»  les représentations supercuspidales.  

On fixe désormais un caractère non trivial~$\Lambda\colon F \to \C^\ast$ du groupe additif de~$F$. 

Si~$F$ est non archimédien, on suppose~$p \neq 2$. Ajoutons deux hypothèses qui simplifient ce numéro : prenons~$\Lambda$ trivial sur l'idéal premier~$\premier_F$ de l'anneau d'entiers~$\entiers_F$, mais pas sur~$\entiers_F$ ;  et  supposons que~$j(\bS)/\Z(\G)$ se déploie sur une extension modérée de~$F$.

Soit donc~$\vartheta$ un caractère spécifique de~$S_{\pm}$.
Associons-lui d'abord une fonction~\mbox{$a_{S, \Lambda, \vartheta}\colon S_{\pm}\to \{\pm 1\}$} en utilisant les constructions du \no\ref{sec:fonctions_signe}. Cette étape est assez technique, sauf si~$\RacinesSym$ est vide : alors~$\Spm=S\times\{\pm 1\}$ et on prend la seconde~projection.

Pour définir~$a_{S, \Lambda, \vartheta}$ dans le cas où~$\RacinesSym$ est non vide, on « linéarise » le caractère~$\vartheta$ le long de chaque racine ; on attache ainsi un élément~$a_\alpha$ de~$F_\alpha^\times$ à chaque~$\alpha \in \RacinesSym$. 
Identifions~$\vartheta^2$ à un caractère de~$S$ et fixons $\alpha$ dans~$\RacinesSym$.
Alors~$\alpha$ vient avec une coracine~$\alpha^\vee \in \Char_\ast(\bS)$, qui induit un morphisme de~$F_\alpha^\times$ dans~$\bS(F_\alpha)$.
Par ailleurs, la norme~$N_{F_\alpha/F}$ induit une application  de~$\bS(F_\alpha)$ dans~$S=\bS(F)$.
On peut donc identifier~\mbox{$\vartheta^2 \circ N_{F_\alpha/F} \circ \alpha^\vee$} à un caractère~$\lambda_{\vartheta, \alpha}$ de~$F_\alpha^\times$.

Si~$\lambda_{\vartheta, \alpha}=1$, alors on pose~$a_\alpha=1$. 
Supposons désormais~$\lambda_{\vartheta, \alpha} \neq 1$ et notons~$\mathrm{Tr}_{F_\alpha/F}$ l'application trace \mbox{de~$F_\alpha^\times$ dans~$F$}. L'élément~$a_\alpha$ se définit en comparant le caractère~$\lambda_{\vartheta,\alpha}$ au  caractère~$\lambda_{\Lambda, \alpha}=\Lambda \circ \mathrm{Tr}_{F_\alpha/F}$ de~$F_\alpha$ :

\begin{itemize}
\item Si~$F=\R$, on prend pour~$a_\alpha$ l'unique~$a \in \C^\times$ tel qu'on ait \mbox{$\lambda_{\vartheta,\alpha}(\exp(z))=\lambda_{\Lambda, \alpha}(2a z)$} pour tout~$z$ dans $F_\alpha=\C$. Explicitement, si on prend~\mbox{$\Lambda\colon x \mapsto \exp(2i\pi x)$} et si on voit~$d\vartheta$ comme un élément de~$\lie{s}_\C^\ast$ (la complexifiée de l'algèbre~$\Lie(\bS)$), puis la coracine~$\alpha^\vee$ comme un élément de~$\lie{s}_\C$, alors \mbox{$a_\alpha = \dual{d\vartheta}{\alpha^\vee}/(2i\pi)$}. Cette quantité apparaît implicitement dans la notion de régularité introduite au~\S\,{}\ref{sec:caractere_HC}. 

\item Si~$F$ n'est pas archimédien, notons~$r_{\alpha}$ la profondeur de $\lambda_{\vartheta, \alpha}$ ;
ce dernier induit donc un caractère non trivial~$\lambda_{\vartheta,\alpha}^\flat$ de $(F_\alpha^\times)_{r_{\alpha}}/(F_\alpha^\times)_{r_{\alpha}^+}$.
Si~$r_\alpha >0$, il existe un unique~$\overline{a}_\alpha$ dans~\mbox{$(F_{\alpha})_{-r_\alpha}/(F_{\alpha})_{-r_\alpha^+}$ }vérifiant
\mbox{$\lambda^\flat_{\vartheta,\alpha}(z+1) =  \lambda_{\Lambda, \alpha}(2\overline{a}_\alpha z)$}
pour tout~\mbox{$z\in (F_\alpha)_{r_{\alpha}}/(F_\alpha)_{r_{\alpha}^+}$} ; on a~$\overline{a}_\alpha \neq 0$ et on en fixe un représentant~$a_\alpha$ dans~$F_\alpha^\times$, choisi de trace nulle. Si~$r_\alpha=0$, on choisit pour~$a_\alpha$ un élément arbitraire de valuation nulle et de trace nulle dans~$F_\alpha^\times$.
\end{itemize}
On définit maintenant la fonction $a_{S, \Lambda, \vartheta}\colon S_{\pm}\to \{\pm 1\}$ promise, par la formule 
\[
a_{S, \Lambda, \vartheta}\left( \crochets{s, (\delta_\alpha)_{\alpha\in\RacinesSym}}\right) = \prod \limits_{\substack{(\Gamma\cdot \alpha) \in \RacinesSym/\Gamma\\ \alpha(s)\neq 1}} \kappa_\alpha\left(\frac{\delta_\alpha-\delta_{-\alpha}}{a_\alpha}\right) \prod \limits_{\substack{(\Gamma \cdot \alpha) \in \RacinesSym/\Gamma \\ \alpha(s)= 1}} \kappa_\alpha\left(\delta_\alpha\right)
\]
où les produits sont pris sur des orbites de~$\Gamma$ dans $\RacinesSym$. (Le choix d'un représentant~\mbox{$\alpha \in \RacinesSym$} pour chaque orbite n'affecte pas le résultat.)
 
La fonction~$a_{S, \Lambda, \vartheta}$ obtenue ne dépend pas des choix faits ci-dessus pour définir les~$a_\alpha$ dans le cas non archimédien.
Cette fonction vaut $(-1)$ sur l'élément non trivial de la fibre de $\Spm \to S$ au-dessus de~$1_S$.
Le produit $a_{S, \Lambda, \vartheta}\cdot \vartheta$ est une fonction de~$S_{\pm}$ dans~$\C$, dont on voit sans peine qu'elle est constante sur les fibres de la projection~$\Spm \to S$.
Si l'on fixe un nombre complexe $c_{G, S, \Lambda}$, on définit donc comme annoncé une fonction~$\Theta_\vartheta^\sharp\colon S\to \C$ en posant 
\[ \Theta_\vartheta^\sharp(s) = c_{G, S, \Lambda} \sum \limits_{w \in \mathrm{Norm}_G(S)/S} \crochets{a_{S, \Lambda, \vartheta}\cdot \vartheta}(w^{-1} \cdot s) \quad \text{pour $s \in S$.}\]
On choisira pour constante~$c_{G,S,\Lambda}$ la valeur \mbox{$e(\G) \epsilon_L(\Char^\ast(\T_\G)_\C-\Char^\ast(\bS)_\C,\Lambda)$},
où~$e(\G)$ est le signe de Kottwitz~\parencite*{Kottwitz_1983}
et~$\epsilon_L(\cdots\negthinspace{})$ est un facteur~$\eps$ local dans la normalisation définie par Langlands  \mbox{(Tate, \cite*{Tate_Corvallis}, dans  \cite*[vol.~2, p.~17, 3.6.1]{Corvallis})}. Alors la fonction~$\Theta_\vartheta^\sharp$ ne dépend pas du choix~de~$\Lambda$.
 

\subsubsection{} \label{sec:annonce_caractere_shallow}
Supposons maintenant que~$\bS$ soit un \emph{tore maximal elliptique} de~$\G$ et utilisons l'inclusion de~$\bS$ dans~$\G$ pour plongement.
Notons~$\Sreg$ l'ensemble des éléments réguliers de~$S$, égal à~\mbox{$S \setminus \Delta^{-1}(0)$} où~$\Delta$ est le discriminant de Weyl (\no\ref{sec:caractere_sur_reg}).
On définit une fonction $\Theta_\vartheta\colon \Sreg \to \C$ en posant 
\begin{equation} \label{def_theta} \Theta_\vartheta(s)= \Delta(s)^{-1/2}\,{} \Theta_\vartheta^\sharp(s) \qquad (s \in \Sreg)\pt\end{equation} 
L'observation suivante laisse peu de doute sur l'intérêt des fonctions~\eqref{def_theta} : 
\begin{prop}[\cite{Kaletha_regular_supercuspidals}, \S\,{}4.11]
Si $F=\R$ et si~$\vartheta$ est régulier \textup{(\S\,{}\ref{sec:caractere_HC})}, alors~$\Theta_\vartheta$ coïncide avec la fonction~\eqref{formule_caractere_2}.
C'est donc la restriction à $\Sreg$ du caractère d'une  représentation de série discrète de~$G$, uniquement déterminée à équivalence près.
\end{prop}

Si~$F$ est non archimédien et si on se place dans le cas bien modéré (\no\hyperref[sec:cas_bien_modere]{1.5}), les travaux de \textcite{Kaletha_regular_supercuspidals, Kaletha_Lpackets} établissent des liens très remarquables
entre la fonction~$\Theta_\vartheta$ et le caractère de représentations supercuspidales.
Pour « presque tout »~$\vartheta$,
la fonction~$\Theta_\vartheta$ coïncide avec le caractère d'une représentation supercuspidale --- non sur~$\Sreg$ tout entier
(nous verrons que le caractère  y prend nécessairement une forme plus compliquée),
mais sur les éléments de~$\Sreg$ qui sont \emph{topologiquement semi-simples modulo le centre de~$\G$}.

Expliquons brièvement cette notion.
Soit~$g$ un élément de~$G$.
On dit que~$g$ est \emph{topologiquement semi-simple modulo $\Z_\G=\Z(\G)$} si $g$ est semi-simple et si, pour tout tore maximal~$\T$ du centralisateur $\mathrm{Cent}_{\G}(g)$ et pour tout caractère algébrique $\beta$ de~$\T/\Z_\G$,
la valeur~$\beta(g)$ est d'ordre fini premier à~$p$.
On note~$\Gtopss$ l'ensemble des éléments topologiquement semi-simples modulo~$\Z_\G$ de~$G$, que Kaletha appelle  «\,{}superficiels\,{}». 
\smallskip

Sous certaines hypothèses mineures\footnote{Il serait souhaitable de disposer de conditions explicites sur le caractère~$\vartheta$ qui garantissent l'existence d'une telle construction. De telles conditions découlent implicitement des constructions qui vont suivre, mais à ma connaissance, une condition  simple formulée directement sur~$\vartheta$ reste à décrire.} sur~$\vartheta$,
Kaletha construit \emph{explicitement} une représentation lisse~$\pi(\vartheta)$ de~$G$ vérifiant : 
\begin{assertions}
\item Le caractère de~$\pi(\vartheta)$ est égal à $\Theta_\vartheta$ sur~$\Sreg \cap \Gtopss$ ;
\item En général,  $\pi(\vartheta)$ est une somme directe finie de supercuspidales irréductibles ;
\item Si~$\vartheta$ vérifie une hypothèse de régularité supplémentaire, alors~$\pi(\vartheta)$ est irréductible.
\end{assertions} 

L'analogie avec le cas réel serait presque parfaite si~\ass{a} suffisait à déterminer la classe d'équivalence de~$\pi(\vartheta)$, comme dans le th.~\ref{HC_parametrage}.
Nous verrons au  \no\ref{sec:chan_oi} que c'est parfois le cas.
Hélas, ce n'est pas la situation générale et $\Sreg \cap \Gtopss$ peut même être vide. 


\subsection{Définition et paramétrage  des représentations non singulières}\label{sec:def_reg_sc}

Les observations du \no\ref{sec:annonce_caractere_shallow} semblent indiquer que beaucoup de représentations supercuspidales peuvent, comme dans le cas réel,
être reliées à des couples $(\bS, \vartheta)$ où~$\bS$~est un tore maximal elliptique de~$\G$ et~$\vartheta$ est un caractère spécifique du revêtement double~$S_{\pm}$.
Mais la classification du \S\,{}\ref{sec:Yu_et_Fintzen} se fait par les données cuspidales du \no\ref{sec:donnees_cuspidales},
apparemment plus sophistiquées et où les caractères de tores n'apparaissent pas.

Cela dit, si~$\mathscr{D}$ est une donnée cuspidale, alors son socle est l'induite compacte (vers un sous-groupe de Levi tordu)
d'une représentation qui provient essentiellement d'une représentation cuspidale d'un groupe réductif défini sur le corps fini~$k_F$.
Or, \textcite{Deligne_Lusztig} ont construit beaucoup de ces représentations cuspidales à l'aide de caractères de tores (finis). 

Si~$\Gfini$ est un groupe réductif connexe défini sur~$k_F$, et si~$\Sfini$ est un tore maximal de~$\Gfini$,
alors Deligne et Lusztig associent à tout caractère~ $\xi$ de~$\Sfini(k_F)$ 
une représentation virtuelle~$\DeLus_{\Sfini, \xi}$ de~$\Gfini(k_F)$.
Si le tore~$\Sfini$ est \emph{elliptique}, c'est-à-dire si tout sous-tore de~$\Sfini$ qui se déploie sur~$k_F$ est contenu dans le centre de~$\Gfini$,
et si le caractère $\xi$ n'est orthogonal à aucune coracine de la paire~$(\Gfini, \Sfini)$,
alors~$\DeLus_{\Sfini, \xi}$ s'identifie (à un signe global près) à une authentique représentation de~$\Gfini(k_F)$,
qui est de plus somme finie d'irréductibles cuspidales.
La représentation~$\DeLus_{\Sfini, \xi}$ est irréductible dès que~$\xi$ est \og en position générale \fg,
ce qui signifie que le seul automorphisme de~$\Sfini$ qui fixe~$\xi$ et qui est induit par un automorphisme intérieur de~$\Gfini(k_F)$ est l'identité. 
Voir l'exposé de~\textcite{Serre_seminaire}.

C'est en utilisant la théorie de Deligne--Lusztig que~\textcite{Kaletha_regular_supercuspidals, Kaletha_Lpackets}, suivant une idée de~\textcite{Murnaghan_2011}, 
montre que «\,{}la plupart\,{}» des représentations supercuspidales peuvent se construire \emph{grosso modo} à partir d'un couple $(\bS, \vartheta)$ comme ci-dessus.

Jusqu'à la fin du \S\,{}\ref{sec:def_reg_sc}, on suppose que~$\G$ se déploie sur une extension modérée de~$F$ et que la caractéristique résiduelle~$p$ est~$\neq 2$.


\subsubsection{Supercuspidales régulières et non~singulières de profondeur nulle} \label{sec:reg_prof_nulle}
Soit~$\pi$ une représentation supercuspidale irréductible de~$G$ de profondeur nulle.
Compte tenu du \no\ref{sec:profondeur_nulle}, fixons un point~$x$ de~$\Imm(G)$ tel que~$\reduit{x}$ soit un sommet de~$\Immred(G)$
et une représentation irréductible~$\sigma$ de~$G_{\reduit{x}}$ vérifiant~$\pi \simeq \cInd_{G_{\reduit{x}}}^G(\sigma)$ et telle que~$\sigma_{|G_{x,0}}$ contienne une représentation cuspidale du groupe fini~$G_{x,0}/G_{x,0^+}$.
Identifions ce dernier à~$\Gfini_x(k_F)$ où~$\Gfini_x$ est un groupe réductif connexe défini sur~$k_F$. 

Soit~$\Sfini$ un tore maximal  elliptique de~$\Gfini_x$. 
Alors il existe un tore maximal elliptique~$S$ du groupe $p$-adique~$G$, minimalement ramifié\footnote{Si~$S=\bS(F)$, on dit que~$S$ est \emph{non ramifié} si~$\bS$ se déploie sur une extension non ramifiée de~$F$, et que~$\bS$ est \emph{minimalement ramifié} si la dimension du plus grand sous-tore non ramifié de~$\bS$ est la plus grande possible parmi les dimensions de sous-tores non ramifiés de~$\G$ ; \cf \textcite{Kaletha_regular_supercuspidals}, \S\,{}3.4.1. }
dans~$G$,
tel que~$\Sfini(k_F)$ s'identifie au quotient~$S_{0}/S_{0^+}$, où~$S_0, S_{0^+}$ sont les groupes issus de la filtration de Moy--Prasad pour l'unique point de~$\Immred(S)$.
Si~$\xi$ est un caractère de~$\Sfini(k_F)$, disons que~$\xi$ est \emph{régulier} (ou $G$-régulier) si,
lorsqu'on l'identifie à un caractère de~$S_0$ trivial sur~$S_0^+$, il n'existe pas d'automorphisme non trivial de~$S$
 soit induit par un automorphisme intérieur de~$G$ et fixe~$\xi$.\footnote{Un tel automorphisme laisse toujours~$S_0$ stable, ce qui donne un sens à la condition « fixer~$\xi$ ». }
Pour que~$\xi$ soit régulier, il est nécessaire qu'il soit en position générale, si bien que~$\DeLus_{\Sfini, \xi}$ est cuspidale irréductible.

On dit alors que~$\pi$ est \emph{régulière} si $\sigma_{|G_{x,0}}$ contient une représentation irréductible triviale sur~$G_{x,0^+}$
qui induit, sur le groupe fini~$\Gfini_x(k_F)$, une représentation de la forme~$\DeLus_{\Sfini, \xi}$
où~$\Sfini\subset \Gfini_x$ est un tore maximal elliptique et~$\xi$ est un caractère régulier de~$\Sfini(k_F)$.

On peut préciser grandement le lien entre les représentations supercuspidales régulières de profondeur nulle et les caractères de tores. 

Si~$\bS$ est un tore maximal elliptique de~$G$ et~$S_{0}/S_{0^+}=\Sfini(k_F)$ comme ci-dessus, la notion de caractère régulier se propage du tore fini~$\Sfini(k_F)$ au groupe~$S$ :
si~$\theta$ est un caractère de profondeur nulle de~$S$, disons que~$\theta$ est \emph{régulier}
si sa restriction à~$S_0$ (qui est triviale sur~$S_{0^+}$) s'identifie à un caractère régulier de~$\Sfini(k_F)$. 

Considérons maintenant un couple~\mbox{$(\bS, \theta)$}, où~$\bS$~est un tore maximal elliptique minimalement ramifié de~$\G$
et~$\theta$ est un caractère régulier de profondeur nulle de~$S$ ; on dira que~$(\bS, \theta)$ est un \emph{couple régulier de profondeur nulle}. 

À partir d'un tel couple,~\textcite[\S\,{}3.4.4]{Kaletha_regular_supercuspidals} construit une représentation supercuspidale régulière~$\pi^G_{(\bS,\theta)}$ de~$G$.
Voici une idée du chemin.  
Par construction de l'immeuble réduit~$\Immred(G)$, le tore~$\bS$ détermine un sommet~$y$ de~$\Immred(G)$. Soit~$x$ un point de l'immeuble étendu~$\Imm(G)$ vérifiant~$\reduit{x}=y$.
On peut identifier $S_0/S_{0^+}$ à~$\Sfini(k_F)$, pour un tore maximal elliptique~$\Sfini$ du groupe connexe~$\Gfini_x$  vérifiant~$G_{x,0}/G_{x,0^+} = \Gfini_x(k_F)$.
Le caractère~$\theta$ induit un caractère~$\xi$ de~$\Sfini(k_F)$, en position générale.
En adaptant la méthode de Deligne et Lusztig, Kaletha  construit  à partir de~$(\Sfini, \xi)$ une représentation lisse irréductible~$\kappa_{(\bS,\theta)}$ du stabilisateur~$G_{\reduit{x}}$ (cette étape est délicate !). 
Enfin, il pose~$\pi^G_{(\bS, \theta)} = \cInd_{G_{\reduit{x}}}^{G} \kappa_{(\bS,\theta)}$. 

\begin{prop}[\cite{Kaletha_regular_supercuspidals}, prop. 3.4.27]\label{prop:regulieres_prof_nulle} \begin{assertions}
\item Pour tout couple régulier de profondeur nulle $(\bS, \theta)$, la représentation $\pi^G_{(\bS, \theta)}$ est irréductible, supercuspidale, régulière et de profondeur nulle \pv
\item L'application $(\bS, \theta) \mapsto \pi^G_{(\bS, \theta)}$ induit une bijection de l'ensemble des classes de \mbox{$G$-conjugaison} de couples réguliers de profondeur nulle sur l'ensemble des représentations supercuspidales irréductibles régulières de profondeur nulle de~$G$\pt
\end{assertions}
\end{prop}

Si l'on part d'un couple~$(\bS, \theta)$ comme ci-dessus et si on ne suppose plus~$\theta$ régulier,
il arrive que la construction ci-dessus fonctionne et permette de construire une représentation de profondeur nulle~$\pi^G_{(\bS, \theta)}$,
qui n'est cependant pas irréductible en général. 

C'est notamment le cas lorsque~$\theta$ est \emph{non~singulier}, au sens suivant.
Si~$\bS'$ est le plus grand sous-tore non ramifié de~$\bS$ et si~$F'$ est une extension non ramifiée de~$F$ sur laquelle~$\bS'$ se déploie,
alors pour tout~$\alpha$ dans l'ensemble~$\Racines$ des racines de~$\bS$ dans~$\G$, la restriction~$\alpha_\mathrm{res}$ de~$\alpha$ à~$S'$
induit un caractère~$\theta \circ N_{F'/F} \circ \alpha_{\mathrm{res}}^\vee$ de $(F')^\times$.
On dit que~$\theta$ est \emph{non singulier} si ce caractère a une restriction non triviale à $\entiers_{F'}^\times$, quel que soit $\alpha \in \Racines$.
La validité de cette condition ne dépend pas du choix de~$F'$.

Si~$\theta$ est non~singulier, la construction évoquée ci-dessus fournit encore une représentation~$\kappa_{(\bS,\theta)}$ du stabilisateur~$G_{\reduit{x}}$,
qu'on peut induire en une représentation lisse~$ \pi^G_{(\bS, \theta)}$ de~$G$, de profondeur nulle.
En général, $\kappa_{(\bS, \theta)}$ n'est pas irréductible et~$\pi^G_{(\bS, \theta)}$ est somme directe finie de représentations supercuspidales irréductibles.

Si~$\pi$ est une représentation irréductible de~$G$ de profondeur nulle,
on dit que~$\pi$ est  \emph{non~singulière} si elle est facteur irréductible de l'une des représentations~$\pi^G_{(\bS, \theta)}$. Pour la représentation~$\pi^G_{(\bS, \theta)}$ elle-même, je parlerai de \emph{représentation standard non~singulière}.

\begin{smallrema}
En vue des applications à la correspondance de Langlands, Kaletha étudie avec grand soin les composantes irréductibles de~$\pi_{(\bS, \theta)}$ et leurs multiplicités.
Cela revient à étudier la décomposition de~$\kappa_{(\bS,\theta)}$ ;
mais cette étude est difficile,  car la construction de~$\kappa_{(\bS, \theta)}$ nécessite d'adapter la théorie de Deligne--Lusztig au groupe~$G_{\reduit{x}}/G_{x,0^+}$.
Or ce dernier,  contrairement à $G_{x,0}/G_{x,0^+}$, ne provient pas d'un groupe \emph{connexe} sur~$k_F$ ; Kaletha étend donc au-delà du cas connexe des résultats de décomposition de~Lusztig. 
 \end{smallrema}

\subsubsection{Supercuspidales régulières et non~singulières de profondeur arbitraire} \label{sec:def_non_singulieres}
Passons au cas général, en supposant désormais que~$p$ ne divise pas l'ordre du groupe de Weyl.
Rappelons que si~$\mathscr{D}$ est une donnée cuspidale pour~$G$, son socle est une représentation de profondeur nulle d'un sous-groupe de Levi tordu de~$G$.
Disons que~$\mathscr{D}$ est \emph{régulière} (\resp \emph{non singulière}) si son socle est régulier (\resp une représentation standard non singulière) au sens du \no précédent.
Soit~$\pi$ une représentation lisse irréductible de~$G$ ; on dit que~$\pi$ est \emph{régulière} (\resp \emph{non~singulière}) si elle est équivalente à~$\pi(\mathscr{D})$
pour une donnée cuspidale régulière~$\mathscr{D}$ (\resp à un facteur irréductible de~$\pi(\mathscr{D})$ pour une donnée non singulière~$\mathscr{D}$).
Si~$\pi$ est régulière, toute donnée~$\mathscr{D}'$ vérifiant~$\pi \simeq \pi(\mathscr{D}')$ est régulière d'après le th.~\ref{th:hakim_murnaghan}. 

À nouveau, les représentations supercuspidales régulières peuvent être paramétrées par certains caractères de tores maximaux de~$G$. Précisons ce paramétrage.
Soit~$\mathscr{D}$ une donnée  non~singulière.
Écrivons $\mathscr{D}=((L^0, \dots, L^n), x, \rho, (r_0, \dots, r_{n-1}), (\phi_0, \dots, \phi_{n-1}))$ en reprenant les notations du \no\ref{sec:donnees_cuspidales},
et~$L^i = \L^i(F)$ où~$\L^i$ est un sous-groupe de Levi de~$\G$.
D'après le \no précédent, il existe un tore maximal minimalement ramifié~$\bS \subset \L^0$ et un caractère $\phi_{-1}$ de~$S=\bS(F)$, non singulier et de profondeur nulle,
tels que~$\rho$ soit équivalente à $\pi^{L^0}_{(\bS, \phi_{-1})}$ si~$\mathscr{D}$ est régulière, et à l'une de ses composantes irréductibles si~$\mathscr{D}$ est non~singulière.
Dans ces conditions,~$\bS$ est un tore maximal de chacun des~$\L^i$, en particulier de~$\G$.
On peut donc considérer le caractère~$\theta=\phi_{-1} (\phi_0)_{|S} \cdots (\phi_{n-1})_{|S}$ de~$S$.

Si~$\mathscr{D}$ est régulière (\resp non~singulière),  le couple~$(\bS, \theta)$ satisfait : 
\begin{conditions}
\item $\bS$ est un tore maximal elliptique de~$\G$ ;
\item soit~$E$ une extension modérée de~$F$ sur laquelle~$\bS$ se déploie, soit $R_{0^+}$ le système de racines formé des~$\alpha \in \Racines$ vérifiant~$\theta(N_{E/F}(\alpha^\vee(E_{0^+}^\times)))=1$, et soit~$\G^0 \subset \G$ l'unique groupe réductif connexe dont~$\bS$ est un tore maximal et~$R_{0^+}$ le système de racines ; alors~$\bS$ est minimalement ramifié dans~$\G^0$ ;
\item le caractère de profondeur nulle~$\theta_{|G^0}$ de~$G^0$ est $G^0$-régulier (\resp non~singulier).
\end{conditions}
Dans~(ii), le système $R_{0^+}$ ne dépend pas du choix de l'extension~$E$, donc~$G^0$ non plus.
De plus, si~$\theta$ provient d'une donnée non~singulière~$\mathscr{D}$ comme ci-dessus, on a~$G^0 = L^0$.

 Si~$(\bS, \theta)$ est un couple vérifiant (i)--(iii), on dit qu'il est \emph{elliptique régulier} (\resp \emph{elliptique non~singulier}) : \cf~\textcite[déf.~3.7.5]{Kaletha_regular_supercuspidals}. Et si le tore maximal elliptique~$\bS$ est donné, on dit qu'un caractère~$\theta$ de~$S$ est \emph{régulier} (\resp \emph{non singulier}) si le couple~$(\bS, \theta)$ l'est.
Ces caractères fournissent un paramétrage simple des représentations régulières :

\begin{prop}[\cite{Kaletha_regular_supercuspidals}, prop. 3.7.8]\label{prop:regulieres_generales}
La construction ci-dessus induit une bijection de l'ensemble des classes d'équivalence de représentations supercuspidales régulières
sur l'ensemble des classes de $G$-conjugaison de couples elliptiques réguliers\pt
\end{prop}

Pour démontrer la prop.~\ref{prop:regulieres_generales}, Kaletha décrit un algorithme qui reconstruit une donnée~$\mathscr{D}$ à partir d'un couple $(\bS, \theta)$.
Les conditions (ii)-(iii) ci-dessus fournissent essentiellement le socle de~$\mathscr{D}$, et il faut reconstruire les groupes~$L^i$ et les caractères~$\phi_{i}$.
L'algorithme généralise les « factorisations de Howe » pour~$\GL(n)$~\parencite{Howe_1977}. 

L'algorithme de factorisation fonctionne aussi dans le cas non~singulier :
à partir d'un couple non~singulier~$(\bS, \theta)$, il construit une donnée non~singulière~$\mathscr{D}$,
de telle façon que la classe d'équivalence de~$\mathscr{D}$ ne dépende que de la classe de conjugaison de~$(\bS, \theta)$. 

Si~$(\bS, \theta)$ est un couple non~singulier, on note~$\pi_{(\bS, \theta)}$ la représentation $\pi(\mathscr{D})$ ainsi obtenue. 
Si~$\theta$ est régulier, elle est irréductible, et en général, c'est une somme finie de supercuspidales irréductibles (avec des multiplicités qui peuvent être~$>1$). 

Il est raisonnable de dire que « presque tout » le dual cuspidal de~$G$ est formé de représentations régulières.
Par exemple, si~$\G$ est le groupe $\GL(n)$ et si~$p$ ne divise pas~$n$, toute représentation supercuspidale irréductible est régulière.
Si~$G = \mathrm{SL}(2, F)$ et si~$p \neq 2$, il n'y a que quatre représentations non régulières ; elles sont non singulières.


\subsection{Formules du caractère}

Dans ce paragraphe, on continue à se placer dans le cas bien modéré (\no\hyperref[sec:cas_bien_modere]{1.5}).

\subsubsection{Formule simple sur le lieu superficiel} \label{def_chi_data}
Les travaux de~\textcite{Kaletha_regular_supercuspidals} puis~\textcite{Fintzen_Kaletha_Spice} donnent des formules du caractère particulièrement simples
pour les représentations supercuspidales régulières, ou pour les représentations supercuspidales (réductibles) attachées à des données non~singulières.
Voyons-en d'abord le cas le plus simple, qui concrétise les annonces du \no\ref{sec:annonce_caractere_shallow}. 

Pour faire le lien avec le \S\,{}\ref{sec:revetement_double}, expliquons comment passer d'une donnée cuspidale non singulière à un caractère de revêtement double de tore maximal elliptique.
Soient~$\mathscr{D}$ une donnée non singulière et~$(\bS, \theta)$ un couple elliptique non singulier vérifiant~\mbox{$\pi_{(\bS, \theta)}=\pi(\mathscr{D})$}, comme construit au \no\ref{sec:def_non_singulieres}.
Considérons le revêtement~$\Spm$ du \no\ref{sec:definition_du_revetement}.
Pour passer de~$\theta$ à un caractère spécifique de~$\Spm$, on utilise une \emph{donnée de transfert}\footnote{Les données de transfert sont les « $\chi$-data » de~\textcite{Langlands_Shelstad}.}.
Il s'agit d'une famille $(\chi_\alpha)_{\alpha \in \Racines}$ de caractères des~$F_\alpha^\times$, équivariante sous l'action de~$\Gamma$, vérifiant \mbox{$\chi_{-\alpha}=\chi_\alpha^{-1}$} pour tout~$\alpha \in \Racines$ et~$(\chi_\alpha)_{|\Fpmalpha^\times}=\kappa_\alpha$ pour tout~$\alpha \in \RacinesSym$.
En pratique, nous utiliserons des données $(\chi_\alpha)_{\alpha \in \Racines}$ vérifiant $\chi_\alpha=1$ pour~$\alpha \in \Racines\setminus\RacinesSym$, que j'appellerai~\emph{réduites}. 

À chaque donnée réduite, on peut associer un caractère spécifique~$\chi$ de~$S_\pm$, de la façon suivante.
Avec les notations du \S\,{}\ref{sec:revetement_double}, si~$\tilde{s}=(s, (\delta_\alpha)_{\alpha \in \RacinesSym})$ est un élément de~$\Sigma_\pm$, alors quand $\alpha$ parcourt une orbite~$\Omega$ de~$\Gamma$ dans~$\RacinesSym$, la valeur de~$\chi_\alpha(\delta_\alpha)$ reste constante.
On note~$\chi_{\Omega}(\tilde{s})$ la valeur commune, puis on pose \mbox{$\chi(\tilde{s})=\prod_{\Omega \in \RacinesSym/\Gamma} \chi_\Omega(\tilde{s})$}.
À cause des contraintes sur les~$\chi_\alpha$, cette formule induit un caractère spécifique de~$\Spm$, comme annoncé.
La multiplication par ce caractère induit une bijection de~$\Char(S)$ sur~$\Charspec(\Spm)$. 

Plusieurs choix de données de transfert (réduites) sont possibles.
Pour écrire les formules du caractère, il y a un « bon choix »  \parencite[\S\,{}4.2]{Fintzen_Kaletha_Spice}  qui dépend de la donnée~$\mathscr{D}$ et dont la description est assez technique.
Une fois ce choix fixé, on obtient à partir de~$\theta$ un caractère spécifique~$\vartheta$ de~$\Spm$.
Si~$\theta_1, \theta_2$ sont deux caractères distincts et non singuliers de~$S$, alors les caractères spécifiques~$\vartheta_1, \vartheta_2$ correspondants sont distincts.
Et on a le résultat suivant \parencite[prop.~4.3.2]{Fintzen_Kaletha_Spice} :

\begin{prop} \label{prop:formule_shallow}
Pour toute donnée non singulière~$\mathscr{D}$, le caractère de~$\pi(\mathscr{D})$ est donné par la formule~\eqref{def_theta}
sur l'ensemble~$\Sreg \cap \Gtopss$ des éléments superficiels de~$S$\pt
\end{prop}


\subsubsection{Formule fine dans le cas où l'exponentielle converge}\label{sec:formule_fine}
On conserve les notations du \no précédent. S'il y a peu (ou pas) d'éléments superficiels dans~$S$, la prop.~\ref{prop:formule_shallow} donne peu d'informations.
Mais on peut aller beaucoup plus loin et donner une formule pour le caractère de~$\pi(\mathscr{D})$ qui est valable pour \emph{tout} élément semi-simple régulier de~$G$,
en supposant cependant que la caractéristique de~$F$ est nulle et que~$p$ est assez grand (car la formule utilise l'exponentielle, qui ne converge agréablement que dans ce cas).

Les ingrédients de la formule sont moins élémentaires.
Les preuves s'appuient sur un immense travail \parencite{DeBacker_Reeder, DeBacker_Spice, Spice_Compositio, Spice_2} sur la structure générale des formules de caractère. 

À côté de l'ensemble~$\Gtopss$ des éléments topologiquements semi-simples, définissons l'ensemble $\Gtopnil$ des éléments \emph{topologiquement nilpotents} de~$G$ (modulo~$\Z_\G$), formé des~$\gamma \in G$ dont l'image $\overline{\gamma}$ dans~$\G/\Z_\G$ satisfait $\overline{\gamma}^{p^n} \underset{n \to \infty}{\longrightarrow} 1$.
Tout élément semi-simple régulier~$\gamma$ de~$G$ peut s'écrire sous la forme~$\gamma_0 \gamma_{0^+}$ où~$\gamma_0 \in \Gtopss$ et~$\gamma_{0^+} \in \Gtopnil$ commutent (« décomposition de~Jordan topologique », unique modulo le centre ; \cf\mbox{\cite{Spice_Jordan}}). 

Rappelons que l'algèbre de Lie~$\lie{s}=\Lie(S)(F)$ est naturellement munie de la filtration de Moy--Prasad attachée à l'unique point de~$\Immred(S)$.
On dispose donc d'un sous-espace~$\lie{s}_{0^+}$ de~$\lie{s}$.
Supposons que~$F$ soit de caractéristique nulle et que la restriction à  $\lie{s}_{0^+}$ de l'application exponentielle soit convergente, ce qui est le cas si~$p$ est assez grand.
Il existe alors un élément~$\xi$ de~$\lie{s}^\ast$ vérifiant $\theta(\exp(Y)) = \Lambda(\dual{\xi}{Y})$ pour tout~$Y\in\lie{s}_{0^+}$. 

La formule du caractère qui suit utilise les transformées de Fourier d'intégrales orbitales et on se contente d'indications grossières sur ces objets.
Soit~$\J$ la composante neutre du centralisateur de~$\gamma_0$ dans~$\G$.
C'est un groupe réductif connexe défini sur~$F$ qui contient~$\bS$.
Notons~$J = \J(F)$ et $\lie{j}=\Lie(\J)(F)$.
Pour tout élément~$\eta$ du dual~$\lie{j}^\ast$,  on peut former  une distribution $\mathscr{O}_{\eta}$ sur~$\lie{j}^\ast$, l'\emph{intégrale orbitale} en~$\eta$ :
par définition, si~$f \in C^\infty_c(\lie{j}^\ast)$, la valeur~$\dual{\mathscr{O}_\eta}{f}$ est l'intégrale de~$f$ sur l'orbite de~$\eta$ pour l'action coadjointe de~$J$.
La transformée de Fourier de~$\mathscr{O}_\eta$ est une distribution~$\widehat{\mathscr{O}}_\eta$ sur~$\lie{j}$.
Cette distribution peut être représentée par une fonction localement constante $\widehat{\mu}_\eta$ sur~$\lie{j}$.
(Définir correctement tous ces objets demande du soin et des choix de mesures ; \cf \cite[\S\,{}4]{DeBacker_Spice} ; \cite[\S\,{}4.2]{Kaletha_regular_supercuspidals}.)

\begin{theo}[\cite{Fintzen_Kaletha_Spice}] \label{th:formule_complete}
Supposons que~$F$ soit de caractéristique nulle et que la restriction à  $\lie{s}_{0^+}$ de l'application exponentielle soit convergente. 
Soient~$\gamma$ un élément semi-simple régulier de~$G$ et $\gamma= \gamma_0 \gamma_{0^+}$ une décomposition de Jordan topologique de~$\gamma$.
La valeur en~$\gamma$ du caractère de~$\pi(\mathscr{D})$ est 
 \begin{equation} \label{formule_complete} c_{G, J, \Lambda}\,{} \Delta(\gamma)^{-1/2} \sum \limits_{\substack{\crochets{g} \in S\backslash{}G/J \\ \mathrm{Ad}(g)(\gamma_0) \in S}} \crochets{a_{S, \Lambda, \vartheta}\cdot \vartheta}(\mathrm{Ad}(g)(\gamma_0)) \cdot \widehat{\mu}_{\mathrm{Ad}^\ast(g^{-1})(\xi)}(\log(\gamma_{0^+}))\end{equation}
 où la constante~$c_{G, J, \Lambda}$ est donnée par $e(G) e(J) \varepsilon_L(\Char^\ast(\T_\G)_{\C}-\Char^\ast(\T_\J)_{\C},\Lambda)$\pt
\end{theo}
Si~$p$ est assez grand, avec une condition explicite (\cf \cite[th.~1.7]{Kaletha_ICM}),  alors l'hypothèse sur l'exponentielle est vérifiée.

Remarquons que sur les éléments de~$\Sreg$, la seule différence avec la formule superficielle~\eqref{def_theta} provient de la contribution des intégrales orbitales.
Ces dernières se rattachent plus ou moins directement à la théorie des groupes finis.
Par exemple, si~$\pi(\mathscr{D})$ est de profondeur nulle, les intégrales orbitales dans~\eqref{formule_complete} sont directement liées aux caractères de représentations cuspidales de groupes finis \parencite{DeBacker_Reeder}.
Comme nous l'avons vu, tous les autres ingrédients sont en ligne directe avec les formules pour les groupes réels.
La formule~\eqref{formule_complete} combine élégamment ingrédients « réels » et  « finis ». 


\subsubsection{Théorème de Chan--Oi}\label{sec:chan_oi}
Ayant vu une forme « complète » de la formule du caractère, revenons à la formule simple~\eqref{def_theta}.
Elle n'est valable que sur les éléments superficiels ; mais dans le cas des groupes réels, nous avons vu qu'elle suffit à déterminer complètement une représentation irréductible de série discrète, via le th.~\ref{HC_parametrage}.
Dans le cas non archimédien, Henniart a montré que certaines représentations de~$\GL(n,F)$ peuvent être déterminées si l'on connaît la restriction de leur caractère, non à l'ensemble~$\Greg$, mais à un ensemble beaucoup plus petit d'éléments \emph{très réguliers} \parencite{Henniart_tresreg_1}.
Cette idée a été étendue récemment par~\textcite{Chan_Oi_2023} aux groupes réductifs généraux : ces auteurs donnent une condition suffisante pour que la formule~\eqref{def_theta} détermine une unique représentation supercuspidale comme dans le th.~\ref{HC_parametrage}. 

Soit~$\gamma$ un élément de~$G$.
Disons que~$\gamma$ est \emph{très régulier} s'il admet une décomposition de Jordan topologique $\gamma_0 \gamma_{0^+}$ (\no\ref{sec:formule_fine}) où~$\gamma_0$ est un élément semi-simple régulier de~$G$. 

Fixons un couple régulier $(\bS, \theta)$ (\no\ref{sec:def_non_singulieres}).
Chan et Oi introduisent une condition sur~$\bS$, qu'ils appellent \emph{inégalité d'Henniart}. La condition garantit l'existence de «\,{}suffisamment\,{}» d'éléments très réguliers, et  s'exprime sur un groupe fini attaché à~$\bS$. 

\begin{smallrema} 
Donnons un énoncé précis lorsque~$\G$ et~$\bS$ sont non ramifiés.
Soit~$\Sfini$ un tore défini sur~$k_F$ tel que~$\Sfini(k_F)$ s'identifie au quotient de Moy--Prasad~$S_0/S_{0^+}$.
Soit~$\mathsf{Z}^\ast$ le groupe des $k_F$-points de l'image dans~$\Sfini$ du centre de~$\G$.
Notons $\Sfini^\ast$ le quotient $\Sfini(k_F)/\mathsf{Z}^\ast$ et $\Sfini_{\textrm{t.reg.}}^\ast$ la partie de $\Sfini^\ast$ image de l'ensemble des éléments très réguliers de~$S_0$ par le quotient $S_0 \to \Sfini(k_F)$.
Soit $W({\Gfini},\Sfini)$  le groupe de Weyl de~$(\Gfini, \Sfini)$.
L'inégalité d'Henniart est l'assertion suivante :
\[ \Card(\Sfini^\ast) > 2\,{} \Card(W({\Gfini},\Sfini)) \,{} \Card(\Sfini^\ast\setminus \Sfini_{\textrm{t.reg.}}^\ast ). \]
Dans ce contexte, elle est toujours vérifiée si~$p$ est assez grand.
La formulation la plus générale utilise le groupe~$S/S_{0^+}$, qui provient d'un groupe réductif non connexe sur~$k_F$.
L'inégalité est alors plus contraignante et peut échouer dans des cas simples si~$\bS$ est totalement ramifié.
Voir \parencite[\S\,{}5.2, 6.2 et Appendice]{Chan_Oi_2023}.
\end{smallrema}

\begin{theo}[\cite{Chan_Oi_2023}]
Soit~$(\bS, \theta)$ un couple elliptique régulier. 
\begin{assertions}
\item Le caractère de~$\pi_{(\bS, \theta)}$ est donné par~\eqref{def_theta} sur les éléments très réguliers de~$S$\pt
\item Si~$\bS$ satisfait l'inégalité d'Henniart, alors $\pi_{(\bS, \theta)}$ est la seule représentation supercuspidale irréductible de~$G$ vérifiant~\textup{a)}\pt
\end{assertions}
\end{theo}

\section{Correspondance de Langlands et paquets supercuspidaux}\label{sec:LLC}


\subsection{Généralités sur la correspondance de Langlands locale}\label{sec:LLC_generalites}

Voir les textes de~Borel~(\cite*{Borel_Corvallis}, dans \cite*{Corvallis}) ou de~\textcite{Kaletha_Taibi}.


\subsubsection{} Il s'agit de relier les représentations irréductibles de~$G$ à des objets de nature arithmétique, ou galoisienne. Je ne donne pas ici de motivation pour la forme de ces objets, et me contente d'indications rapides pour introduire les notations. Les textes de Borel, Springer et Tate dans~\parencite*{Corvallis} forment une introduction canonique.

Soit~$\dualG$ le groupe réductif complexe dual de~$\G$.
C'est un 
groupe réductif complexe connexe dont les racines sont les coracines de~$\G$, et les coracines sont les racines de~$\G$.
Si l'on fixe un épinglage de~$\dualG$, l'action de~$\Gamma$ sur la donnée radicielle de~$\G$ induit une action de~$\Gamma$ sur~$\widehat{G}$
(action qui préserve l'épinglage donné et se factorise par le groupe de Galois d'une extension finie de~$F$).

Définissons~$\lgr{G}$ comme le produit semi-direct $\widehat{G}\rtimes\Gamma$ associé à cette action («\,{}forme galoisienne du $L$-groupe\,{}») ;
il vient donc avec un morphisme surjectif $\lgr{G}\to\Gamma$.

Soit~$\Weil_F$ le groupe de Weil de~$F$.
Si~$F$ n'est pas archimédien, c'est un sous-groupe du groupe de Galois~$\Gamma$ ; il contient le sous-groupe d'inertie $\Inertie_F$, et~\mbox{$\Weil_F/\Inertie_F$} est isomorphe à~$\Z$.
On munit~$\Weil_F$ de l'unique topologie qui fait de~$\Inertie_F$ un sous-groupe ouvert sur lequel la topologie induite est celle qu'induit la topologie de~$\Gamma$.
Si~$F = \R$, le groupe de~Weil est une extension de~$\C^\ast$ : on a $\Weil_\R = \langle \C^\ast, j\rangle$ où~$j^2=1$ et $jzj^{-1}=\overline{z}$ pour~$z \in \C^\ast$.

Une version du paramétrage proposé par Langlands utilise le \emph{groupe de {Weil--Deligne}} :
on pose~$\WeilDeligne_F =\Weil_F \times \SL(2,\C)$ pour~$F$ non archimédien, et \mbox{$\WeilDeligne_F = \Weil_F$} pour~$F=\R$. 

Un \emph{L-homomorphisme} est un homomorphisme  $\phi\colon \WeilDeligne_F\to\lgr{G}$ qui est continu, commute aux morphismes de~$\WeilDeligne_F$ et~$\lgr{G}$ vers~$\Gamma$, et pour lequel $\phi(\Weil_F) \cap \widehat{G}$ est formé d'éléments semi-simples.
Lorsque~$F$ n'est pas archimédien,  on demande de plus que la restriction de~$\phi$  au facteur $\SL(2,\C)$ de~$\WeilDeligne_F$ soit algébrique.

Le groupe~$\widehat{G}$ agit sur l'ensemble des $L$-homomorphismes par conjugaison au but.
On note~$\Langlands(G)$ l'ensemble des classes de conjugaison sous~$\dualG$ de $L$-homomorphismes.


\subsubsection{}
Sous sa forme la plus simple, la «~conjecture de Langlands locale pour~$G$~» est l'espoir qu'il soit possible de décrire une application de~$\Irr(G)$ dans~$\Langlands(G)$ qui soit à fibres finies, surjective quand~$G$ est quasi-déployé, et possède de nombreuses propriétés «\,{}naturelles\,{}» dans ce contexte. 
Autrement dit, à chaque $L$-homomorphisme~$\phi$, on espère attacher un ensemble~$\Paquet(\phi)$ de représentations irréductibles de~$G$, de façon~que
\begin{conditions}
\item le «\,{}paquet\,{}»~$\Paquet(\phi)$ soit fini et ne dépende que de la classe de~$\phi$ dans~$\Langlands(G)$ ;
\item l'ensemble~$\Irr(G)$ soit la réunion des paquets~$\Paquet(\phi)$ ;
\item et que la correspondance $\phi \leftrightsquigarrow \Paquet(\phi)$ vérifie une liste détaillée de conditions.
\end{conditions}
Les propriétés que l'on espère ont été précisées au fil des ans ; je ne les décrirai pas en détail.
On ne connaît d'ailleurs pas de liste de propriétés qui garantisse l'unicité de la correspondance, sauf pour $\GL(n)$ \parencite{Henniart_car}.
Je renvoie pour ces questions aux textes d'\textcite{Harris} et de~\textcite{Kaletha_Taibi}, ainsi qu'au rapport classique de Borel~(\cite*{Borel_Corvallis}, dans~\cite*{Corvallis}).
Énonçons cependant trois \emph{desiderata} classiques. 

Le premier concerne le cas où~$\G$ est un tore~$\T$. Dans ce cas, $\Irr(G)$ s'identifie à~$\Char(T)$.
Par ailleurs, le tore complexe $\widehat{T}$ peut être pris égal à $\Char^\ast(\T) \otimes_\Z \C$, un $L$-homomorphisme s'identifie à un $1$-cocycle continu de~$\Weil_F$ dans~$\widehat{T}$, et $\Langlands(T)$ s'identifie au groupe de cohomologie continue $H^1(\Weil_F, \widehat{T})$.
Comme l'a montré Langlands, la théorie du corps de classes permet alors de construire une bijection entre~$\Char(T)$ et $\Langlands(T)$, que l'on adopte pour la correspondance
(avec des paquets réduits à un caractère de~$T$). 

Le deuxième concerne les paquets de série discrète.
Si~$\M$ est un sous-groupe de Levi de~$\G$ défini sur~$F$, on dispose d'un plongement \mbox{$\iota_\M\colon \lgr{M}\to\lgr{G}$,} qui est canonique à conjugaison près par~$\widehat{G}$ au but. 
On espère que les représentations de série discrète forment la réunion des paquets~$\Pi(\phi)$ associés aux paramètres~$\phi$ qui ne se factorisent par aucun des plongements~\mbox{$\iota_\M$, $\M \neq \G$}. 
On qualifie de \emph{discrets} les $L$-homomorphismes~$\phi$ ayant cette propriété. On espère aussi que les paquets~$\Pi(\phi)$ associés soient tous non~vides.

Le troisième concerne la compatibilité avec l'induction parabolique lorsque~$F$ n'est pas archimédien.
Notons $\norm{\cdot}\colon \Weil_F \to \R_+^\ast$ l'homomorphisme décrit par Tate dans \parencite*[vol.~2, p.~7, \no1.4.6]{Corvallis} et~$\iota_{\Weil_F}\colon \Weil_F\to \WeilDeligne_F$ l'application $w \mapsto (w, \diag(\norm{w}, \norm{w}^{-1}))$.
Soit~$M$ un sous-groupe de Levi de~$G$, facteur de Levi d'un sous-groupe parabolique~$P$ ; et soit~$\sigma$ une représentation irréductible de~$M$. On souhaite que pour tout sous-quotient irréductible~$\pi$ de~$\mathrm{Ind}_{P}^{G}(\sigma)$, les paramètres associés $\phi_\sigma\colon \WeilDeligne_F \to \lgr{M}$ et $\phi_\pi \colon \WeilDeligne_F \to \lgr{G}$  vérifient $\iota_M \circ \phi_\sigma \circ \iota_{\Weil_F} = \phi_\pi \circ \iota_{\Weil_F}$.
Cela implique que si~$\phi$ est discret et trivial sur le facteur~$\SL(2,\C)$ de~$\WeilDeligne_F$, alors~$\Pi(\phi)$ ne contient que des représentations supercuspidales. On dit qu'un paramètre~$\phi$ est \emph{supercuspidal} s'il est discret et trivial sur~$\SL(2,\C)$.
Si~$F$ est archimédien, on dira «\,{}supercuspidal\,{}» pour «\,{}discret\,{}».


\subsubsection{} \label{sec:caractere_stable}
Disons un mot de la notion de \emph{caractère stable}, qui permet de caractériser la correspondance de Langlands supercuspidale construite par~\textcite{Kaletha_regular_supercuspidals, Kaletha_Lpackets}.

Soit~$\phi$ un paramètre discret. 
Chaque élément du paquet~$\Paquet(\phi)$ vient avec un caractère~$\Theta_\pi$, que l'on peut voir comme une fonction de $\Greg$ dans~$\C$.
On espère qu'il existe une  combinaison linéaire $S\Theta^G_{\phi}$ des fonctions $\Theta_\pi$, $\pi \in \Paquet(\phi)$, qui vérifie : 
\begin{conditions}
\item La fonction $S\Theta^G_{\phi}$ est invariante par \emph{conjugaison stable}, c'est-à-dire par les automorphismes de~$G$ de la forme~$h \mapsto g h g^{-1}$ où~$g$ est un élément du groupe $\G(\Fsep)$  ;
\item Si~$\Sigma$ est une partie de~$\Paquet(\phi)$ distincte de~$\Paquet(\phi)$, alors aucune combinaison linéaire des fonctions~$\Theta_\pi$, $\pi \in \Sigma$, n'est invariante par conjugaison stable.
\end{conditions}
Dans ce cas, on dit que $S\Theta^G_{\phi}$ est \emph{atomiquement stable}.
Une telle combinaison est unique à une constante multiplicative près ; c'est « le » \emph{caractère stable} du paquet~$\Pi(\phi)$.

Ce caractère stable est lié à la structure interne que l'on espère voir vérifiée par les paquets.
Je me limiterai à quelques remarques dans le cas où~$G$ est quasi-déployé, car le cas général nécessite beaucoup de soin.
En particulier, il est naturel d'y faire intervenir une notion de « forme intérieure rigide » qui n'a été bien comprise pour~$F \neq \R$ que récemment ; je renvoie au texte de~\textcite{Kaletha_Taibi} pour une description détaillée.

 On peut considérer, dans le groupe connexe~$\dualG$, le centralisateur $\mathrm{Cent}_{\widehat{G}}(\im(\phi))$.
 Il contient le groupe $Z(\widehat{G})^\Gamma$ des éléments $\Gamma$-invariants du centre de~$\widehat{G}$.
 Considérons le quotient $\mathrm{Cent}_{\widehat{G}}(\im(\phi))/Z(\widehat{G})^\Gamma$ et notons~$\mathbb{S}_\phi$ son groupe de composantes.
 Si~$\phi$ est un paramètre discret, le groupe~$\mathbb{S}_\phi$ est fini  (et si~\mbox{$F=\R$}, il  est abélien). 
 
 Dans le cas où~$F$ n'est pas archimédien, on s'attend à ce qu'il existe une bijection naturelle entre $\Paquet(\phi)$ et l'ensemble fini $\Irr(\mathbb{S}_\phi)$ (cette bijection doit dépendre d'un choix de « donnée de Whittaker » pour~$G$).
 Si l'on note~$\rho_\pi$ la représentation de~$\mathbb{S}_\phi$ correspondant à un élément~$\pi$ de~$\Paquet(\phi)$, on s'attend alors à ce que $S\Theta^G_\phi$ soit égal à $\sum_{\pi \in \Pi(\phi)} \dim(\rho_\pi) \Theta_\pi$.

De plus, on pense que des relations explicites et précises, les \emph{identités endoscopiques}, sont vérifiées par le caractère stable $S\Theta^G_{\phi}$ et par les caractères stables $S\Theta^H_{\phi}$ attachés à certains groupes quasi-déployés~$H$ pour lesquels~$\phi$ induit canoniquement un paramètre de Langlands discret. Je ne peux décrire ici précisément ces identités (\cf \cite[\S\,{}2.2]{Kaletha_ICM}); disons simplement que lorsqu'elles sont vérifiées, il est possible de reconstituer les éléments~$\pi$ du paquet~$\Pi(\phi)$, et leur paramétrage par $\Irr(\mathbb{S}_\phi)$, à partir des fonctions~$S\Theta^H_\phi$.


\subsection{Groupe dual et correspondance de Langlands pour le revêtement~$\Spm$}\label{sec:LLC_pour_revetement}

Supposons donnés, comme au \S\,{}\ref{sec:revetement_double}, un tore~$\bS$  défini sur~$F$ et un $F$-plongement $\bS \to \G$ dont l'image est un tore maximal.
Le revêtement double~$\Spm$ construit au \S\,{}\ref{sec:definition_du_revetement} n'est pas le groupe des points rationnels d'un groupe réductif algébrique connexe.
Cependant, nous allons définir un groupe $\lgr{S}_\pm$ qui joue le rôle de groupe dual pour~$\Spm$.

 Soit~$E$ une extension de~$F$ sur laquelle~$\bS$ se déploie.
 Notons~$\Gamma_{E|F}$ son groupe de Galois : c'est un quotient de~$\Gamma$.
 Le groupe $\lgr{S}_\pm$ à définir provient d'une extension, non scindée en général, de~$\Gamma_{E|F}$ par le tore complexe~$\widehat{S}$.
 Spécifier une telle extension revient à spécifier un $2$-cocycle de~$\Gamma_{E|F}$ dans~$\widehat{S}$ (\cf \cite[\S\,{}16, \no5]{Bourbaki_A8}). 
 
 Choisissons pour cela une \emph{jauge}, c'est-à-dire une fonction $p \colon \Racines \to \{\pm 1\}$ qui attribue un signe à chaque racine.
 On définit alors une fonction $t_p$ de $\Gamma_{E|F} \times \Gamma_{E|F}$ dans~$\widehat{S}$ de la manière suivante.
 Si~$\sigma, \tau$ sont des éléments de~$\Gamma_{E|F}$, notons $\lambda(\sigma, \tau)$ la somme des éléments~$\alpha$ de~$\Racines$ vérifiant~$p(\alpha)=1$, $p(\sigma^{-1}\alpha)=-1$ et $p\crochets{(\sigma\tau)^{-1}\alpha}=1$.
 Alors $\lambda(\sigma,\tau)$ est un élément de $\Char^\ast(S)$, qui s'identifie canoniquement à $\Char_{\ast}(\widehat{S})=\mathrm{Hom}(\C^\ast, \widehat{S})$.
 On peut donc évaluer  $\lambda(\sigma,\tau)$ en $(-1)$ et on obtient un élément de~$\widehat{S}$, que l'on note~$t_{p}(\sigma, \tau)$. 
 
 La fonction~$t_{p}$ est un $2$-cocycle. Définissons $\lgr{S}_{\pm}^{E|F}$ comme l'extension de~$\Gamma_{E|F}$ par~$\widehat{S}$ associée à ce cocycle : il s'agit du groupe d'ensemble sous-jacent $\widehat{S}\times \Gamma_{E|F}$ muni de la loi \mbox{$(s_1, \sigma_1) \cdot (s_2, \sigma_2) = (s_1 \,{} s_2 \,{} t_p(\sigma_1,\sigma_2), \sigma_1 \sigma_2)$}. En utilisant la projection $\Gamma\to \Gamma_{E|F}$, on  déduit de~$\lgr{S}_{\pm}^{E|F}$ une extension~$\lgr{S}_{\pm}$ de~$\Gamma$.
 Cette extension dépend du choix de~$p$ ; mais si l'on choisit une autre jauge, les extensions associées sont canoniquement isomorphes.

Appelons \emph{paramètre de Langlands pour $\Spm$} un morphisme continu $\phi\colon \Weil_F \to\lgr{S}_\pm$ qui commute aux applications naturelles de ces deux groupes vers~$\Gamma$, et disons que des paramètres de Langlands sont \emph{équivalents} s'ils sont conjugués sour~$\widehat{S}$.
 
\textcite[\S\,{}3.4]{Kaletha_Covers} décrit une bijection naturelle entre  caractères spécifiques de~$\Spm$ et classes d'équivalence de paramètres de Langlands pour $\Spm$.
Sans entrer dans les détails, supposons un instant que l'on dispose d'une donnée de transfert réduite $(\chi_\alpha)_{\alpha\in\RacinesSym}$.
Nous avons vu au \no\ref{def_chi_data} qu'elle détermine un caractère spécifique~$\chi$ de~$\Spm$.
Kaletha lui associe un morphisme $\phi_{\chi}\colon \Weil_F \to\lgr{S}_{\pm}$ : la construction est relativement technique mais reprend des ingrédients déjà présents chez~\textcite{Langlands_Shelstad}.
La multiplication par~$\phi_\chi$ fournit alors une bijection entre morphismes $\Weil_F \to \lgr{S}$ et  $\Weil_F \to \lgr{S}_\pm$.
Or, du côté du groupe, la multiplication par~$\chi$ induit une bijection entre  caractères spécifiques de~$\Spm$ et caractères continus de~$S$.
En composant ces deux bijections avec la correspondance de Langlands pour~$S$, on obtient une bijection entre caractères spécifiques et classes d'équivalence de paramètres de Langlands.
Elle ne dépend pas du choix de donnée de transfert ; c'est la correspondance de Langlands pour~$\Spm$.


\subsection{Construction des paquets supercuspidaux}\label{sec:construction_LLC}

Si~$F$ est non archimédien, on se place toujours dans le cas bien modéré (\no\hyperref[sec:cas_bien_modere]{1.5}).

Soit~$\phi\colon \WeilDeligne_F \to \lgr{G}$ un paramètre supercuspidal. Kaletha décrit un paquet fini~$\Paquet(\phi)$ de représentations de série discrète  de~$\G$ ~\mbox{\parencite*{Kaletha_regular_supercuspidals,  Kaletha_Covers, Kaletha_Lpackets}} ;
si~$F$ n'est pas archimédien, le paquet~$\Paquet(\phi)$ est entièrement formé de représentations supercuspidales non~singulières.
Toute représentation supercuspidale non~singulière est contenue dans un tel paquet, et Kaletha montre que ces paquets ont toutes les propriétés que l'on espère.

Décrivons cette construction.
Lorsque~$F=\R$, elle coïncide avec celle de \textcite{Langlands_CIRRAG}, simplifiée par l'usage des revêtements doubles \parencite{AV1, AV_Contragredient}. 


\subsubsection{Du paramètre au tore} 
Identifions le paramètre  \mbox{$\phi\colon \WeilDeligne_F \to \lgr{G}$} à sa restriction à~$\Weil_F$.
La première observation cruciale est que~$\phi$ détermine un tore maximal~$\widehat{S}$ de~$\dualG$, défini comme suit dans le cas non archimédien.
Considérons la composante neutre~$\mathrm{Cent}_{\widehat{G}}(\phi(\Inertie_F))^0$ du centralisateur de l'image du sous-groupe d'inertie.
Cette composante neutre est contenue dans le centralisateur~$\widehat{M}=\mathrm{Cent}_{\widehat{G}}(\phi(\InertieSauvage_F))$ de l'image du sous-groupe d'inertie sauvage.
Le groupe~$\widehat{M}$  est un sous-groupe de Levi de~$\widehat{G}$. Posons
\[ \widehat{S}=\mathrm{Cent}_{\widehat{M}}\crochets{\mathrm{Cent}_{\widehat{G}}(\phi(\Inertie_F))^0}.\]
Kaletha observe que~$\widehat{S}$ est un tore maximal de~$\dualG$. 
L'image de~$\phi$ est contenue dans le normalisateur $\mathrm{Norm}_{\,\lgr{G}}(\widehat{S})$ (\cite*[lemme~5.2.2]{Kaletha_regular_supercuspidals}; \cite*[lemme~4.1.3]{Kaletha_Lpackets}).
Quitte à remplacer~$\phi$ par un conjugué sous~$\widehat{G}$, on peut supposer que~$\Gamma$ normalise le tore~$\widehat{S}$ ainsi qu'un sous-groupe de Borel~$\widehat{B}_0$ de~$\widehat{G}$ ayant~$\widehat{S}$ comme facteur de Levi.
Pour l'analogue dans le cas archimédien, on remplace $\phi(\Inertie_F)$ par $\phi(\C^\ast)$ ; \cf \cite{AV_Contragredient}, \S\,{}6.
 
Dans ces conditions, il existe un tore maximal elliptique de~$\G$ dont le tore dual est~$\widehat{S}$, mais il n'est déterminé qu'à conjugaison stable près.
Comme nous allons le voir, cette ambiguïté est l'une des raisons pour l'existence de contributions distinctes  au paquet~$\Pi(\phi)$.
Soit~$\bS$ le tore défini sur~$F$ dont~$\widehat{S}$ est le tore dual ; on peut construire non seulement un plongement~$\bS \to \G$, mais  une \emph{classe de conjugaison sous~$G(\Fsep)$} de tels plongements, cette classe~$\Plongements$ étant stable par l'action de~$\Gamma$ à la source et au but. 

\begin{smallrema}
Plus précisément, supposons donnés des épinglages $(\T, \B, \{X_\alpha\})$ et~$(\widehat{T}, \widehat{B}, \{Y_{\ch{\alpha}}\})$ de~$\G$ et~$\widehat{G}$ respectivement.
Alors~$\widehat{T}$ et~$\widehat{S}$ sont conjugués dans~$\widehat{G}$, donc il existe un automorphisme intérieur de~$\widehat{G}$ qui se restreint en un isomorphisme~$\widehat{j}\colon \widehat{T}\to \widehat{S}$.
Soit~$j\colon \bS \to \T$ le morphisme de $F$-tores dual de~$\widehat{j}$, et $j\colon \bS \to \G$ sa composée avec le plongement canonique de~$\T$ dans~$\G$.
La classe de conjugaison de~$j$ sous~$\G(\Fsep)$ est stable par~$\Gamma$ et indépendante des choix d'épinglage.
C'est elle que l'on prend pour~$\Plongements$.
(Pour tout cela, voir \cite[\S\,{}5.1]{Kaletha_regular_supercuspidals}.)
\end{smallrema}

Parmi les plongements appartenant à~$\Plongements$, il en est qui sont définis sur~$F$ ; ces $F$-plongements envoient~$\bS$ sur un tore maximal elliptique de~$\G$.
Si~$j\colon \bS\to\G$ appartient à~$\Plongements$ et est défini sur~$F$, alors on peut reprendre la construction du \S\,{}\ref{sec:revetement_double} pour construire un revêtement double~$\Spm$ du groupe $S=\bS(F)$, indépendant du choix de~$j$. 


\subsubsection{Plongement de~$\lgr{\Spm}$ dans~$\lgr{G}$} \label{sec:construction_L_plongement}
Nous allons voir que~$\lgr{\Spm}$ se plonge dans~$\lgr{G}$ de façon à peu près canonique, et que le paramètre~$\phi\colon \Weil_F \to\lgr{G}$ se factorise en un paramètre à valeurs dans~$\lgr{\Spm}$.
Pour construire le paquet~$\Paquet(\phi)$, on applique   la correspondance du \S\,{}\ref{sec:LLC_pour_revetement} à ce paramètre  \mbox{$\Weil_F\to\lgr{\Spm}$}, ce qui fournit un caractère spécifique de~$\Spm$ ;
puis les constructions du \S\,{}\ref{sec:supercuspidales_nonsingulieres}, qui donnent les représentations à incorporer à~$\Paquet(\phi)$.

Pour avoir une réalisation concrète de~$\lgr{\Spm}$ comme en~\ref{sec:LLC_pour_revetement}, précisons un choix de jauge.
Nous avons fixé ci-dessus un sous-groupe de Borel~$\widehat{B}_0\subset \widehat{G}$ contenant~$\widehat{S}$.
Il détermine un ensemble de racines positives dans le système~$\Racines$ des racines de~$\bS$ dans~$\G$.
On prend pour~$p$ la jauge  qui attribue le signe~$+1$ aux racines positives et à elles seules. 

La construction d'un plongement~$\lgr{\Spm}\to\lgr{G}$ repose sur deux propriétés du groupe de Weyl $\Omega = \mathrm{Norm}_{\widehat{G}}(\widehat{S})/\widehat{S}$.
La première est qu'il permet de réaliser l'action du groupe de Galois~$\Gamma$ sur~$\widehat{S}$ : pour tout~$\gamma \in \Gamma$, il existe un unique élément~$\omega_\gamma$ de~$\Omega$ tel que les automorphismes de~$\widehat{S}$ induits par~$\gamma$ et~$\omega_\gamma$ coïncident.
La seconde est que si l'on se donne un épinglage  $(\widehat{S}, \widehat{B}, \{X_\alpha\})$ de~$\widehat{G}$, alors on dispose d'une section \emph{canonique} \mbox{$\tau \colon \Omega \to~\mathrm{Norm}_{\widehat{G}}(\widehat{S})$} de l'application de passage au quotient de $ \mathrm{Norm}_{\widehat{G}}(\widehat{S})$ dans $\Omega= \mathrm{Norm}_{\widehat{G}}(\widehat{S})/\widehat{S}$ : l'application~$\tau$ est la \emph{section de Tits} (\cite*{Tits_section} ; voir aussi \cite[p.~17]{Kaletha_Covers}).
Or on peut raffiner le couple~$(\widehat{S}, \widehat{B}_0)$ en un épinglage $\Gamma$-invariant de~$\widehat{G}$ \parencite[Fact~3.2.2]{Kaletha_trace}, unique à conjugaison près par le groupe $\widehat{G}^\Gamma$.
Considérons alors l'application 
\begin{eqnarray*} &  \lgr{\Spm} &\to~ \widehat{G} \rtimes \Gamma \\ 
&(s, \gamma) & \mapsto~ (s \,{} \tau(\omega_\gamma), \gamma) \pv \end{eqnarray*}
elle est injective, d'image contenue dans~$\mathrm{Norm}_{\widehat{G}}(\widehat{S})$ et commute aux projections vers~$\Gamma$.
En réinterprétant des observations de~\textcite{Langlands_Shelstad}, Kaletha signale que c'est un morphisme de groupes.
On obtient donc bien un plongement de~$\lgr{\Spm}$ dans~$\lgr{G}$, unique à conjugaison près par l'action de~$\widehat{G}^\Gamma$ au but.


\subsubsection{Construction du paquet~$\Paquet(\phi)$}
Déduisons du plongement~$\lgr{\Spm}\to \lgr{G}$ et du paramètre~$\phi$ un paramètre de Langlands~$\Weil_F \to\lgr{\Spm}$,  puis appliquons la correspondance de Langlands pour~$\Spm$ (\S\,{}\ref{sec:LLC_pour_revetement}) pour obtenir un caractère spécifique~$\vartheta$ de~$\Spm$. 

Soit~$j\colon \bS\to \G$ un $F$-plongement dans la classe stable~$\Plongements$ déterminée par le paramètre~$\varphi$.
Il identifie~$j(\bS)$ à un tore maximal elliptique de~$\G$, induit un isomorphisme de~$\Spm$ sur le revêtement $\Spm^j$ de~$S^j=j(\bS)(F)$, et on en déduit un caractère spécifique~$\vartheta^j$ de~$\Spm^j$.
Si~$F=\R$, ce caractère spécifique est régulier et on sait lui associer une représentation de série discrète~$\pi_{(j(\bS), \theta^j)}$ (th.~\ref{HC_parametrage}).
Si~$F$ n'est pas archimédien, on peut trouver explicitement un caractère~$\theta^j$ de~$S$ qui est non~singulier (\no\ref{sec:def_non_singulieres}) et à partir duquel on obtient le caractère spécifique~$\vartheta$ \emph{via} le \no\ref{def_chi_data} (\cf \cite[\S\,{}4.1--4.3]{Kaletha_Covers}).
Nous  avons vu au  \no\ref{sec:def_non_singulieres} comment lui associer  une représentation~$\pi_{(j(\bS), \theta^j)}$, somme finie de représentations supercuspidales irréductibles.

\begin{enonce*}{Définition}
Le paquet~$\Paquet(\phi)$ est la réunion des classes d'équivalence des  représentations que l'on obtient comme composantes irréductibles des représentations~$\pi_{(j(\bS), \theta^j)}$, lorsque~$j$ parcourt l'ensemble des $F$-plongements appartenant à~$\Plongements$\pt
\end{enonce*}

Si l'on part de  $F$-plongements~$j_1, j_2\colon \bS\to \G$ conjugués par l'action de~$G=\G(F)$, alors les représentations $\pi_{(j_1(\bS), \theta^{j_1})}$ et $\pi_{(j_2(\bS), \theta^{j_2})}$ sont équivalentes.
Mais il peut arriver : 
\begin{assertions}
\item que~$j_1$ et~$j_2$ aient même image et bien que conjugués sous~$\G(\Fsep)$, ne soient pas conjugués sous~$G$ ; alors $j_1(\bS)(F)=j_2(\bS)(F)$, mais les représentations \mbox{$\pi_{(j_1(\bS), \theta^{j_1})}$ et $\pi_{(j_2(\bS), \theta^{j_2})}$} n'ont aucun constituant commun ;
\item  que   $j_1(\bS)(F)$ et~$j_2(\bS)(F)$ soient stablement conjugués sans être rationnellement conjugués, et $\pi_{(j_1(\bS), \,{} \theta^{j_1})}$ et $\pi_{(j_2(\bS),\,{} \theta^{j_2})}$ n'ont alors aucun constituant commun. 
\end{assertions}
Ce sont ces deux situations, ainsi que le fait que les représentations  $\pi_{(j(\bS),\,{} \theta^j)}$ puissent être réductibles dans le cas où~$\theta^j$ est non~singulier sans être régulier, qui fournissent des éléments distincts du paquet~$\Pi(\phi)$. Dans le cas~$F=\R$, seul le cas~\ass{a} peut se présenter.

Si~$F$ n'est pas archimédien, le paquet~$\Paquet(\phi)$ est entièrement formé de représentations supercuspidales non~singulières.
Réciproquement, si~$\pi$ est une représentation supercuspidale non~singulière, on voit en inversant la construction ci-dessus qu'il existe un paramètre supercuspidal~$\phi$ tel que~$\pi$ soit contenue dans le paquet~$\Paquet(\phi)$. 

Si le paquet~$\Paquet(\phi)$ contient une représentation supercuspidale régulière, alors il est entièrement formé de supercuspidales régulières.
Les paramètres « réguliers »~$\phi$ pour lesquels cela se produit ont une caractérisation simple \parencite[déf. 5.2.3]{Kaletha_regular_supercuspidals}.


\subsection{Quelques propriétés des paquets supercuspidaux}


\subsubsection{}
Pourquoi penser que le \S\,{}\ref{sec:construction_LLC} donne « la bonne » construction des paquets~$\Paquet(\phi)$ associés aux paramètres supercuspidaux ?
Donnons de brèves indications sur les propriétés de la construction de Kaletha.
Nous le disions, elle coïncide avec celle de~\textcite{Langlands_CIRRAG} lorsque~$F = \R$.
Lorsque~$\G = \GL(n)$, on sait qu'elle coïncide avec la constructions équivalentes de Harris--Taylor et d'Henniart~\parencite{Oi_Tokimoto}.
Pour le rapport avec la construction de~\textcite{Moy_1986}, voir~\textcite{Bushnell_Henniart}. 


\subsubsection{}  \label{sec:endoscopie} Une propriété importante et délicate est que les paquets du \S\,{}\ref{sec:construction_LLC} ont «\,{}la bonne\,{}» structure interne et que leurs caractères stables vérifient les identités endoscopiques attendues.
Évoquons brièvement cela en se limitant au cas où~$G$ est quasi-déployé. 

Soit~$\phi$ un paramètre de Langlands supercuspidal et~$\mathbb{S}_\phi$ le groupe des composantes du \no\ref{sec:caractere_stable}.
Ce groupe est abélien si~$\Pi(\phi)$ est formé de représentations régulières, mais pas en général.
Kaletha décrit une bijection entre~$\Irr(\mathbb{S}_\phi)$ et le paquet~$\Paquet(\phi)$ (en supposant fixée une donnée de Whittaker pour~$G$).
Dans le cas d'un paquet supercuspidal  non régulier, cette description nécessite une compréhension fine de la décomposition des représentations $\pi_{(j(\bS),\,{} \theta^j)}$ ; voir le \no\ref{sec:reg_prof_nulle} pour cet aspect difficile de~\parencite{Kaletha_Lpackets}. 

On dispose d'une formule simple pour la combinaison $S\Theta^G_\phi=\sum_{\pi \in \Pi(\phi)} \dim(\rho_\pi) \Theta_\pi$.
Pour en donner une idée dans le cas des éléments superficiels, fixons un élément~$g$ de~\mbox{$\Greg \cap \Gtopss$}.
Alors $S\Theta_\phi(g)=0$ si~$g$ n'est dans l'image d'aucun $F$-plongement de la classe~$\Plongements$.
Si on identifie~$\bS$ au tore maximal elliptique de~$\G$ donné par l'image d'un tel plongement, on dispose d'un caractère spécifique~$\vartheta$ de~$S = \bS(F)$ ; et si~$g \in S$, on~a 
\[ S\Theta^G_\phi(g)=e(G)^{-1} c_{G, \bS, \Lambda} \cdot \Delta(s)^{-1/2}\cdot \sum_{w \in W(\G,\bS)(F)} \crochets{a_{\bS, \Lambda, \theta} \cdot \vartheta} (w^{-1}\cdot g). \]
Outre  la constante, la seule différence avec~\eqref{def_theta} est que la somme est prise sur les $F$-points du groupe de Weyl absolu (\cf \no\ref{sec:weyl}), ce qui permet à la valeur~$S\Theta_\phi(g)$ d'être indépendante du plongement choisi et constante sur la classe de conjugaison \emph{stable} de~$g$. 

Lorsque~$F$ est de caractéristique~0 et~$p$ est assez grand pour que l'exponentielle converge, on dispose d'une formule valable hors du lieu superficiel, analogue au th.~\ref{th:formule_complete} (\cf \cite[p.~14]{Kaletha_ICM}).
Elle fournit une fonction \emph{invariante par conjugaison stable}. Et surtout, si l'on considère les fonctions $S\Theta_\phi^H$ associées à $\phi$ et aux groupes endoscopiques~$H$  évoqués au \no\ref{sec:caractere_stable}, on constate qu'elles  \emph{vérifient les identités endoscopiques attendues} --- au moins pour~$\phi$ régulier \parencite[th.~4.4.4]{Fintzen_Kaletha_Spice}.
Lorsque~$\Pi(\phi)$ est un paquet supercuspidal non régulier, une partie significative des identités endoscopiques est vérifiée (\loccit). Démontrer les identités restantes est ouvert. 

\begin{smallrema}
On connaît d'autres propriétés des paquets supercuspidaux.
Par exemple, si un paquet~$\Pi(\phi)$ est formé de séries discrètes, \textcite{Hiraga_Ichino_Ikeda} proposent une formule pour le degré formel \parencite[p.~423, déf.~4]{Bourbaki_TSV} des représentations du paquet~$\Paquet(\phi)$, en termes de facteurs $\eps$ et~$L$ associés à~$\phi$.
On connaît désormais le degré formel des représentations de~Yu \parencite{Schwein} et les paquets supercuspidaux vérifient la conjecture d'Hiraga--Ichino--Ikeda \parencite{Schwein, Ohara}.
\end{smallrema}


\section{Au-delà du cas modéré}\label{sec:beyond}


\subsection{Constructions de supercuspidales en petite caractéristique résiduelle}


\subsubsection{Les constructions disponibles pour les groupes classiques}
Les représentations supercuspidales de la plupart des groupes classiques sont connues en petite caractéristique résiduelle, grâce à des raffinements de la méthode introduite par Bushnell et Kutzko pour $\GL(n)$ ($n$ et~$p$ arbitraires) --- pour laquelle je renvoie à l'exposé d'Henniart~\parencite*{Henniart_1991_seminaire}.
Le cas de $\SL(n)$,~$n$ et~$p$ arbitraires, est traité par Bushnell et Kutzko~\parencite*{Bushnell_Kutzko_SLN}.
En caractéristique résiduelle~$p \neq 2$, les représentations supercuspidales des groupes classiques sont construites par \textcite{Stevens_2008} en s'appuyant, notamment, sur une série de travaux de Morris.
Les formes intérieures de~$\GL(n)$ sont traitées par Sécherre et Stevens~\parencite*{Secherre_Stevens} ; le cas des formes intérieures anisotropes avait été résolu par Broussous~\parencite*{Broussous_1998}, généralisant un travail de Zink en caractéristique nulle. Plus récemment,  \textcite{Skodlerack} aborde les formes quaternioniennes des groupes classiques ($p\neq 2$).


\subsubsection{Représentations épipélagiques}
Vers 2008, Gross et Reeder ont introduit une construction très simple de représentations supercuspidales à partir de données assez différentes de celles de Yu.
Leur construction~\parencite*{Gross_Reeder} fut étendue par Reeder et~Yu~\parencite*{Reeder_Yu} à l'aide de la filtration de Moy--Prasad.
Cela permet de définir une nouvelle classe de représentations supercuspidales qui sont toujours de «\,{}petite\,{}» profondeur, baptisées pour cette raison \emph{épipélagiques}.
Reeder et Yu les étudiaient «\,{}pour~$p$ grand\,{}», mais \textcite{Fintzen_Romano}, puis \textcite{Fintzen_JEMS}, ont montré qu'hors du cas bien modéré, on obtient ainsi des représentations supercuspidales irréductibles qui échappent à la construction de Yu.
Voir aussi \textcite{Romano}, qui décrit  les équivalences entre représentations de Reeder--Yu et certains de leurs paramètres de Langlands. 


\subsection{Représentations supercuspidales unipotentes}
Il s'agit de représentations de profondeur nulle (\no\ref{sec:profondeur_nulle}), donc obtenues à partir de représentations cuspidales de groupes réductifs finis ; mais loin des représentations~$\DeLus_{\Sfini, \xi}$ du~\mbox{\S\,{}\ref{sec:def_reg_sc}}, on utilise des représentations unipotentes des groupes finis.
C'est en quelque sorte le cas « le plus éloigné » des représentations non singulières, étudié par Lusztig, et par Morris, dans les années 1990.
Pour des progrès récents, notamment sur la correspondance de Langlands, voir \textcite{Solleveld_Unipotent, Solleveld_Unipotent_2} ; Feng, Opdam et Sollevled~\parencite*{Feng_Opdam_Solleveld}.


\subsection{Cuspidalité pour les paramètres de Langlands enrichis}
Supposons qu'on dispose d'une construction de la correspondance de Langlands pour~$G$ qui coïncide avec celle de Kaletha sur les paramètres supercuspidaux.
Si~$\pi$ est une représentation supercuspidale irréductible de~$G$ et si $\pi$ est singulière, alors elle se trouve  dans un paquet~$\Paquet(\phi)$ qui ne peut être entièrement supercuspidal.
Étant donné un paramètre~$\phi\colon \WeilDeligne_F \to \lgr{G}$ discret non supercuspidal, comment savoir si~$\Paquet(\phi)$ contient des représentations supercuspidales ?
Et si l'on connaît la structure interne du paquet~$\Paquet(\phi)$ en termes de groupes de composantes, comment y distinguer les éléments supercuspidaux ?
Aubert, Moussaoui et Solleveld ont proposé une réponse~\parencite*{AMS}.


\subsection{Structure du dual lisse et de l'ensemble des paramètres de Langlands}
Compte tenu des théorèmes de Jacquet mentionnés au \no\hyperref[sec:reduction_jacquet]{1.2}, on peut imaginer qu'une étude fine de l'induction parabolique permette de trouver des renseignements précis sur la structure de~$\Irr(G)$, modulo la connaissance des représentations supercuspidales de sous-groupes de Levi de~$G$.
En dehors du cas de $\GL(n)$, ce problème  est difficile.
Mais il y a eu récemment des avancées remarquables, utilisant généralement la structure de certaines algèbres de Hecke.
Mentionnons des travaux de  \textcite{Heiermann_selecta} et la conjecture d'Aubert--Baum--Plymen--Solleveld sur la structure du dual lisse \parencite*{ABPS_Takagi, ABPS_Orsay}, dont l'essentiel est maintenant résolu (voir~\cite{Solleveld_2022}, pour une des dernières~étapes).

Ces renseignements ont des implications pour la correspondance de Langlands.
Dans le cas de $\GL(n)$, la construction de la correspondance de Langlands se ramène à celle de la partie supercuspidale de la correspondance pour les sous-groupes de Levi de~$\GL(n)$, qui sont des produits de groupes $\GL(r)$ avec~$r \leq n$.
Pour~$G$ arbitraire, on peut se demander jusqu'où il est possible d'aller en combinant  la correspondance supercuspidale et une étude fine de l'induction parabolique.
On peut avancer dans cette voie en comparant certaines des algèbres de Hecke ci-dessus et des algèbres analogues associées aux paramètres de Langlands.
Par exemple, en s'appuyant sur les travaux ci-dessus mentionnés et sur la correspondance supercuspidale de Kaletha, \textcite{Aubert_Xu} ont proposé  une construction explicite de la correspondance de Langlands pour le groupe exceptionnel~$\mathbf{G_2}$ sur un corps local non archimédien de caractéristique résiduelle $>3$. 
\footnote{Des idées très différentes autour de la correspondance de Howe ont mené~\textcite{Gan_Savin} à une autre proposition pour ce même groupe, pour~$F$ de caractéristique nulle et~$p$ quelconque.} 

\printbibliography

\end{document}
